\documentclass[mnsc,nonblindrev]{informs4} 

\OneAndAHalfSpacedXI 

\usepackage{natbib}
 \bibpunct[, ]{(}{)}{,}{a}{}{,}%
 \def\bibfont{\small}%

\usepackage{enumitem}
\usepackage{amsmath}
\usepackage{mathrsfs}
\usepackage{mathtools}
\usepackage{enumitem}
\usepackage{xcolor}
\usepackage{tikz}
\usepackage{array,graphicx}
\usepackage{float}
\usepackage{pifont}
\usepackage{bm}
\usetikzlibrary{positioning}
\usepackage{url}

\usepackage{breakurl}
\usepackage[breaklinks=true]{hyperref}
\usepackage{pgfplots}
\usepackage{xcolor}
\usepackage{siunitx}
\usepackage{upgreek}
\pgfplotsset{compat=newest}
\pgfplotsset{/pgfplots/error bars/error bar style={very thick}}
\pgfplotsset{
  every axis plot/.append style={very thick, black},
}
\usepackage{accents}
\usepackage{bbm}
\usepackage{comment}
\usepackage{ctable}
\usepackage{clipboard}
\newclipboard{letter}

\usepackage{hyperref}
\hypersetup{
	colorlinks,
	linkcolor={red!50!black},
	citecolor={blue!50!black},
	urlcolor={blue!80!black}
}

\newcommand*{\pr}{\mathsf{P}}
\newcommand{\bigo}{\mathcal{O}}
\newcommand{\indicator}{\mathbb{I}}
\newcommand{\proj}{\textnormal{Proj}}

\newcommand{\numproducts}{\ensuremath{n}}
\newcommand{\product}{\ensuremath{j}}
\newcommand{\productset}{\mathscr{P}}
\newcommand{\numcustomers}{\ensuremath{m}}
\newcommand{\customer}{\ensuremath{i}}
\newcommand{\customerset}{\mathscr{C}}
\newcommand{\distfunc}{\ensuremath{F}}
\newcommand{\pdist}{\ensuremath{G}}
\newcommand{\revfunc}{\ensuremath{R}}
\newcommand{\mrevfunc}{\ensuremath{\Upsilon}}
\newcommand{\surfunc}{\ensuremath{S}}

\newcommand{\histprice}{\ensuremath{P}}
\newcommand{\price}{\ensuremath{p}}
\newcommand{\vecprice}{\mathbf{\price}}

\newcommand{\vechistprice}{\mathbf{\histprice}}
\newcommand{\choice}[1]{\ensuremath{c}_{#1}}

\newcommand{\valuationset}{\mathscr{V}}
\newcommand{\valuation}{\ensuremath{v}}
\newcommand{\vecvaluation}{\mathbf{\valuation}}

\newcommand{\wcfun}{\ensuremath{f}}
\newcommand{\optrevenue}{\tau^*}
\newcommand{\revenue}{\ensuremath{r}}
\newcommand{\eps}{\epsilon}
\newcommand{\histpriceexample}{\ensuremath{P}^*}

\newcommand{\assignVar}{\ensuremath{x}}
\newcommand{\vecassignVar}{\mathbf{\assignVar}}
\newcommand{\valVar}{\valuation^{\product}}

\newcommand{\nopurchase}{\varnothing}
\newcommand{\novalVar}{\valuation^{\nopurchase}}
\newcommand{\vecnovalVar}{\mathbf{\valuation}^{\nopurchase}}

\newcommand{\valoptionset}{\mathscr{W}}

\newcommand{\dualJ}{\mu}
\newcommand{\dualS}{\tau}
\newcommand{\vecdualJ}{\bold{\dualJ}}
\newcommand{\vecdualS}{\bold{\dualS}}

\newcommand{\maxprice}{\ensuremath{\histprice}^{\max}}

\newcommand{\allocVar}{\ensuremath{y}}
\newcommand{\vecallocVar}{\mathbf{\allocVar}}

\newcommand{\approxfun}{\ensuremath{g}}
\newcommand{\precision}{\delta}

\definecolor{mycolor}{rgb}{0.0, 0.0, 0.0}

\newcommand{\rv}{\mathcal{X}}
\newcommand{\cdf}{F_{\rv}}
\newcommand{\maxnum}{R}

\newcommand{\I}{{i}}

\newcommand{\ubar}[1]{\underaccent{\bar}{#1}}

\usepackage{lipsum}

\newcommand{\rev}[1]{\textcolor{black}{#1}}

\TheoremsNumberedThrough     
\ECRepeatTheorems

\EquationsNumberedThrough    

\MANUSCRIPTNO{MS-RMA-21-00031.R1}

\begin{document}


\RUNAUTHOR{Chen et al.} 

\RUNTITLE{Model-free Assortment Pricing}

\TITLE{\Large Model-Free Assortment Pricing with Transaction Data}

\ARTICLEAUTHORS{%
\AUTHOR{Ningyuan Chen}
\AFF{Rotman School of Management, University of Toronto, Toronto, ON, Canada, \EMAIL{ningyuan.chen@utoronto.ca}} 
\AUTHOR{Andre A. Cire}
\AFF{Rotman School of Management, University of Toronto, Toronto, ON, Canada, \EMAIL{andre.cire@rotman.utoronto.ca}} 
\AUTHOR{Ming Hu}
\AFF{Rotman School of Management, University of Toronto, Toronto, ON, Canada, \EMAIL{ming.hu@rotman.utoronto.ca}} 
\AUTHOR{Saman Lagzi}
\AFF{Rotman School of Management, University of Toronto, Toronto, ON, Canada, \EMAIL{saman.lagzi@rotman.utoronto.ca}} 
} 

\ABSTRACT{We study the problem when a firm sets prices for products based on the transaction data, i.e., which product past customers chose from an assortment and what were the historical prices that they observed. Our approach does not impose a model on the distribution of the customers' valuations and only assumes, instead, that purchase choices satisfy incentive-compatible constraints. The individual valuation of each past customer can then be encoded as a polyhedral set, and our approach maximizes the worst-case revenue assuming that new customers' valuations are drawn from the empirical distribution implied by the collection of such polyhedra. We study the single-product case analytically and relate it to the traditional model-based approach. Moreover, we show that the optimal prices in the general case can be approximated at any arbitrary precision by solving a compact mixed-integer linear program. We also design three approximation strategies that are of low computational complexity and interpretable. In particular, the cut-off pricing heuristic has a competent provable performance guarantee. Comprehensive numerical studies based on synthetic and real data suggest that our pricing approach is uniquely beneficial when the historical data has a limited size or is susceptible to model misspecification.
}%



\maketitle

\vspace{-2em}
\section{Introduction}

Online retailing has seen steady growth in the last decade.
According to the survey by \cite{gaubys_2020}, the global market share of e-commerce is expected to surpass 20\% in 2022, and the trend is only accelerating.
Retailers, however, face various challenges when transitioning to an online business model.
First, the granularity of the data gathered from past customers, where firms
observe what customers bought, the products they viewed, and their prices at the time of the purchase, far exceeds that of the offline setting.
Second, the firm is able to roll out products rapidly and many new products (or new configurations of old products) can be displayed each day \citep{caro2020future}. Seemingly contradictory to the first point, there may be little purchase information for a firm to take timely actions on these new products, such as price correction and adjustments.
 

In an attempt to address this challenge, we study the pricing decision of an assortment of products for a firm based on its historical transaction data.
The data records the assortment of products viewed by past customers and their prices, which may vary across customers due to promotions.
\Copy{introcensor}{The data also records the decision made by the customers: the purchase of one of the products. \rev{However, those customers who do not make any purchase may not need to be observed and recorded in the transaction data.}}

The common approach to handle this situation is what we refer to as ``model-estimate-optimize.''
The firm first builds a discrete choice model to characterize how customers form their utilities and make choices.
For example, in the multinomial logit (MNL) model, the probability of a customer choosing product $j$ from a set of $\numproducts$ products priced at $(\price_1,\dots,\price_\numproducts)$ can be expressed as
\begin{equation*}
    \frac{\exp(\alpha_\product-\beta \price_\product)}{1+\sum_{k=1}^n \exp(\alpha_k-\beta \price_k)}
\end{equation*}
for some parameters $\{\alpha_k\}_{k=1}^\numproducts$ and $\beta$ representing the average attractiveness of the products and the price sensitivity, respectively.
It is equivalent to a random utility model in which a random customer has utility $\alpha_k-\beta \price_k+\epsilon_k$ for product $k$, where $\epsilon_k$ follows an independent and identically distributed (i.i.d.) Gumbel distribution and customers choose the product, or no purchase, with the highest utility.
The second step is to ``estimate.''
Based on the historical data, the firm can estimate the parameters (see \citealt{train2009discrete} for the detailed steps).
Finally, the firm ``optimizes'' its pricing decisions by setting the optimal prices for future customers using the estimated model.
This approach has enjoyed wide popularity because of its simplicity and computational efficiency.

However, in the setting we consider, the model-estimate-optimize approach falls short in \rev{three} aspects.
\Copy{two-benefit}{First, when the data size is small relative to the number of products, the estimation of the parameters (such as $\alpha_\product$ in the MNL model) may be noisy and unstable.
It is not clear whether the firm can rely on the estimation to make adequate pricing decisions.
Second, the discrete choice model may not be able to capture the behavioral pattern in the data, which is referred to as model misspecification.
In this case, the firm can choose a more complex model, potentially with more parameters, to minimize the misspecification error.}
This will inevitably exacerbates the first issue, as pointed out in \cite{abdallah2020demand}. \Copy{othercensor}{\rev{Third, the model-based approach would be sensitive to whether the no-purchases are being observed and recorded in the data, in contrast to our robust approach.}}

In this paper, we propose a data-driven approach to the assortment pricing problem.
As opposed to the model-estimate-optimize approach, we consider no prior models on how customers form their utilities for the products.
Instead, when a customer is observed to choose a product from a set of products with given prices, 
we assume that a (model-free) incentive-compatibility condition specifies her potential valuations for the products.
For example, when a product is purchased by a customer, the utility of that product must be necessarily at least as high as any of the products that have not been purchased.
Leveraging this condition, the data for a transaction defines a polyhedral set containing the vector of valuations $(\valuation_1,\dots,\valuation_{\numproducts})$ of the particular customer, which is referred to as the incentive-compatible (IC) polyhedron.
When a new customer arrives, without knowing anything about her preferences,
we uniformly sample from the IC polyhedra from customers in the historical data, and set prices to maximize the revenue for valuations in the sampled polyhedron. More precisely, we set prices such that the revenue derived from the arriving customer is maximized, when her valuations, drawn from the sampled IC polyhedron, lead to the least possible revenue. 

The contribution of this work is fourfold. First, to the best of our knowledge, this is the first study that leverages incentive-compatibility constraints to define a valuation polyhedron for pricing. With this novel approach, we circumvent the need for pre-specified customer models as is customary in the random utility literature, investigating the optimal pricing directly based on the data. We believe this approach sheds new insights into the formulation of other related data-driven problems.

Second, by exploiting the structure of the polyhedra, we present a disjunctive model of the pricing problem and several structural results associated with the optimal prices. We show that the disjunctive model can be approximated at any desired precision by a compact bilinear program, itself solvable by a compact mixed-integer integer programming model after an appropriate reformulation. For scalability purposes, we also present low-complexity, interpretable approximation algorithms that are shown to achieve strong theoretical and numerical performance. We also study special cases of practical interest where the optimal prices can be obtained efficiently.

\color{mycolor}
Third, we build the connection between our model-free approach and the traditional model-based approach when the firm sells only one product.
In particular, the revenue garnered by our approach within a given underlying model
can be expressed in a closed form.
\Copy{halfresult}{As a result, \rev{when the number of historical customers approaches infinity and prices shown to past customers are randomly drawn from a uniform distribution, under mild technical conditions, we are able to establish a tight constant bound of $1/2$ for the performance of our model-free pricing relative to the optimal model-based revenue. This result showcases the robustness of our approach under a model-based framework. Moreover, we show that if the no-purchases are not recorded in the data, our approach can lead to strictly higher revenue than the best the model-based approach could do.}}


\color{black}
Fourth, we conduct a comprehensive numerical study based on synthetic and real data.
In particular, we use the well-studied IRI dataset \citep{Kruger2009} and fit different choice models including the linear, MNL and mixed logit models.
We then generate a small number of customer purchases based on the models.  After applying our approach to the data, the generated revenues significantly outperform the incumbent prices in the data, for all the models, demonstrating the benefit of the data-driven model-free approach:  the insensitivity to model misspecification and stable performance when the data size is limited.

\label{pg:cutoffcomp}\Copy{cutoffcomp}{Our findings in the numerical study suggest that (1) our best approximation algorithm consistently recovers at least $96\%$ of the optimal robust revenue in various settings and \rev{its computational complexity} scales linearly with the number of samples.}
(2) When the historical data is generated from commonly used discrete choice models such as the MNL model,
our data-driven approach performs well, compared to the optimal prices under the correctly specified estimated model, especially for a limited data size. On the other hand, when the discrete choice model is misspecified, our approach is more robust and outperforms the misspecified optimal prices.
(3) Applying our data-driven approach to real datasets leads to increase in revenues over the incumbent prices in the dataset.


\section{Related Work}
\label{sec:relatedwork}
Our study is broadly related to three streams of literature.
The first steam is the papers studying the estimation and optimization of the discrete choice models with prices, which provide the basis for the model-estimate-optimize approach.
For the estimation of popular discrete choice models including the MNL model, \citet{train2009discrete} provides an excellent review.
Recently, two new choice models have drawn attention of scholars in the Operations Research community, the rank-based model \citep{farias2013nonparametric}
and the Markov chain choice model \citep{blanchet2016markov}.
Due to their flexibility, the estimation is not as straightforward as others such as the MNL model.
Several studies propose various algorithms to address the issue \citep{van2015market,van2017expectation,csimcsek2018expectation}.
As for optimal pricing, the MNL model is studied in, e.g., \citet{hopp2005product} and \citet{dong2009dynamic}.
For the nested logit model, which is a generalization of the MNL model, \citet{li2011pricing,gallego2014multiproduct,li2015d} investigate its optimal pricing problem.
\citet{zhang2018multiproduct} study the optimal pricing problem of the generalized extreme value random utility models with the same price sensitivity.
In contrast, \cite{mai2019robust} use a robust framework to tackle the same problem for extreme value utilities.

Although we also adopt a robust approach, we do not rely on a particular distribution and the decision is fully guided by the data.
Similarly, \cite{rusmevichientong2012robust} and \cite{jin2020price} study the robust optimization problem of pricing or assortment planning in the MNL model.
Several papers incorporate pricing into the rank-based model \citep{rusmevichientong2006nonparametric,jagabathula2017nonparametric} and the Markov chain model \citep{dong2019pricing} and study the optimal pricing problem. More recently, \cite{yan2019representative} use the transaction data to fit a general discrete choice model of a representative consumer and solve the optimal pricing based on the fitted model by mixed-integer linear programming. While the problem they study is essentially similar to ours, their approach relies on  the model-estimate-optimize paradigm and hence is conceptually distinct from the proposed approach here.

In the context of revenue management, the definition of the term ``data-driven'' typically depends on the problem context.
When the agent makes decisions in the process of data collection, the data-driven approach is usually associated with a framework that integrates such a process with decision making, so that the agent is learning the unknown parameters or environment while maximizing revenues.\footnote{This is referred to as the online problem, which should not be confused with the notion of online retailing.}
The previous works by \citet{bertsimas2017data,zhang2020data,cohen2018dynamic,cohen2020feature,ettl2020data,ban2020personalized} fall into this category.
It is connected to a large body of literature on demand learning and dynamic pricing. We refer to \citet{den2015dynamic} for a comprehensive review of earlier papers in this setting.

In contrast, our study essentially handles an offline setting, in which the data has been collected and given. Two recent papers also provide alternatives to the model-estimate-optimization approach.
\citet{bertsimas2020predictive} leverage statistical methods such as $k$-nearest neighbors and kernel smoothing to integrate past observations into the current decision making problem given a covariate.
It circumvents the step of estimating a statistical model.
\citet{elmachtoub2017smart} achieve a similar goal by skipping the minimization of the estimation error and directly focusing on the decision error.
They design a new loss function that combines the errors in both stages: estimation and optimization.
Both papers study a setting where the optimization is conditioned on a covariate.
Our problem is less general than their formulation and does not have a covariate, which allows us to utilize the special structure of the multi-product pricing problem (i.e., the incentive-compatibility of customers) that does not apply to their general approach.
\cite{ban2019big} propose an algorithm integrating historical demand data and newsvendor optimal order quantity, without estimating the demand model separately.

There are a few papers using model-free approaches in revenue management.
\citet{allouah2019sample} study the problem that the seller observes samples from the buyer's valuation but is agnostic to the underlying distribution.
The optimal price is solved for a family of possible distributions in the maximin sense.
\citet{chen2019use} adopt the estimate-then-model approach, which is model-free, to estimate choice models using random forests.
They do not consider the optimization stage.
Our problem is similar to the work by \citet{ferreira2016analytics}, who study the demand forecasting and price optimization of a new product by a retailer.
However, it focuses on a single product and it is unclear how to adapt it to multiple products.
A few studies apply analytics to promotion planning \citep{cohen2017impact,cohen2020optimizing}.
They consider more practice-based assumptions than ours and use various approximations for the demand model.

This paper uses the incentive-compatibility of past customers as a building block, and hence is related to the literature on auction designs, especially those papers using a data-driven or robust formulation.
\cite{bandi2014optimal} study the multi-item auction design with budget constrained buyers and use a robust formulation for the set of valuations.
The uncertainty set is constructed using the historical information, such as means and covariance matrices.
Because of the polyhedral structure of the uncertainty sets, the robust optimization problem is tractable.
Our formulation investigates the worst-case revenue in the IC polyhedron (see Section~\ref{sec:problem}), which is similar to the idea of robust optimization for valuations drawn from uncertainty sets. However, there are three distinctions.
First, in \cite{bandi2014optimal}, the uncertainty set is constructed on the valuations of all bidders for a single product, while in our problem
the valuations of a single customer for all the products fall into a polyhedron.
This is because of the different applications.
Second, our data-driven approach averages over the empirical distribution of historical customers, while the historical information is used to construct the uncertainty set in \cite{bandi2014optimal}.
Third, the optimization problem in our problem cannot be solved efficiently and we resort to approximation algorithms.
\label{pg:derakhshani}\Copy{derakhshani}{\cite{derakhshan2021linear} consider a data-driven optimization framework to find the optimal personalized reservation price of buyers, when past bids are input to the algorithm.
They do not impose assumptions on the valuation distributions and maximize the revenue when the valuation of the future customer is drawn from the empirical distribution of the historical data. 
This is similar to the motivation of this paper.
However, because of the different contexts, the reduction and approximations have little in common with ours.}
Similarly, \cite{koccyiugit2020distributionally} study the problem of designing an auction to sell an indivisible good to a group of bidders when the bidders valuations come from an ambiguous distribution and the bidders attitude towards this ambiguity is unknown. They design optimal mechanisms that are robust to the worst-case realization of bidders valuation and their attitudes. 
\cite{koccyiugit2018robust} apply the robust approach to multiproduct pricing. The authors assume the valuation of a single customer is drawn from a rectangular uncertainty set and design mechanisms to maximize the worst-case revenue ensuring incentive compatibility.
They show that the optimal robust selling mechanism is to sell products separately with randomized posted prices.
Finally, \cite{allouah2020prior} study the single-item auction for a group of buyers
when the seller does not have access to the valuation distribution of the buyers and buyers do not have any information about their competitors. 
Instead, the auction is designed for a general class of distributions with a competitive ratio objective, which is conceptually similar to our model-free approach.

\section{Problem Description}
\label{sec:problem}

We consider a firm that has observed the transaction data of $\numcustomers$ customers $\customerset = \{1,\dots,\numcustomers\}$ with respect to an assortment of $\numproducts$ products $\productset = \{1,\dots,\numproducts\}$. Specifically, the firm's historical data includes the product prices $\vechistprice_\customer = (\histprice_{\customer 1}, \dots, \histprice_{\customer \numproducts} ) > \mathbf{0}$ viewed by each customer $\customer \in \customerset$, as well as the product $\choice{\customer} \in \productset$ that the customer chose from that assortment. 
\Copy{censor}{We consider that customers viewed all products $\productset$ and purchased one product from that assortment \rev{(customers without purchases are censored in the data)}; both assumptions are made without loss of generality (Remarks~\ref{rem:nopurchasecustomers} and \ref{rem:observabilityOfProducts} below). That is, we can incorporate the situation where historical customers may have not observed some of the products in the assortment during their transactions, and removing those customers who made no purchase during their shopping session does not change our results. } The goal of the firm is to set the product prices $\vecprice = (\price_1, \dots, \price_\numproducts)$ for newly arriving customers that leverages this offline historical data.

In this study, we investigate a pricing approach that operates on the full set of possible customer utilities under incentive-compatibility constraints.
More precisely, let
$\valuation_{\customer \product} \ge 0$ be the \emph{unknown} valuation that customer $\customer \in \customerset$ assigns to each product $\product \in \productset$. We assume that the utility of purchasing $\product$ is given by $\valuation_{\customer \product} - \histprice_{\customer\product }$ (and hence quasilinear in $\valuation_{\customer \product}$), and that such utilities must be compatible with observations from the data; i.e., the set of all possible valuations of customer $\customer$ is
\begin{align}
    \label{eq:i-valuation}
    \valuationset_\customer \equiv
    \big\{
        (\valuation_{\customer 1}, \dots, \valuation_{\customer \numproducts}) \in \mathbb{R}^{\numproducts}_{+}
        \; \colon \;
            &\valuation_{\customer \choice{\customer}} - \histprice_{\customer \choice{\customer}} \ge 0, 
            \, \valuation_{\customer \choice{\customer}} - \histprice_{\customer \choice{\customer}} \ge \valuation_{\customer \product'} - \histprice_{\customer \product'}, \;\; \forall \product' \in \productset \setminus \{ \choice{\customer} \}
    \big\}.
\end{align}
The first inequality in \eqref{eq:i-valuation} indicates that the utility of purchasing product $\choice{\customer}$ is non-negative. The second inequality specifies that valuations are incentive-compatible, i.e., that the utility of purchasing product $\choice{\customer}$ must be either the same or larger than the utility for the remaining products $\productset \setminus \{ \choice{\customer} \}$. We note that $\valuationset_\customer$ is defined by a finite set of closed halfspaces and hence is a polyhedral set.
We refer to $\valuationset_\customer$ as the \textit{incentive-compatible} (\textit{IC}) \textit{polyhedron} of customer $\customer$.

For a newly arriving customer, however, her valuation is not known to the firm. Based on historical data, it is reasonable to use the empirical distribution to form such an estimate. That is, we assume that her valuation is equally likely to fall into one of the IC polyhedra $\valuationset_1, \valuationset_2, \dots, \valuationset_\numcustomers$.
Given the new prices $\vecprice$ and considering that the new customer's valuation falls into $\valuationset_\customer$ for some $\customer$,
we apply a robust approach where the arriving customer picks product $\product \in \productset$ that yields the lowest possible revenue and such a choice is consistent with $\valuationset_\customer$ under prices $\vecprice$.
Such a revenue is described by the program
\begin{align}
    \wcfun_\customer(\vecprice)
    \equiv
    &&\min_{\vecvaluation_\customer \in \valuationset_\customer, \revenue \ge 0}
        &\quad
        \revenue
            \label{model:DP} \tag{DP} \\
    &&\textnormal{s.t.}
        &\quad
            \bigvee_{\product \in \productset}
            \left(
                \begin{array}{l}
                    \revenue = \price_\product \\
                    \valuation_{\customer \product} - \price_\product \ge 0 \\
                    \valuation_{\customer \product} - \price_\product \ge \valuation_{\customer \product'} - \price_{\product'}, \;\; \forall \product' \in \productset
                \end{array}
            \right)
            \vee
            \left(
                \begin{array}{l}
                    \revenue = 0 \\
                    \valuation_{\customer \product} - \price_{\product} \le 0, \;\; \forall \product \in \productset
                \end{array}
            \right).
            \label{eq:disjConditions}
\end{align}
Model \eqref{model:DP} (\rev{short for Disjunctive Program}) is a classical disjunctive program \citep{balas1998disjunctive}, i.e., the set of feasible solutions is defined by a disjunction of polyhedra representing the feasible valuations for each product.
Specifically, with the objective $\min r$ we select the product with the lowest revenue for which there exists a feasible valuation
$\vecvaluation_\customer \in \valuationset_\customer$ that is also incentive-compatible under the new prices $\vecprice$.
This is modeled by the left-hand side disjunction of \eqref{eq:disjConditions}, where some product $\product$ (and hence price $\price_\product$) is selected only if its net utility (i.e., $\valuation_{\customer \product} - \price_\product$) is non-negative and at least as large as that of choosing other products.
The right-hand side disjunction of \eqref{eq:disjConditions} formulates the case where the customer chooses no product, only possible if the utility of choosing any product under $\vecprice$ is non-positive.

\smallskip
\begin{remark}
    \label{rem:nopurchasecustomers}
    As noted earlier, we assume that \rev{the data is censored, and} each historical customer $\customer \in \customerset$ has purchased one product from the assortment. In particular, if some customer $\customer$ has not purchased any products, we can assume that $\customer$ belongs to the polyhedron with zero revenue in model \eqref{model:DP}, i.e., the right-hand side term in disjunction \eqref{eq:disjConditions}. Thus, revenue driven by customer $\customer$ is always zero for any prices, and therefore she can be removed from the model. 
    \hfill $\square$
\end{remark}

\smallskip
Based on \eqref{model:DP}, we assign a weight $1/\numcustomers$ and sum over $\customer\in \customerset$ to take the expectation with respect to the empirical distribution.
To maximize over $\vecprice$, the objective function is expressed as
\begin{align*}
\optrevenue \equiv \sup_{\vecprice \ge 0} \;\;& \frac{1}{\numcustomers} \sum_{\customer \in \customerset} \wcfun_{\customer}(\vecprice).
\label{model:OP} \tag{OP}
\end{align*}

The problem \eqref{model:OP} (\rev{short for Optimal Pricing}) highlights important distinctions from existing pricing approaches.
First, it does not rely on a parametric discrete-choice model that explicitly specifies the distribution of customer utilities (see, e.g., \citealt{train2009discrete}).
Instead, we adopt a model-free approach and consider the worst-case valuation for a given $\vecprice$ under quasilinearity and incentive-compatible customer preferences, which are arguably weaker and more justifiable than existing parametric models to the best of our knowledge.
Second, we do not attempt to estimate a non-parametric model (e.g., \citealt{jagabathula2017nonparametric}).
Instead, the goal is to investigate the structure of \eqref{model:OP}, trade-offs, and scenarios where prices derived from \eqref{model:OP} are beneficial in comparison to existing models, given its data-driven nature and emphasis on the historical data.

\smallskip
\begin{remark}
    \label{rem:observabilityOfProducts}
    The choice model \eqref{model:DP} assumes that each customer $\customer \in \customerset$ has seen the same assortment $\productset$
and their historical prices $\vechistprice_\customer$.
While in practice customers may have seen different assortments, which can be subsets of $\productset$, this assumption can be relaxed by setting a sufficiently large $\histprice_{\customer \product}$ if customer $\customer$ has not been offered product $\product$ (e.g., the sum of all historical prices). In that case, the structural results below can be rewritten accordingly without loss of generality. $\square$
\end{remark}

\rev{We next discuss two characteristics of the pricing model.}

\noindent \textit{Attaining the maximum of \eqref{model:OP}.} We note that the maximum of \eqref{model:OP} may not be attainable because of the discontinuity of the customer choices in \eqref{eq:disjConditions}, hence the use of supremum in the objective. For example, consider an instance with one customer ($\numcustomers=1$) and one product ($\numproducts=1$), where $\histprice_{1 1} = \histpriceexample$ for some $\histpriceexample > 0$ and $\choice{1} = 1$. Thus, from \eqref{eq:i-valuation}, the customer valuation satisfies $\valuation_{1 1} \ge \histpriceexample$. Suppose now we assign a price $\price_1 \ge \histpriceexample$. It follows from \eqref{eq:disjConditions} that the optimal revenue is zero, since any valuation $\histpriceexample \le \valuation_{1 1} \le \price_1$ is incentive-compatible with the no-purchase option. However, for any $\price_1 = \histpriceexample - \eps$ with $\eps > 0$, we have $\optrevenue \rightarrow \histpriceexample$ as $\eps \rightarrow 0.$ Thus, $\wcfun_\customer(\vecprice)$ is discontinuous in $\vecprice$.

\Copy{interpretation}{
\noindent \textit{Interpretation as a choice model.}
One may take the worst-case choices embedded in \eqref{model:DP} of all historical customers and 
assign an equal probability to them.
This construction could potentially be interpreted as the choice behavior of the new customer. 
However, the essence of our approach is to take a conservative mapping from the price vector directly to the revenue, bypassing the choice model. We do not claim that \eqref{model:OP} can be used as a proper predictive model for consumer choices.}

\section{Model-Free Pricing for a Single Product}
\label{sec:single-prod}
In this section, we investigate the special case of a single product, i.e., $\numproducts=1$.
The traditional model-based approach for optimal pricing relies on a demand function that characterizes the distribution of customers' valuations. 
Given that such a notion is absent in our framework, our goal in this section is to explore the connection between the proposed pricing scheme and the model-based approach when a known demand function specifies consumer purchase behavior.
That is, we attempt to rigorously characterize the performance of model-free pricing in a model-based setting. \label{pg:roadmap}\Copy{roadmap}{\rev{To that avail, in this section we will provide a sufficient condition under which the model-free pricing achieves the optimal revenue in the model-based setting. Moreover, we provide tight performance guarantees for the model-free price, in the model-based setting, when historical prices are drawn independently from a uniform distribution. Finally, we characterize the sample complexity of our model-free approach. Namely, we show the rate, in terms of number of samples, at which the model-free optimal price converges to its asymptotic optimal price.}}

Let $\distfunc: [0,+\infty)\to [0,1]$ denote the cumulative distribution function (CDF) of the random customer valuation for the product.
Thus, given a price $p\ge 0$, the probability that a customer purchases the product is $1-\distfunc(p)$. 
The expected revenue at $p$, which is the target of model-based approaches, is $\revfunc(p)\triangleq p(1-\distfunc(p))$.
Our objective is to compare the price from \eqref{model:OP} with $\argmax_{p\ge 0} \revfunc(p)$. Note that the latter requires an infinite number of data points to estimate $\distfunc$ exactly \rev{and is susceptible to data censoring, which we discuss in Lemma~\ref{lemma:censlemma} in Section \ref{sec:censresults} of the e-companion}.

To generate the historical data that can be used by our framework, 
we have to specify how the historical prices are drawn. 
We assume that the seller picks a price \emph{independently} from a distribution with CDF $\pdist(\cdot)$.\footnote{\textcolor{mycolor}{We argue that this is a reasonable assumption. To estimate $\distfunc(\cdot)$ accurately in a model-based approach, the historical prices need to span the whole price range, presumably by conducting an experiment with random prices.}}
Observing a price $p\sim \pdist$, the customer makes a purchase with probability $1-\distfunc(p)$.
Because of the censoring,
the historical data in our framework $\{(\histprice_{i1}, \choice{i})\}_{i=1}^\numcustomers$ can be 
thought of being generated independently from the distribution  $\pr(\histprice_{i1}\ge p) =\int_p^{+\infty} (1- \distfunc(x)) \mathrm{d}\pdist(x) $ and $\choice{i}=1$. 
\rev{As will be shown in Section~\ref{reformulation},  
the optimal price for our model-free approach satisfies
$\price^*_m\in \argmax_{\price\ge 0} \price \sum_{i=1}^\numcustomers \indicator(\histprice_{i1}\ge \price)$.}

To compare the fundamentals of the model-based approach and our framework without the statistical noise, 
we let $\numcustomers\to\infty$ and define $p^*$ as $\price^*\in \argmax_{\price\ge 0} \price \left(\lim_{\numcustomers\to\infty}\frac{\sum_{i=1}^\numcustomers \indicator(\histprice_{i1}\ge \price)}{\numcustomers}\right)$.
Moreover, for the simplicity of analysis and exposition, we assume
\begin{assumption}\label{asp:cdf-F-G}
Both CDFs $\distfunc(\cdot)$ and $\pdist(\cdot)$ have a probability density function (PDF), denoted by $f(\cdot)$ and $g(\cdot)$, respectively.
\end{assumption}
Our next result characterizes the optimal data-driven price $\price^*$ in our framework asymptotically and compares it to the optimal price out of the model-based approach assuming the model is known.
\begin{proposition}\label{prop:compare-opt-price}
\Copy{prop-compare-price}{Under Assumption~\ref{asp:cdf-F-G}, we have
\begin{equation}\label{eq:opt-price-single-prod}
    p^* \in \argmax_{\price\ge 0} \price\int_\price^{+\infty} g(x)(1-\distfunc(x))\mathrm{d}x.
\end{equation}
Moreover, when $g(\price) \propto f(\price)/(1-\distfunc(\price))$, the optimal price of our framework is also optimal for the model-based approach:
\begin{equation*}
    p^* \in \argmax_{\price\ge 0} \revfunc(\price)=\argmax_{\price\ge 0} p(1-F(p)).
\end{equation*}}
\end{proposition}

Proposition~\ref{prop:compare-opt-price} reveals connections and differences between the two approaches.
First, the optimal price of our framework depends on the distribution of the historical prices, which itself does not affect the optimal price of the model-based approach.
The reason that the historical prices impact the expected revenue in our data-driven setting is that our approach needs to recover the true valuation distribution from censored demand. To estimate  $\distfunc(\cdot)$ well, one needs to know the fraction of no-purchases for the historical prices since the observed demand is censored by the offered price. 

\Copy{limitation}{Second, even when the data size is sufficiently large, the model-free pricing is not expected to converge to the optimal model-based prices in general.} When the historical prices in our framework have a probability density function specified by the \emph{hazard rate} function of the valuation distribution, the optimal prices of the model-based and model-free approaches coincide.  
However, one would not know the hazard rate function exactly, otherwise it would be possible to infer  $\distfunc(\cdot)$ from such a hazard rate and use the optimal price of the model-based approach, assuming that the model is correctly specified.

\subsection{Model-free Pricing under Uniformly Distributed Historical Prices}

Next, we investigate a simple price distribution that does not require the information of $\distfunc$---the uniform distribution in the historical data---and study the performance of our framework. We consider the following assumption.
\begin{assumption}\label{asp:uniform}
The distribution $f(\cdot)$ is supported on $[0,a]$ and $g(\price) = 1/b$ for $\price\in [0,b]$, where $0\le a \le b$.
\end{assumption}

Under Assumption~\ref{asp:uniform}, the first-order condition for  $\price^*$ in \eqref{eq:opt-price-single-prod} implies that
\begin{equation}\label{eq:focsing}
    \revfunc(p^*)=p^*(1-\distfunc(p^*)) = \int_{p^*}^a (1-\distfunc(x))\mathrm{d}x \triangleq \surfunc(p^*).
\end{equation}
\Copy{divider}{In other words, the optimal model-free price from our framework is an equal divider of the welfare: the revenue earned by the firm, $\revfunc(p^*)$, is equal to the surplus gained by the consumers, $\surfunc(p^*)$.} Our question therefore can be framed as follows: if we plug the optimal price from our framework $p^*$ (or equivalently, an equal divider of the welfare) into $R(p^*)$, how large would the gap be between the optimal model-free profit $R(p^*)$ and the optimal model-based profit $\max_p R(p)$?
Once equipped with Assumption \ref{asp:unique-max}, which we will show later is not restrictive, Theorem~\ref{prop:rev-bound} answers this question definitively.
\begin{assumption}\label{asp:unique-max}
The maximizer of \eqref{eq:opt-price-single-prod}, $p^*$, is the only solution to \eqref{eq:focsing} in $[0,a]$ such that $p^*\int_{p^*}^{a}(1-\distfunc(x))\mathrm{d}x>0$;
$\revfunc(p)$ is unimodal in $[0,a]$ and has a unique maximizer $\hat p$.
\end{assumption}
\begin{theorem}
\label{prop:rev-bound}
Suppose Assumptions~\ref{asp:cdf-F-G},~\ref{asp:uniform} and~\ref{asp:unique-max} hold.
We have
$
    \revfunc(p^*)/\revfunc(\hat p)\ge \min\left\{\hat p/p^*, 1/2\right\}.
$
Moreover, the bound is asymptotically tight: for any $\epsilon>0$, we can construct $\distfunc(\cdot)$ such that $\revfunc(p^*)/\revfunc(\hat p)\leq 1/2+\epsilon$; for any $0<x_1<x_2<1/e$, we can construct $\distfunc(\cdot)$ such that $\hat p=x_1$, and $|p^*-x_2|<\epsilon$ such that $\revfunc(p^*)/\revfunc(\hat p)\leq \hat p^{1-\epsilon}/p^{*1-\epsilon}$.
\end{theorem}

\label{pg:prosketch1}\Copy{prosketch1}{
\begin{proof}{Proof Sketch.}
Based on $\distfunc(\cdot)$, we may have $p^*\le \hat p$ or $p^*> \hat p$.
Hence we consider both cases separately.
When $p^*\le\hat p$, we show that $\revfunc(p^*)/\revfunc(\hat p)\ge1/2$. Then, for any $\epsilon>0$, we construct a distribution $\distfunc(\cdot)$ that leads to $\revfunc(p^*)/\revfunc(\hat p)\leq1/2+\epsilon$. Intuitively, the worst-case performance of $p^*$ is achieved when customers' valuation distribution is heavily concentrated at $a$.

When $p^*>\hat p$, we leverage the fact that $p^*$ is the  divider of welfare into equal portions of revenue and consumer surplus.
At $\hat p$, consumer surplus is strictly larger than the firm's revenue.
These observations lead to $\revfunc(p^*)/\revfunc(\hat p)\ge \hat p/p^*$. 
Moreover, for any $\epsilon>0$, we construct a distribution $\distfunc(\cdot)$ that achieves $\revfunc(p^*)/\revfunc(\hat p)\leq \hat p^{1-\epsilon}/p^{*1-\epsilon}$. To elaborate on the intuition behind this worst-case performance, beyond $\hat p$, both $S(p)$ and $R(p)$ deteriorate at a nearly identical rate, leading to $p^*$ to be as far away from $\hat p$ as possible. \hfill $\square$
\end{proof}}

We note that our bound is asymptotically tight, as $\epsilon$ in Theorem~\ref{prop:rev-bound} can be chosen to be arbitrarily close to zero. 
Theorem~\ref{prop:rev-bound} establishes a performance bound when we use model-free pricing in the model-based setting with a demand function.
Since the revenue maximizer $\hat p$ is typically larger than the welfare divider $p^*$, using our optimal price guarantees at least a half of the optimal revenue in this case.
Note that our framework is completely model-free and does not assume a behavioral pattern of consumers. Theorem~\ref{prop:rev-bound} suggests the robust performance of our approach even when there is an embedded model. Moreover, it also establishes the revenue bound for the welfare divider $p^*$, which may be of independent interest. The proof for the bound and the construction of tight instances are highly non-trivial.
We defer it to the e-companion. 

We would like to point out that the bound established in Theorem~\ref{prop:rev-bound} is applicable to a very wide range of customer valuation distributions and can be improved significantly when customers' valuation distribution follows a specific distribution such as the uniform distribution (see Proposition~\ref{prop:unif-bound} in the e-companion for more details). The next lemma illustrates the applicability of the bound in Theorem~\ref{prop:rev-bound} by establishing that Assumption~\ref{asp:unique-max} is not restrictive and is satisfied for a large array of distributions.

\begin{lemma}\label{lem:IGFR}
Suppose the distribution $f(\cdot)$ is supported on $[0,a]$ and has a strictly increasing generalized failure rate, i.e., $xf(x)/(1-F(x))$ is strictly increasing on $[0,a]$. Then Assumption~\ref{asp:unique-max} holds.
\end{lemma}

Lemma~\ref{lem:IGFR} shows that Assumption~\ref{asp:unique-max} holds for distributions with a finite support that have a strictly increasing generalized failure rate. As \cite{ziya2004relationships} point out, assuming a strictly increasing generalized failure rate for customer valuation distributions is widely used in the literature of revenue management. 

Another implication is that when $\hat p\ge p^*$, Theorem~\ref{prop:rev-bound} guarantees that $\revfunc(p^*)/\revfunc(\hat p)\ge 1/2$. The following lemma puts forward a sufficient condition for $\hat p\ge p^*$. 

\begin{lemma}\label{lem:fixedbound}
Suppose Assumptions~\ref{asp:cdf-F-G},~\ref{asp:uniform} and~\ref{asp:unique-max} hold. Moreover, let $H(p)={\int_{0}^{p}(1-\distfunc(x))\mathrm{d}x}/{\int_{0}^{a}(1-\distfunc(x))\mathrm{d}x}$ and assume that $H(\cdot)$ has a failure rate that is no less than that of $\distfunc(\cdot)$, i.e., ${h(p)}/{1-H(p)}\ge {f(p)}/{1-\distfunc(p)}$ for any $p\in [0,a]$, where $h(p)={1-\distfunc(p)}/{\int_{0}^{a}(1-\distfunc(x))\mathrm{d}x}$. Then, we have $\hat p\ge p^*$ and thus ${\revfunc(p^*)}/{\revfunc(\hat p)}\ge {1}/{2}$ by Theorem~\ref{prop:rev-bound}.
\end{lemma}

As pointed out by \cite{chen2020digital}, the conditions of Lemma~\ref{lem:fixedbound} are arguably mild and are satisfied for many distributions that have increasing failure rates and are commonly used in the marketing and operations management literature, e.g., uniform and right-truncated exponential distributions.

\label{pg:singlepcens}\Copy{singlepcens}{\rev{Finally, if the no-purchases are not recorded in the data, our approach can lead to strictly higher revenue than the best the model-based approach could do (see Section \ref{sec:censresults} of the e-companion). This result demonstrates the advantage of our model-free robust approach in requiring no information on observing the no-purchases, in contrast to the traditional model-based approach which is sensitive to the information on the no-purchases (see Figure~\ref{fig:censor} in the e-companion).}}

\subsection{Sample Complexity of Model-Free Pricing}\label{samplecomplexity}

\rev{
We investigate the sample complexity of model-free pricing, i.e.,
how fast the optimal price of our framework approaches its asymptotic optimal price $p^*$, given the increasing number of historical customers, $\numcustomers$. The following proposition establishes the rate of convergence under some mild conditions.}

\begin{proposition}\label{prop:convrate}
Suppose Assumption~\ref{asp:cdf-F-G} holds and $g(\cdot)$ is continuous and supported on $[0,b]$.
Assume that $p^*$ is the unique maximizer of $\max_{\price\ge 0}  \price\int_\price^{b} g(x)(1-\distfunc(x))\mathrm{d}x$ and there exists $\alpha>0$ such that $p^* \int_{p^*}^{b} g(x)(1-\distfunc(x))\mathrm{d}x-p \int_\price^{b} g(x)(1-\distfunc(x))\mathrm{d}x\ge \alpha (p-p^*)^2$ for all $p\in [0,b]$.
Then with $m$ historical customers and for any $\epsilon>0$ we have $\pr(|p^*_m-p^*|\ge \epsilon)\le 4e^{\frac{-\alpha^2\epsilon^4 m}{2b^2}}$. This suggests a sample complexity of $O\big((\frac{b^2}{\alpha^2\epsilon^4})\log(\frac{1}{\delta})\big)$ to ensure that  $|p_m^*-p^*|\ge \epsilon$ with probability at most $\delta$. 

\end{proposition}

\color{black}

\section{Reformulations of the Pricing Model \eqref{model:OP}}\label{reformulation}

In this section, we investigate linear reformulations of \eqref{model:OP} that serve as the basis of our exact and tractable approximation strategies.
We start in \S\ref{sec:revenuePolyhedron} with an analysis of the optimal policy structure for the problem \eqref{model:DP}, presenting a small linear programming model to compute $\wcfun_\customer(\vecprice)$ that is compact with respect to the number of products and customers. We leverage this model in \S \ref{sec:milpOP} to derive an alternative, compact mixed-integer linear program that approximates $\optrevenue$ at any desired absolute error with respect to the optimal solution of \eqref{model:OP}. Finally, we use this reformulation in \S\ref{sec:special_cases} to show the optimal price structure of two commonly seen special cases.

\subsection{Optimal Revenue from an IC Polyhedron}\label{sec:revenuePolyhedron}
We first characterize the IC polyhedron $\valuationset_\customer$ and $\wcfun_\customer(\vecprice)$, for given $\customer \in \customerset$ and $\vecprice$.
It can be interpreted as how customer $\customer$ behaves under the new price $\vecprice$ based on her historical choice.
Proposition \ref{prop:revenueLP} uses classical polyhedral results to transform \eqref{model:DP} into an equivalent linear program.


\begin{proposition}
\label{prop:revenueLP}
The formulation \eqref{model:DP} is equivalent to the linear program
\begin{align}
    \wcfun_\customer(\vecprice)
    =
    \min_{\vecvaluation^1_\customer, \dots, \vecvaluation^{\numproducts}_\customer, \vecnovalVar_\customer, \vecassignVar \ge 0}
        &\quad
        \sum_{\product \in \productset} \price_\product \assignVar_\product
            \label{model:DPLP} \tag{DP-LP} \\
    \textnormal{s.t.}
        &\quad
            \valVar_{\customer \product} - \price_\product \assignVar_\product \ge 0,
            &\forall \product \in \productset,
            \label{eq:purchaseH-1} \\
        &\quad
            \valVar_{\customer \product} - \price_\product \assignVar_\product
            \ge
            \valVar_{\customer \product'} - \price_{\product'} \assignVar_\product,
            &\forall \product, \product' \in \productset,
            \label{eq:purchaseH-2} \\
        &\quad
            \novalVar_{\customer \product} \le \price_\product \assignVar_{\nopurchase},
            &\forall \product \in \productset,
            \label{eq:nopurchase} \\
        &\quad
            \sum_{\product \in \productset} \assignVar_{\product} + \assignVar_{\nopurchase} = 1,
            \label{eq:onechoiceCustomer} \\
        &\quad
            \valuation^{\product}_{\customer \choice{\customer}} - \histprice_{\customer \choice{\customer}} \assignVar_{\product} \ge 0,
            &\forall \product \in \productset \cup \{ \nopurchase \}, \label{eq:feasibleValuation-1} \\
        &\quad
            \valuation^{\product}_{\customer \choice{\customer}} - \histprice_{\customer \choice{\customer}} \assignVar_{\product}
            \ge
            \valVar_{\customer \product'} - \histprice_{\customer \product'} \assignVar_{\product},
            &\forall \product, \product' \in \productset \cup \{ \nopurchase \},
            \label{eq:feasibleValuation-2}
\end{align}
where $\vecvaluation^\product_\customer = (\vecvaluation^\product_{\customer 1} , \dots, \vecvaluation^{\numproducts}_{\customer\numproducts})\ \forall\ \product \in \productset \cup \{ \nopurchase \}$ and $\vecassignVar = (\assignVar_1, \dots, \assignVar_{\numproducts}, \assignVar_{\nopurchase})$.
\end{proposition}

The model \eqref{model:DPLP} (\rev{short for Disjunctive Program LP reformulation}) is a linear program with $\bigo(\numproducts^2)$ variables and
constraints. In particular, the variables $\assignVar_{1}, \dots, \assignVar_{\numproducts}, \assignVar_{\nopurchase}$ represent which product is picked by an arriving customer who draws the IC polyhedron $\valuationset_\customer$, where $\assignVar_{\nopurchase}$ encodes the no-purchase option.
Constraint \eqref{eq:onechoiceCustomer} ensures that the customer selects either one product or the no-purchase option.
Constraints \eqref{eq:purchaseH-1}-\eqref{eq:purchaseH-2} imply that $\vecvaluation^\product_\customer$ is incentive-compatible with the choice $\assignVar_\product$.
They are equivalent to the $\product$-th disjunctive term of \eqref{eq:disjConditions}, while the other inequalities $\product' \neq \product$ of the same family become redundant whenever $\assignVar_{\product} = 0$.
This same reason applies analogously to the no-purchase option and inequality \eqref{eq:nopurchase}.
Constraints \eqref{eq:feasibleValuation-1}-\eqref{eq:feasibleValuation-2} ensure that $\vecvaluation^\product_\customer\in \valuationset_\customer$, i.e., the valuation must be incentive-compatible with the historically chosen product $\choice{\customer}$.
The generated revenue $\sum_{\product\in \productset} \price_\product \assignVar_\product$ is minimized over $\vecvaluation$ and $\vecassignVar$.


We now analyze the structure of \eqref{model:DPLP} to draw insights into the optimal customer choices and reduce the size of the formulation. To this end, consider the set of valuations from \eqref{model:DPLP} that are incentive-compatible with product $\product\in\productset$, i.e.,
$
    \valoptionset^{\product}_{\customer}(\vecprice)
    \equiv
    \big\{
        \vecvaluation^\product_\customer \in \valuationset_\customer
        \colon
        \valVar_{\customer \product} - \price_\product \ge 0,
        \;\;
        \valVar_{\customer \product} - \price_\product
        \ge
        \valVar_{\customer \product'} - \price_{\product'},
        \;\; \forall \product' \in \productset
    \big \}
$
and those incentive-compatible with the no-purchase option:
$
    \valoptionset^{\nopurchase}_{\customer}(\vecprice)
    \equiv
    \left \{
        \vecvaluation^{\nopurchase}_\customer \in \valuationset_\customer
        \colon
        \novalVar_{\customer \product} \le \price_\product,
        \;\; \forall \product \in \productset
    \right \}.
$
We now characterize in Lemma \ref{lem:feasibleOptions} when such valuation sets have at least one feasible point,
that is, there exists a valuation vector for the products that is incentive-compatible with both the historical choice $\choice\customer$ under the historical price $\vechistprice_\customer$ and product $\product$ being chosen under the new price $\vecprice$.
\begin{lemma}
    \label{lem:feasibleOptions}
    For any price $\vecprice \ge 0$, the following statements \ref{it:feasibleNopurchase}-\ref{it:feasibleH} hold:
    \begin{enumerate}[label=(\alph*)]
        \item The no-purchase option is feasible to the $\customer$-th customer class  ($\valoptionset^{\nopurchase}_{\customer}(\vecprice) \neq \emptyset$) if and only if $\price_{\choice{\customer}} \ge \histprice_{\customer \choice{\customer}}$, i.e., the new price of the historically chosen product $\choice{\customer}$ remains the same or increases.
        \label{it:feasibleNopurchase}

        \smallskip
        \item The purchase of the historically chosen product $\choice{\customer}$ by the $\customer$-th customer class is always feasible ($\valoptionset^{\choice{\customer}}_{\customer}(\vecprice) \neq \emptyset$) for all $\vecprice \ge 0$. \label{it:feasibleChoice}

        \smallskip
        \item The purchase of $\product \in \productset \setminus \{\choice{\customer}\}$
        by the $\customer$-th customer class is feasible ($\valoptionset^{\product}_{\customer}(\vecprice) \neq \emptyset$) if and only if $\price_{\product} - \price_{\choice{\customer}} \le \histprice_{\customer \product} - \histprice_{\customer \choice{\customer}}$, i.e.,
        the price difference of $\product$ w.r.t. $\choice{\customer}$ remains the same or decreases.
        \label{it:feasibleH}
    \end{enumerate}
\end{lemma}

\smallskip

Lemma~\ref{lem:feasibleOptions} provides easy-to-check conditions for whether the choice of a specific product or none is feasible.
It also leads to a more compact formulation of \eqref{model:DPLP}.
Intuitively, one can simply screen all products $j\in \productset$ and the no-purchase option according to Lemma~\ref{lem:feasibleOptions} for feasible options, and choose the one with the lowest revenue.
Formally, let $\indicator(C)$ be the indicator function of the logical condition $C$, i.e., it is equal to 1 if $C$ is true and 0 otherwise. Proposition \ref{prop:compactDPLP} applies Lemma \ref{lem:feasibleOptions} to a reformulation of \eqref{model:DPLP} via a projective argument.

\begin{proposition}
    \label{prop:compactDPLP}
    The formulation \eqref{model:DPLP} is equivalent to the linear program
    \begin{align}
        \wcfun_\customer(\vecprice)
        =
        \min_{\vecassignVar \ge 0}
            &\quad
            \sum_{\product \in \productset} \price_\product \assignVar_\product
                \label{model:DPC} \tag{DP-C} \\
        \textnormal{s.t.}
            &\quad
                \sum_{\product \in \productset} \assignVar_{\product} = \indicator(\price_{\choice{\customer}} < \histprice_{\customer \choice{\customer}}),
                \label{eq:onechoiceCustomer_compact} \\
            &\quad
                \assignVar_{\product} \le \indicator(\price_{\product} - \price_{\choice{\customer}} \le \histprice_{\customer \product} - \histprice_{\customer \choice{\customer}}),
                &\forall \product \in \productset,
                \label{eq:indFeasibleOption}
    \end{align}
    where $\vecassignVar = (\assignVar_1, \dots, \assignVar_\numproducts)$.
\end{proposition}

\smallskip
The formulation \eqref{model:DPC} (\rev{short for Disjunctive Program-Combinatorial}) reveals the combinatorial structure of the problem for prices $\vecprice$ and the IC polyhedron of customer $\customer$.
Specifically, it is the minimum price $\price_\product$ among the feasible products according to Lemma \ref{lem:feasibleOptions}, with priority to the no-purchase option if available.
Notice that \eqref{model:DPC} is always feasible for any $\vecprice$ since either $\sum_{\product \in \productset} \assignVar_{\product}=0$ or
$\assignVar_{\choice{\customer}}=1$ is always a viable purchase option according to Lemma \ref{lem:feasibleOptions}-\ref{it:feasibleNopurchase} and \ref{it:feasibleChoice}. Furthermore, of particular importance to our methodology is the dual of \eqref{model:DPC}:
\begin{align}
    \wcfun_\customer(\vecprice)
    =
    \max_{\vecdualJ_\customer \ge 0, \vecdualS_\customer}
        &\quad
        \indicator(\price_{\choice{\customer}} < \histprice_{\customer \choice{\customer}}) \dualS_{\customer}
        -
        \sum_{\product \in \productset}
        \indicator(\price_{\product} - \price_{\choice{\customer}} \le \histprice_{\customer \product} - \histprice_{\customer \choice{\customer}}) \dualJ_{\customer \product}
            \label{model:DPDual} \tag{DP-C-Dual} \\
    \textnormal{s.t.}
        &\quad
            \dualS_\customer - \dualJ_{\customer \product} \le \price_\product,
                \quad\quad\quad\quad\quad\quad \forall \product \in \productset.
        \label{eq:feasibleDualOriginal}
\end{align}

To draw insights on the above dual problem, suppose prices are ordered as $\price_1 \le \price_2 \le \cdots \le \price_\numproducts$. If $\price_{\choice{\customer}} < \histprice_{\customer \choice{\customer}}$, then the no-purchase option is not feasible by Lemma~\ref{lem:feasibleOptions}-\ref{it:feasibleNopurchase} and
the customer necessarily purchases one product from $\productset$.
In this case, the worst-case revenue is $\price_{\product^*}$, where
$
    \product^*
    \equiv
    \min_{\product \in \productset}
    \left \{
    \product
    \colon
    \price_{\product} - \price_{\choice{\customer}} \le \histprice_{\customer \product} - \histprice_{\customer \choice{\customer}}
    \right \}.
$
In an optimal solution $(\dualS^*_\customer, \vecdualJ^*_\customer)$, variable
$\dualS^*_\customer = \price_{\product^*}$ yields the revenue obtained when prices are set to $\vecprice$. The solution $\dualJ^*_{\customer \product}$ for each $\product \in \productset$ captures the lost objective value if the product is feasible, i.e.,
\begin{align*}
\dualJ^*_{\customer \product} =
    \begin{cases}
        \price_{\product^*} - \price_\product,  & \textnormal{if $\product < \product^*$},\\
        0,                   & \textnormal{otherwise.}
    \end{cases}
\end{align*}
We note that these solutions are optimal since all are non-negative (due to the ascending order of prices), feasible to \eqref{eq:feasibleDualOriginal}, and equal to the same solution value of \eqref{model:DPC}, as only terms indexed by $\product \ge \product^*$ have non-zero objective coefficients in the dual problem. The structure of the optimal duals also implies two immediate properties that we will leverage in our reformulations.

\begin{proposition}
    \label{prop:propOptimalityDuals}
    For any $\vecprice \ge 0$, the following statements hold:
    \begin{enumerate}[label=(\alph*)]
        \item $\wcfun_\customer(\vecprice) \le \min \{\price_{\choice{\customer}},\histprice_{\customer \choice{\customer}}\}$, i.e., the optimal revenue is bounded by the new price of the historically chosen product $\choice{\customer}$.
        \label{it:boundOptPrice}

        \item Let $\maxprice \equiv \max_{\customer \in \customerset} \histprice_{\customer \choice{\customer}}$ be the maximum historical price paid by any customer. There exists some $\vecprice' \ge 0$ such that $\wcfun_\customer(\vecprice') \ge \wcfun_\customer(\vecprice)$ and $\price'_\product < \maxprice$ for all $\product \in \productset$, i.e., pricing any product higher that $\maxprice$ does not lead to any increase in revenue.
        \label{it:boundPrices}
    \end{enumerate}
\end{proposition}

Given the dual interpretation above and Proposition \ref{prop:propOptimalityDuals}, we obtain an equivalent version of the dual that is written in terms of the single variable $\dualS_\customer$:
\begin{align}
    \wcfun_\customer(\vecprice)
    =
    \max_{\dualS_\customer \ge 0}
        &\quad
        \indicator(\price_{\choice{\customer}} < \histprice_{\customer \choice{\customer}}) \dualS_{\customer}
        \label{model:DPD} \tag{DP-D} \\
    \textnormal{s.t.}
        &\quad
            \dualS_\customer \le \price_\product + \indicator(\price_{\product} - \price_{\choice{\customer}} > \histprice_{\customer \product} - \histprice_{\customer \choice{\customer}}) \histprice_{\customer \choice{\customer}},
                \quad\quad\quad\quad\quad\quad \forall \product \in \productset.
        \label{eq:feasibleDual}
\end{align}
If $\price_{\product} - \price_{\choice{\customer}} \le \histprice_{\customer \product} - \histprice_{\customer \choice{\customer}}$ for some $\product \neq \choice{\customer}$, then the $\product$-th product is feasible to purchase and inequality \eqref{eq:feasibleDual} reduces to $\dualS_\customer \le \price_\product$, i.e., the maximum revenue is bounded by $\price_\product$. Otherwise, the inequality becomes redundant given Proposition \ref{prop:propOptimalityDuals}-\ref{it:boundOptPrice}. We also note that $\dualS_\customer \le \price_{\choice{\customer}}$ for $\product = \choice{\customer}$ in \eqref{eq:feasibleDual}, i.e., $\choice{\customer}$ is always feasible
and the model is consistent with Proposition \ref{prop:propOptimalityDuals}-\ref{it:boundOptPrice}.

\subsection{Robust Pricing Reformulation}
\label{sec:milpOP}

We now develop reformulations to address our original pricing problem
\eqref{model:OP}. We consider the inner optimization of \eqref{model:OP} obtained by dropping the constant term $1/\numcustomers$ and replacing $\wcfun_\customer(\cdot)$ by \eqref{model:DPD}, which has the same objective sense as the outer problem (hence, moving from a $sup\,min$ problem to a $sup\,max$ problem, which is simply $sup$):
\begin{align}
    \sup_{\vecprice, \bm{\vecdualS} \ge 0}
        &\quad
        \sum_{\customer \in \customerset}
        \indicator(\price_{\choice{\customer}} < \histprice_{\customer \choice{\customer}}) \dualS_{\customer}
        \label{model:OPC} \tag{OP-C}
    \\
    \textnormal{s.t.}
        &\quad
        \dualS_\customer \le \price_\product + \indicator(\price_{\product} - \price_{\choice{\customer}} > \histprice_{\customer \product} - \histprice_{\customer \choice{\customer}}) \histprice_{\customer \choice{\customer}},
        \quad\quad\quad\quad\quad\quad \forall \customer \in \customerset, \product \in \productset. \label{eq:ctDualCondition}
\end{align}
Multiplying the optimal value of \eqref{model:OPC} (\rev{short for Optimal Pricing-Combinatorial}) by $1/\numcustomers$ yields the optimal revenue $\optrevenue$. The difficulty in the above problem is its non-concave and discontinuous objective function, which is defined in terms of indicator functions on both open and closed half-spaces. We will, however, exploit the constraint structure of \eqref{model:OPC} to price products to any absolute error of $\optrevenue$ by a two-step process, the first of which involves a significantly more computationally tractable model.

In particular, the indicator terms of \eqref{model:OPC} can be reformulated in several ways (see, e.g., \citealt{belotti2016handling}). For our purposes, we study the following equivalent bilinear mixed-integer model:
\begin{align}
    \sup_{\vecprice, \bm{\vecdualS} \ge 0, \vecallocVar}
        &\quad
        \sum_{\customer \in \customerset}
             \allocVar_{\customer \choice{\customer}} \dualS_\customer
        \label{model:OPB} \tag{OP-B}
    \\
    \textnormal{s.t.}
        &\quad
            \dualS_\customer \le \price_\product + (1 - \allocVar_{\customer \product}) \histprice_{\customer \choice{\customer}},
                &\forall \customer \in \customerset, \product \in \productset,
                \label{eq:OPBdualfeasible} \\
        &\quad
            \price_{\choice{\customer}} < \histprice_{\customer \choice{\customer}} + (\maxprice- \histprice_{\customer \choice{\customer}})(1 - \allocVar_{\customer \choice{\customer}}),
            &\forall \customer \in \customerset,
            \label{eq:indicatorA} \\
        &\quad
            \price_{\product} - \price_{\choice{\customer}} > \histprice_{\customer \product} - \histprice_{\customer \choice{\customer}} - (\maxprice+\histprice_{\customer \product}-\histprice_{\customer \choice{\customer}}) \allocVar_{\customer \product},
            &\forall \customer \in \customerset, \product \in \productset \setminus \{ \choice{\customer} \},
            \label{eq:indicatorB} \\
        &\quad
            \vecallocVar \in \{0,1\}^{\numcustomers \times \numproducts}.
\end{align}

\color{mycolor}
\Copy{opb-clarification}{
To see the correspondence between the variables of  \eqref{model:OPB} (\rev{short for Disjunctive Program Bilinear reformulation}) and \eqref{model:OPC},
consider a set of new prices $\{p_1,\cdots,p_n\}$. The condition $\indicator(\price_{\product} - \price_{\choice{\customer}} > \histprice_{\customer \product} - \histprice_{\customer \choice{\customer}})=0$ for $\product \in \productset \setminus \{ \choice{\customer} \}$ implies $\allocVar_{\customer \product} = 1$ due to~\eqref{eq:indicatorB}, which may be interpreted as product $\product$ being feasible for customer $\customer$ under the new prices, as per Lemma~\ref{lem:feasibleOptions}-\ref{it:feasibleH}. 
Moreover, when $\indicator(\price_{\product} - \price_{\choice{\customer}} > \histprice_{\customer \product} - \histprice_{\customer \choice{\customer}})=1$, both assignments $\allocVar_{\customer \product} = 0$ and $\allocVar_{\customer \product} = 1$ are feasible for \eqref{model:OPB}.
However, $\allocVar_{\customer \product} = 0$ leads to a (weakly) higher objective function value. 
Similarly, $\indicator(\price_{\choice{\customer}} < \histprice_{\customer \choice{\customer}})=0$ implies $\allocVar_{\customer \choice{\customer}} = 0$ because of~\eqref{eq:indicatorA}, indicating that no-purchase option is available for customer $\customer$. Further, $\indicator(\price_{\choice{\customer}} < \histprice_{\customer \choice{\customer}})=1$ implies that both assignments of $\allocVar_{\customer \choice{\customer}} = 1$ and $\allocVar_{\customer \choice{\customer}} = 0$ are feasible, while $\allocVar_{\customer \choice{\customer}} = 1$ always leads to a (weakly) higher objective function value.}

\color{black}
To elaborate more on the equivalency of \eqref{model:OPC} and \eqref{model:OPB}, we note that for both programs, the objective and constraints for each customer $\customer$ reduce to problem \eqref{model:DPD}. Further, we also remark that inequalities \eqref{eq:indicatorA}-\eqref{eq:indicatorB} are big-M constraints that rely on $\maxprice$ defined in Proposition \ref{prop:propOptimalityDuals}. For completeness, we formalize the validity of model \eqref{model:OPB}.

\begin{proposition}
    \label{prop:optimalityReformulation}
    At optimality, the solution values of \eqref{model:OPC} and \eqref{model:OPB} match.
\end{proposition}
The set of feasible solutions to \eqref{model:OPB}, however, is not polyhedral because the linear inequalities \eqref{eq:indicatorA}-\eqref{eq:indicatorB} are defined by open halfspaces. 
We propose to solve a parameterized version of \eqref{model:OPB} by setting a precision parameter $\eps$ on the constraint violation of \eqref{model:OPB}, which allow us to replace the supremum by a maximum:
\begin{align}
    \approxfun(\eps) \equiv \max_{\vecprice, \bm{\vecdualS} \ge 0, \vecallocVar}
        &\quad
        \sum_{\customer \in \customerset}
             \allocVar_{\customer \choice{\customer}} \dualS_\customer
        \label{model:OPCL} \tag{OP-$\eps$}
    \\
    \textnormal{s.t.}
        &\quad
            \dualS_\customer \le \price_\product + (1 - \allocVar_{\customer \product}) \histprice_{\customer \choice{\customer}},
                &\forall \customer \in \customerset, \product \in \productset,
                    \label{eq:OPBdualfeasible-cl} \\
        &\quad
           \price_{\choice{\customer}} \leq \histprice_{\customer \choice{\customer}} + (\maxprice- \histprice_{\customer \choice{\customer}})(1 - \allocVar_{\customer \choice{\customer}}) -  \eps,
            &\forall \customer \in \customerset,
            \label{eq:indicatorA-cl} \\
        &\quad
            \price_{\product} - \price_{\choice{\customer}} \ge \histprice_{\customer \product} - \histprice_{\customer \choice{\customer}} - (\maxprice+\histprice_{\customer \product}-\histprice_{\customer \choice{\customer}}) \allocVar_{\customer \product} + \eps,
            &\forall \customer \in \customerset, \product \in \productset \setminus \{ \choice{\customer} \},
            \label{eq:indicatorB-cl} \\
        &\quad
            \vecallocVar \in \{0,1\}^{\numcustomers \times \numproducts}
            \label{eq:integrality-cl}.
\end{align}
\label{pg:epsilon}\Copy{epsilon}{Because \eqref{model:OPCL} restricts the feasible space of \eqref{model:OPB} for $\eps>0$, the optimal value must be no greater than the one from the original problem, i.e., $\approxfun(\eps) \le \numcustomers \optrevenue$ for $\eps > 0$.
Conversely, when $\eps=0$, \eqref{model:OPCL} serves as a relaxation and therefore $\approxfun(0)\ge \numcustomers \optrevenue$.}

The challenge in solving $\approxfun(\eps)$ is to choose an appropriate $\eps > 0$ that leads to a sufficiently close approximation to the real supremum $\numcustomers \optrevenue$. Large values of $\eps$ may lead to poor approximations due to the discontinuity of \eqref{model:OP}, while small values are not computationally tractable due to numerical limitations of solvers. Theorem \ref{thm:precisionOpt}, however, shows that it suffices to solve $\approxfun(0)$ to obtain a sufficiently close value to $\numcustomers \optrevenue$ at any desired (absolute) error. It also prescribes a set of final prices with a formal guarantee that can be used by the firm.

\begin{theorem}
    \label{thm:precisionOpt}
    Consider the formulation \eqref{model:OPCL} with $\eps = 0$. The following statements hold.
    \begin{enumerate}[label=(\alph*)]
        \item The model has a finite optimal value attainable by some $\vecprice^0 \ge 0$. \label{it:thmFinite}
        \item If a product $\product \in \productset$ is priced at zero by the optimal vector $\vecprice^0$, we can reprice it at $\price^0_\product = \min_{\customer \in \customerset, \product \in \productset} \histprice_{\customer \product} > 0$ without losing optimality. \label{it:thmPositive}
        \item \label{it:thmBound} Let $\vecprice^0 > 0$ be an optimal vector of positive prices after ordering, i.e., $0 < \price^0_1 \le \price^0_2 \le \dots \le \price^0_\numproducts$. For any desired error $\precision > 0$, let $\vecprice'$ be a vector of prices such that
        \begin{align*}
            \vecprice' =
            \left(
                \price^0_1 - \frac{\precision'}{\numcustomers \numproducts}, \;
                \price^0_2 - 2\frac{\precision'}{\numcustomers \numproducts}, \;
                \price^0_3 - 3\frac{\precision'}{\numcustomers \numproducts}, \;
                \dots,
                \price^0_\numproducts - \numproducts \frac{\precision'}{\numcustomers \numproducts}
            \right)
        \end{align*}
        for any sufficiently small $0 < \precision' \le \precision$ so that $\vecprice' > 0$. Then,
        \begin{align*}
            0 \le \approxfun(0) - \numcustomers \optrevenue \le \approxfun(0) - \sum_{\customer \in \customerset} \wcfun_{\customer}(\vecprice') \le \precision,
        \end{align*}
        where $\optrevenue$ is the optimal solution of the original pricing problem \eqref{model:OP}.
    \end{enumerate}
\end{theorem}

\medskip
By Theorem~\ref{thm:precisionOpt}, we are able to solve \eqref{model:DP} for any desired error $\precision$ using the reformulation \eqref{model:OPCL} with $\eps=0$.
More precisely, the set of prices $\vecprice^*$ can be obtained by
\begin{enumerate}
    \item Solving \eqref{model:OPCL} for $\vecprice$ with $\eps = 0$, which is guaranteed to exist due to Theorem \ref{thm:precisionOpt}-\ref{it:thmFinite};
    \item If any price is zero, increasing it to a positive value according to Theorem \ref{thm:precisionOpt}-\ref{it:thmPositive}; and
    \item Applying the transformation from Theorem \ref{thm:precisionOpt}-\ref{it:thmBound} to any desired error $\precision$.
\end{enumerate}

\medskip
While modern commercial solvers can address nonlinear models of the form \eqref{model:OPCL}, we also present an equivalent mixed-integer linear program by a standard big-M reformulation of the quadratic constraints. The program will be the key to the analysis of our approximation algorithms.
\begin{align}
    \approxfun(0) = \max_{\vecprice, \bm{\vecdualS}, \bar{\bm{\vecdualS}} \ge 0, \vecallocVar}
        &\quad
        \sum_{\customer \in \customerset}
             \bar{\dualS}_{\customer}
        \label{model:OPMIP} \tag{OP-MIP}
    \\
    \textnormal{s.t.}
        &\quad
            \textnormal{\eqref{eq:OPBdualfeasible-cl}, \eqref{eq:indicatorA-cl}, \eqref{eq:indicatorB-cl}, \eqref{eq:integrality-cl} with $\eps = 0$} \\
        &\quad
            \bar{\dualS}_\customer \le \allocVar_{\customer \choice{\customer}} \histprice_{\customer \choice{\customer}},
            &\forall \customer \in \customerset, \\
        &\quad
            \bar{\dualS}_\customer \le \dualS_\customer,
            &\forall \customer \in \customerset, \\
        &\quad
            \bar{\dualS}_\customer \ge \dualS_\customer - (1 - \allocVar_{\customer \choice{\customer}})\histprice_{\customer \choice{\customer}},
            &\forall \customer \in \customerset.
\end{align}

\color{mycolor}
\noindent \label{pg:market-share}\Copy{market-share}{\textit{Prioritizing market share.} We note that in some market environments, the firm may want to ensure a certain  purchase probability for the new customer through its pricing, in addition to revenue maximization. 
The optimization program~\eqref{model:OPMIP} can be modified accordingly to tackle such a setting with market penetration considerations. 
For example, to guarantee a purchase probability of at least $\rho$ for the new customer, it suffices to add the constraint $\sum_{\customer\in \customerset} \allocVar_{\customer \choice{\customer}}\ge \numcustomers\rho$ to \eqref{model:OPMIP}.
Thus, our approach can accommodate the case when the firm prioritizes the market share consideration in its revenue maximization objective.}
\color{black}

\smallskip

%
%
%

\subsection{Special Cases}
\label{sec:special_cases}

We now discuss two cases of practical interest where \eqref{model:OPMIP} can be solved analytically. First, we consider the scenario where each individual customer is offered the same price for all products in the assortment. The prices, however, can be different per customer. This occurs when products are similar in nature and customers are offered personalized promotions over time. Proposition \ref{prop:same-price} states the structure of the optimal solution for this case.
\begin{proposition}
    \label{prop:same-price}
    Suppose that, for each customer $\customer \in \customerset$, all products $\product \in \productset$ have the same historical price $\histprice_{\customer \product} = \histprice_\customer$. Furthermore, without loss of generality, assume prices are ordered, i.e., $\histprice_1 \le \histprice_2 \le \dots \le \histprice_{\numcustomers}$. The price vector $\vecprice^*$ defined by $\price^*_\product = \histprice_{\customer^*}$ for all $\product \in \productset$, where
    $
        \customer^* = \argmax_{\customer \in \customerset} \left\{ (\numcustomers - \customer + 1)\histprice_\customer \right\},
    $
    is optimal to \eqref{model:OPMIP}.
\end{proposition}

\smallskip
Proposition \ref{prop:same-price} reveals a connection with the classical pricing literature. Specifically, if we perceive
the set $\{\histprice_1, \dots, \histprice_\numcustomers\}$ as the empirical distribution of the customer valuations, then our problem, when all products are offered at the same price, reduces to the well-studied revenue maximization problem $\max_\price \{ \price \cdot d(\price) \}$, where $d(\price)$ is the demand function of price $\price$.

As our second practical case, we consider a setting where the price for any product is fixed over time. More precisely, prices may differ per product but not per customer. Proposition \ref{prop:constant-price} also indicates that the optimal prices have a simpler structure, i.e., 
it suffices to set them to their historical prices.

\begin{proposition}
    \label{prop:constant-price}
    Suppose that, for each product $\product \in \productset$, all customers observe the same historical price $\histprice_\product$.
    The price vector $\vecprice^*$ such that $\price^*_\product = \histprice_\product$ for all $\product \in \productset$ is optimal to \eqref{model:OPMIP}.
\end{proposition}

\section{Approximate Pricing Strategies and Analysis}
\label{sec:heuristics}

 The formulation \eqref{model:OPMIP}, albeit amenable to state-of-the-art commercial solvers, may still be challenging to solve due to its difficult constraint structure (e.g., the presence of big-M constraints). Moreover, \cite{abdallah2020demand} point out that in retail, the number of SKUs in a family of products could be on the order of several hundred. Thus, even with a few data points per product, the size of \eqref{model:OPMIP} could be significantly large. While recent works have focused on choice model estimation in such high-dimensional settings \citep{jiang2020high}, to the best of our knowledge, model-free pricing approaches are yet to be developed for such settings.

In this section, we analyze three interpretable and intuitive approximation algorithms that are of low-polynomial time complexity in the input size of the problem, and hence scalable to large problem sizes. We discuss their benefits and worst-case revenue performance in comparison to the optimal solution of \eqref{model:OPMIP}. Specifically, we evaluate in \S\ref{sec:heu-conservative} and \S\ref{sec:heu-LP-relax} two standard heuristics based on historical prices and linear programming (LP) relaxations, respectively. In \S\ref{sec:V0-conservative}, we propose an approximation based on a simplified version of \eqref{model:OPMIP}, which provides the strongest approximation factor among the three policies, is efficient to compute, and also interpretable. 

We note that in practice the firm might want to use simpler pricing algorithms, such as offering the products at their average historical prices or at the prices observed by a random historical customer.
However, it can be shown that such pricing algorithms can lead to an arbitrarily poor revenue performance in comparison to the optimal solution of \eqref{model:OPMIP}.
See Propositions~\ref{prop:average prices} and~\ref{prop:random prices} in the e-companion.
Thus, we do not discuss them in details in this section.

\subsection{Conservative Pricing}
\label{sec:heu-conservative}

A simple approximation a conservative firm may consider is to price all the products at their historically lowest purchase price. That is,
\begin{align}
    \label{eq:conservative-price}
    \price_\product
        = \min_{\customer \in \customerset \colon \choice{\customer} = \product} \histprice_{\customer \product},
    \quad\quad \forall \product \in \productset.
\end{align}
By inequality \eqref{eq:OPBdualfeasible-cl} for $\product = \choice{\customer}$, the above prices guarantee that each customer $\customer$ would purchase at least one product (i.e., the no-purchase option will not be chosen). Furthermore, it also follows that the final total revenue is at least $\numcustomers \ubar{\histprice}$, where $\ubar{\histprice} \equiv \min_{\customer \in \customerset} \histprice_{\customer \choice{\customer}}$.

Such a pricing policy maximizes the demand at the cost of a lower profit margin, and is thus referred to as ``conservative pricing.''
We show below that it has a straightforward worst-case performance bound, which is also tight.
\begin{proposition}
    \label{prop:consv-pricing}
    Let $\ubar{\histprice} \equiv \min_{\customer \in \customerset} \histprice_{\customer \choice{\customer}}$ and
    $\bar{\histprice} \equiv \max_{\customer \in \customerset} \histprice_{\customer \choice{\customer}}$ be the minimum and maximum historical purchase prices, respectively. The revenue from the conservative pricing \eqref{eq:conservative-price} is at least
    $\ubar{\histprice} / \bar{\histprice}$ of the optimal value of \eqref{model:OPMIP}. Futhermore, this ratio is asymptotically tight, \textcolor{mycolor}{as the worst-case is achieved when the number of historical customers grows to infinity}.
\end{proposition}

As one may expect, when products are not similar in nature and their prices vary in a wide range, that this approximation is too conservative and would not perform well. However, the ratio $\ubar{\histprice}/\bar{\histprice}$ provides a useful benchmark which can be used to gauge the performance of other approximations.

\subsection{LP Relaxation Pricing}
\label{sec:heu-LP-relax}
It is a natural practice to consider the LP relaxation of the program \eqref{model:OPMIP}, i.e., to replace the integrality constraint \eqref{eq:integrality-cl} by the continuous domain $\vecallocVar \in [0,1]^{\numcustomers \times \numproducts}$. The resulting LP model can be solved in (weakly) polynomial time, from which we can extract a candidate price vector $\vecprice^{\textnormal{LP}}$. We wish to evaluate the quality of $\vecprice^{\textnormal{LP}}$ with respect to the optimal prices of \eqref{model:OPMIP}. Notice that the remaining variables of \eqref{model:OPMIP} can be easily determined when $\vecprice$ is fixed to $\vecprice^{\textnormal{LP}}$.

In Example~\ref{exp:lp-fail} in the e-companion we show that the resulting LP prices, unfortunately, have a worse worst-case performance than the conservative pricing approach. However, even though the worst-case revenue is not superior relative to conservative pricing, we show in \S\ref{sec:numerics} that the LP relaxation pricing usually outperforms conservative pricing numerically in most cases. 
Intuitively, conservative pricing is more concerned with the worst-case, while LP relaxation can be close to maximizing the expected performance, specifically since we have tuned the big-M constrained to be very tight with respect to the input parameters.

\subsection{Cut-off Pricing}
\label{sec:V0-conservative}
In this section, we provide a heuristic by manipulating and reshaping the IC polyhedra.
In particular, consider the historical customer $\customer$ in \eqref{model:DPC}.
Suppose we omit constraint \eqref{eq:indFeasibleOption} when
deciding which product to buy under $\vecprice$ in the worst case.
In other words, the customer does not fully follow the IC constraints; rather, as long as the price of the historically chosen product $\choice{\customer}$ is lower than its historical price, i.e.,  $\price_{\choice{\customer}} \le \histprice_{\customer \choice{\customer}}$,
all products are eligible for purchase.
Therefore, if the customer decides to purchase a product (i.e., $\price_{\choice{\customer}} \le \histprice_{\customer \choice{\customer}}$), then she chooses the product with the lowest price, i.e., product $\argmin_{\product \in \productset}\price_\product$, in the worst case.

We formulate this setting in the following model:
\begin{align}
\max_{\vecprice, \bm{\vecdualS} \ge 0}
\left\{
    \sum_{\customer \in \customerset} \indicator(\price_{\choice{\customer}} \le \histprice_{\customer \choice{\customer}}) \dualS_\customer
    \; \colon \;
    \dualS_\customer \le \price_\product, \;\; \forall \product \in \productset
\right\}.
\label{model:OPCP} \tag{OP-CP}
\end{align}
Under the assumed purchase rule, all historical customers would choose the same product with the minimum price, as long as they decide to purchase under $\vecprice$. 
Thus, the key is to determine the minimum price $p^*$. 
Then setting all products at this price ($\price_\customer\equiv p^*$) maximizes the objective in \eqref{model:OPCP}. 
As we vary the value of $p^*$, the indicator functions only change values when $p^*=\histprice_{\customer\choice{\customer}}$.
Letting $p^*=\histprice_{\customer\choice{\customer}}$ for some $\customer \in \customerset$,
formulation \eqref{model:OPCP} simplifies to 
\begin{align}
    \label{eq:modelOCPCsimplification}
    \max_{\customer \in \customerset} \sum_{\customer' \in \customerset} \indicator(\histprice_{\customer' \choice{\customer'}} \ge \histprice_{\customer \choice{\customer}}) \; \histprice_{\customer \choice{\customer}},
\end{align}
which can be solved in $\bigo(\numcustomers)$ time complexity by inspecting one $\histprice_{\customer\choice{\customer}}$ at a time. We denote the optimal solution to \eqref{model:OPCP} by $\price^* = \histprice_{\customer \choice{\customer}}$ for some customer $\customer$. 
In particular, $\price^*$ can be perceived as a cut-off price: 
the historical customer $\customer'$ would not purchase any product if and only if the historical price paid by customer $\customer'$, $\histprice_{\customer'\choice{\customer'}}$, is below $\price^*$. 
We leverage this to propose the cut-off pricing approximation $\vecprice^{\textnormal{CP}}$ as follows. For 
every $\product \in \productset$,
\begin{align}
    \label{eq:cutoff-price}
    \price^{\textnormal{CP}}_\product 
    =
    \begin{cases}
    \min_{\customer \in \customerset}
        \left\{
            \histprice_{\customer \choice{\customer}}
            \colon
            \choice{\customer} = \product, \; 
            \histprice_{\customer \choice{\customer}} \ge \price^*
        \right\},
        & \textnormal{if $\{\customer \in \customerset \colon \choice{\customer} = \product, \histprice_{\customer \choice{\customer}} \ge \price^* \} \neq \emptyset$},
        \\
        \min\{\max\{\max_{\customer\in\customerset}\histprice_{\customer\product},\price^*\},\bar{\histprice}\}, & \textnormal{otherwise},
    \end{cases}
\end{align}
where each product $\product$ is priced at its lowest historical purchase price that was greater than or equal to the cut-off price $p^*$.
We note that not all products are priced at the cut-off price.
Rather, they are often priced slightly higher than $p^*$.
Such modifications do not impact the indicator functions and thus the optimal value of \eqref{model:OPCP} but lead to better and less conservative empirical performances.

\label{pg:prop10}\Copy{prop10}{We next show its performance in Proposition \ref{prop:cutoff-pricing-bound}, recalling that $\ubar{\histprice} \equiv \min_{\customer \in \customerset} \histprice_{\customer \choice{\customer}}$ and
    $\bar{\histprice} \equiv \max_{\customer \in \customerset} \histprice_{\customer \choice{\customer}}$ are the minimum and maximum historical purchase prices, respectively, and introducing $med(P)$ and $avg(P)$ as the median and mean of $\{\histprice_{\customer \choice{\customer}}\}_{\customer=1}^\numcustomers$.}

\begin{proposition}
\label{prop:cutoff-pricing-bound}
The cut-off pricing \eqref{eq:cutoff-price} generates a revenue that is at least $\max\{\frac{1}{1+\log(\bar{P}/\ubar{P})}, \frac{med(P)}{2avg(P)}\}$ of the optimal value of \eqref{model:OPMIP}. Furthermore, this bound is asymptotically tight, \textcolor{mycolor}{as it is achieved when both the number of products and number of historical customers grow to infinity}.
\end{proposition}
\color{mycolor}
\label{pg:prosketch2}\Copy{prosketch2}{
\begin{proof}{Proof sketch.}
We create a random variable $\rv$ whose support is on the set $\{\histprice_{1\choice{1}},\dots, \histprice_{\numcustomers\choice{\numcustomers}}\}$.
We show that the performance bound of cut-off pricing
can be translated to $\max_x x(1-\cdf(x))/\mathbb{E}[\rv] $, where $\cdf(x)$ is the CDF of $\rv$.
This is a classic problem and leads to the bound $\max \left \{\frac{1}{1+\log(\bar{P}/\ubar{P})}, \frac{med(P)}{2avg(P)} \right\}$.

To show that the performance bound is asymptotically tight for any $\bar{P}$ and $\ubar{P}$, we construct an example with $\numproducts$ products and $\numcustomers$ customers that achieves a ratio arbitrarily close to $\frac{1}{1+\log(\bar{P}/\ubar{P})}$, when both $\numproducts$ and $\numcustomers$ grow arbitrarily large. Intuitively, this worst-case performance happens when a small decrease in the price of historically lower priced products leads to all the historical customers with high historical purchase prices to choose them, hence showcasing the limitation of cut-off pricing in assuming the customer does not follow the IC constraints. \hfill $\square$ 
\end{proof}}
\color{black}
\smallskip

Compared to Proposition~\ref{prop:consv-pricing}, the cut-off pricing dramatically improves upon the worst-case scenario of the conservative pricing, especially when $\bar P\gg \ubar P$. \color{mycolor} \label{pg:disperevidence}\Copy{disperevidence}{While Proposition~\ref{prop:cutoff-pricing-bound} implies a strong performance bound for this pricing policy, we have shown that for any $\ubar{P}$ and $\bar{P}$ we can construct an asymptotic example such that the revenue from cut-off pricing is arbitrarily close to $\frac{1}{1+\log(\bar{P}/\ubar{P})}$. Thus, we note that when $\ubar{P}/\bar{P}$ is near zero, in theory, the performance of cut-off pricing could be poor. However, this is unlikely to happen in practical settings, 
primarily because Proposition~\ref{prop:cutoff-pricing-bound} suggests that cut-off pricing always generates a revenue with a factor that is at least a half of the ratio of the median to the mean of the historical purchase prices.
Unless the purchase prices are highly skewed, this ratio is likely to be close to $1/2$. 

Empirical observations of the price dispersion $\ubar{P}/\bar{P}$ (see, e.g.,  \citealt{hosken2004patterns,anania2014price,dubois2015price}) often lead to a reasonable bound by Proposition~\ref{prop:cutoff-pricing-bound}. 
For example, in a study of grocery prices across the United States, \cite{hosken2004patterns} investigate the frequency distribution of scaled prices for $20$ categories of goods (the price of each good is scaled by its annual modal price) and show that the entirety of this distribution lies within $[0.6,1.4]$.
Moreover, in the e-companion~\ref{app:tables}, we demonstrate bounds computed from  $31$ different categories in the IRI academic dataset.
We show that for almost all the categories, the bound implied by Proposition~\ref{prop:cutoff-pricing-bound} falls into the range of 30\% to 50\%.} 
Finally, we also observe in \S\ref{sec:numerics} that cut-off pricing is superior numerically to the the two earlier proposed heuristics. 
\rev{
\begin{remark}\label{cutoffs:signleprod}
In the case of model-free pricing with a single product, cut-off pricing recovers the optimal price. This follows from Proposition~\ref{prop:same-price} and the fact that the cut-off price achieves the optimum of \eqref{eq:modelOCPCsimplification}. 
\end{remark}}

\color{black}  


\section{Numerical Analysis}
\label{sec:numerics}

We now present a numerical study of the proposed approaches on both synthetic and real datasets, evaluating our methodologies with respect to classical model-based methodologies in scenarios of practical interest. We start in \S
\ref{sec:approximationPerformance} with an analysis of the empirical performance of the approximation pricing strategies from \S\ref{sec:heuristics}. In \S \ref{sec:small-data-num}, we evaluate our data-driven pricing approach on ``small-data'' regimes, which are typically challenging to classical model-based methods. Next, we consider a scenario in \S\ref{sec:model-misspec} where the firm misspecifies the pricing model. Finally, in \S\ref{sec:real-data} we compare all approaches on a large-scale dataset from the U.S. retail industry.



\subsection{Approximation Performance}
\label{sec:approximationPerformance}

In this section, we use synthetic data to investigate the performance of the approximation algorithms developed in Section~\ref{sec:heuristics}. We generate instances with $\numcustomers \in \{50, 100, 150, 200\}$, $\numproducts\in\{10, 15, 20, 25\}$, and historical prices $\histprice_{\customer\product}$ drawn uniformly at random from the interval $(0,10)$. 
Each customer chooses a product in the assortment with equal probability $1/\numproducts$.
Based on the synthetic data, we compute the optimal value $\approxfun(0)$ from \eqref{model:OPMIP}, and the objective of the three approximations in Section~\ref{sec:heuristics}.
We report the objective values of the approximations relative to the optimal value $\approxfun(0)$ for 200 independent instances.

\begin{figure}
\centering
\caption{The performance of the three approximation algorithms.}
\label{fig:appcomparecust}
\vskip1em
\begin{minipage}[b][5.5cm][s]{.45\textwidth}
\centering
\begin{tikzpicture}
\pgfplotsset{compat=newest}
\begin{axis}[width=5.5cm, xlabel={Number of customers $\numcustomers$}, ylabel={Relative performance}, ymax=1, ymin=0, xmin=50, xmax=200, yticklabel=\pgfmathparse{100*\tick}\pgfmathprintnumber{\pgfmathresult}\,\%,
        yticklabel style={/pgf/number format/.cd,fixed,precision=2}, legend pos=outer north east]
\addplot+[only marks]
table[x =n, y=value3]{data2.txt};
\addlegendentry{Cut-off pricing}
\addplot+[only marks]
table[x =n, y=value2]{data2.txt};
\addlegendentry{LP Relaxation}
\addplot+[mark=x, only marks]
table[x =n, y=value]{data2.txt};
\addlegendentry{Conservative}
\end{axis}
\end{tikzpicture}
\end{minipage}\hspace{.5cm}
\begin{minipage}[b][5.5cm][s]{.45\textwidth}
\centering
\begin{tikzpicture}
\begin{axis}[width=5.5cm, xlabel={Number of products $\numproducts$}, ymax=1, ymin=0, xmin=10, xmax=25, yticklabel=\pgfmathparse{100*\tick}\pgfmathprintnumber{\pgfmathresult}\,\%,
        yticklabel style={/pgf/number format/.cd,fixed,precision=2}, legend pos=outer north east]
\addplot+[only marks]
table[x =n, y=value3]{data3.txt};
\addplot+[only marks]
table[x =n, y=value2]{data3.txt};
\addplot+[only marks, mark=x]
table[x =n, y=value]{data3.txt};
\end{axis}
\end{tikzpicture}
\end{minipage}\qquad
\end{figure}
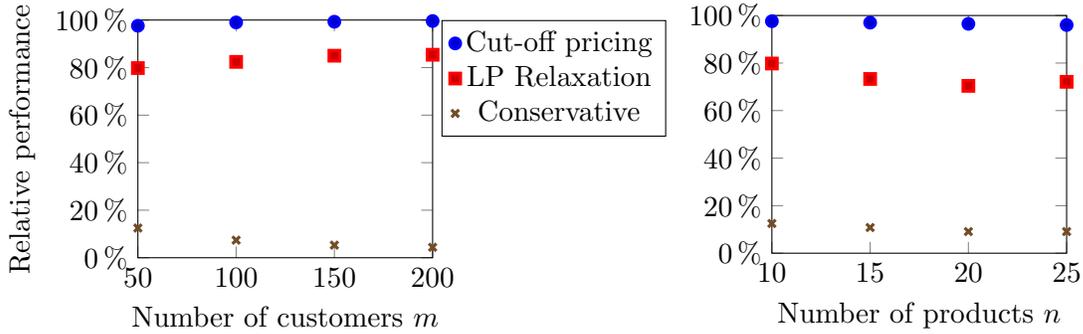

Figure~\ref{fig:appcomparecust} depicts the relative performance ratio in percentage (i.e., $100 \times$optimal revenue/heuristic revenue) of the conservative, LP relaxation, and cut-off pricing when the number of customers $\numcustomers$ and the number of products $\numproducts$ vary in the historical data.
The figures suggest that conservative pricing performs poorly (achieving less than 15\% of the optimal revenues on average) and cut-off pricing does the best among the three, obtaining at least 96\% of the optimal value in all the cases. Increasing the number of customers or decreasing the number of products improves the performance of LP relaxation pricing and cut-off pricing. Moreover, the numerical results suggest that cut-off pricing significantly outperforms  conservative pricing, as expected from its theoretical performance guarantee (Proposition~\ref{prop:cutoff-pricing-bound}), as cut-off pricing generally has an average performance that is an order of magnitude higher than that of conservative pricing. 
Finally, we observe that the numerical  performance by cut-off pricing is far above that of its theoretical performance; we include additional tables in the e-companion \ref{app:tables} with the expected theoretical performance. 

Table~\ref{tab:soln-time} shows the solution time of the conservative, LP relaxation, cut-off pricing, and optimal problem when the number of customers $\numcustomers$ and the number of products $\numproducts$ vary in the historical data. The results suggest that while the optimal solution time for mid-size problems is not too egregious, however it does not scale well when the problem size grows, as is expected. However, cut-off pricing pricing solves all problem sizes in under one second.


\ctable[
  notespar,
  caption={The average solution time of the approximation strategies and optimal solution. Standard errors are reported in parentheses.}
  \label{tab:soln-time},
doinside =\small
]{r c @{\hspace*{6mm}} c @{\hspace*{6mm}} c@{\hspace*{6mm}} c@{\hspace*{6mm}} r}
{{} }{
\toprule
  &\multicolumn{4}{c}{Solution time in seconds}  & \\\cmidrule(r){2-6}
$(\numcustomers,\numproducts)$ & \text{Conservative} & \text{LP Relaxation}  & \text{Cut-off}  & \text{Optimal}   \\
\midrule

$(50,10)$         &  $0.041$  ($0.001$)  & $0.453$   ($0.001$)   & $0.040$ ($0.001$) & $1.477$ ($0.078$) \\

$(50,15)$ &  $0.065$  ($0.001$)  & $1.049$   ($0.005$)   & $0.062$ ($0.001$) & $4.250$ ($0.195$) \\

$(50,20)$ &    $0.080$  ($0.001$)  & $1.762$   ($0.006$)   & $0.083$ ($0.001$) & $5.691$ ($0.254$) \\

$(50,25)$ &   $0.103$  ($0.001$)  & $2.767$   ($0.007$)   & $0.109$ ($0.001$) & $8.003$ ($0.313$) \\

$(100,10)$ & $0.081$ (0.001)  &   $0.929$ (0.002)  & $0.078$ ($0.001$) & $17.260$ ($0.800$)   \\

$(150,10)$     & $0.117$ (0.001)  &   $1.419$ (0.003)  & $0.132$ (0.002) & $89.050$ ($5.400$) \\

$(200,10)$   & $0.165$ (0.002)  &   $2.018$ (0.006)  & $0.181$ (0.003) & $525.800$ ($42.930$) \\
\bottomrule}

\subsection{Small Sample Size}
\label{sec:small-data-num}
When the data size is small, model-based methods risk unstable estimations even when the model is correctly specified. In this section, we evaluate revenues obtained from the proposed pricing approaches in such cases, comparing with a classical multinomial logit (MNL) model \citep{train2009discrete}. 

We generate synthetic instances with $\numproducts=10$ products and varying number of customers $\numcustomers \in \{20,30,\dots,110\}$. 
The probability of customer $\customer$ choosing product $\product$ from the assortment is given by
\begin{equation*}
    \frac{\exp(\alpha_\product-\beta \histprice_{\customer\product})}{1+\sum_{k=1}^\numproducts \exp(\alpha_k-\beta \histprice_{\customer k})},
\end{equation*}
where $\beta = 0.5$. 
So the number of data points $m$ is merely enough to estimate 10 parameters.
We consider two sets of experiments for the remaining parameters $\alpha$ and the historical prices, representing low and high customer utility. 
Specifically, for the high-utility experiment, historical prices $\histprice_{\customer\product}$ are drawn uniformly at random from the interval $[5.5,8.5]$, and  
$\{\alpha_\product\}_{\product=1}^{10}$ are independent and drawn uniformly at random from the interval $[1,3]$. In the low-utility experiment, historical prices are drawn uniformly at random from the interval $[2.5,4.5]$ and 
$\{\alpha_\product\}_{\product=1}^{10}$ are independent and drawn uniformly at random from the interval $[-2,0]$. 
Table~\ref{tab:scenarios} summarizes the setup of the two experiments.
Note that, in both cases, this choice of parameters guarantees that the optimal price of the MNL model is within the range of the historical prices. \rev{In Section~\ref{app:tables} of the e-companion we also study the effect of the dispersion of historical prices in the data on the comparative performance of our data-driven framework.}

\begin{center}
\begin{minipage}{.4\textwidth}
\ctable[
  notespar,
  caption={Instance parameters for study in \S\ref{sec:small-data-num}.}
  \label{tab:scenarios}, pos=H,
doinside =\small
]{l c c c }
{\tnote[]{}\\ \tnote[]{}\\}{
\toprule
\text{Experiment} & \text{$\alpha$} & \text{$P_{ij}$} \\ 
 
\midrule
Low utility         &  $[-2,0]$  & $[2.5,4.5]$       \\

High utility & $[1,3]$  &   $[5.5,8.5]$ \\
\bottomrule}

\end{minipage}
\begin{minipage}{.4\textwidth}
\ctable[
  notespar,
  caption={Instance parameters for the studies in \S\ref{sec:model-misspec}.}
  \label{tab:scenariosmis},  pos=H,
doinside =\small
]{l c c c }
{\tnote[]{}\\ \tnote[]{}\\}{
\toprule
\text{Experiment} & \text{$\alpha$} & \text{$P_{ij}$} \\ 
 
\midrule
Low utility         &  $[-2,0]$  & $[2.5,4.5]$       \\

High utility & $[1,3]$  &   $[5.5,8.5]$ \\

\bottomrule}
\end{minipage}
\end{center}

We simulate 200 independent instances for each $\numcustomers$. In each instance, we calculate the optimal solution to \eqref{model:OPMIP} and cut-off pricing.
Moreover, we use the BIOGEME package developed by \cite{bierlaire2003biogeme} to estimate the parameters of the MNL model with the historical data, and then calculate the optimal price based on the fitted model.
Note that they are the prices obtained from the model-estimate-optimize approach under the correct model specification.
We then evaluate the three sets of prices with respect to the ground-truth model and compare their expected revenues.

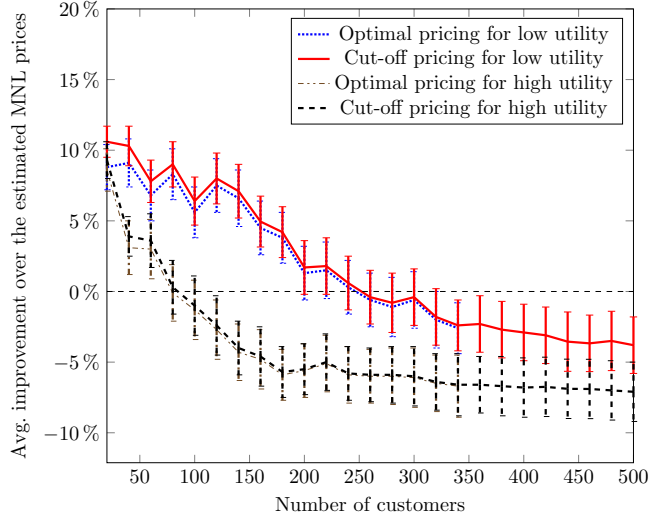
\begin{figure}
\centering
\caption{The performance of the data-driven optimal pricing \eqref{model:OPMIP} and cut-off pricing relative to the optimal MNL prices estimated from the data. The difference in the revenues is converted to percentage by dividing it by the optimal revenue of the model before being averaged.}
\label{fig:lowdata}
    \scalebox{.7}{\begin{tikzpicture}
    \pgfplotsset{compat=newest, ymax=.3}
    \begin{axis}[width=10cm, scale only axis,  xlabel={Number of customers}, ylabel={Avg. improvement over the estimated MNL prices}, xmin = 20, xmax=500, ymax = 0.2, yticklabel=\pgfmathparse{100*\tick}\pgfmathprintnumber{\pgfmathresult}\,\%, clip=false,
            yticklabel style={/pgf/number format/.cd,fixed,precision=2}]
    \addplot+[mark=none, densely dotted]
    plot [error bars/.cd, y dir = both, y explicit]
    table[x=n, y=value3, y error =error3]{dataoptim.txt};
    \addlegendentry{Optimal pricing for low utility}
    \addplot+[mark=none]
    plot [error bars/.cd, y dir = both, y explicit]
    table[x =n, y=value4, y error =error4]{data.txt};
    \addlegendentry{Cut-off pricing for low utility}
    \addplot+[dash dot dot, mark=none, thin] 
    plot [error bars/.cd, y dir = both, y explicit]
    table[x =n, y=value, y error =error]{dataoptim.txt};
    \addlegendentry{Optimal pricing for high utility}
    \addplot+[dashed, mark=none]  
    plot [error bars/.cd, y dir = both, y explicit]
    table[x =n, y=value2, y error =error2]{data.txt};
    \addlegendentry{Cut-off pricing for high utility}
    
    \addplot[ultra thin, dashed, black] coordinates {(20,0) (500,0)};
    \end{axis}
    \end{tikzpicture}}
\end{figure}

In Figure~\ref{fig:lowdata}, we illustrate the average difference of the revenues of cut-off (data-driven optimal) pricing and the estimated MNL prices, relative to the optimal MNL prices when the parameters are known. 
For the optimal pricing policy, we do not compute it for more than 350 customers due to the prohibitive computational cost.
Note that if the quantity is positive, then it implies that our approach outperforms the estimated MNL prices. 
From the figure, when the number of customers is less than 70, both data-driven pricing schemes outperform the estimated MNL prices, in both experiments.
Note that the MNL prices are estimated based on the correct specification of the model.
This experiment further suggests that data-driven approaches may be beneficial with respect to model-based approaches when the data size is small.
It is also expected that as the data size grows, the model-estimate-optimize approach eventually converges to the optimal prices of the true model, as the model is correctly specified.
In this regime, data-driven approaches may not be beneficial. \rev{We note the MNL prices are estimated from uncensored data, favoring the model based approaches and not affecting the data-driven prices (Remark~\ref{rem:nopurchasecustomers}). In Figure~\ref{fig:censor} and Section~\ref{sec:censresults} in the e-companion we investigate the case when the MNL prices are estimated from censored data.}

\subsection{Model Misspecification}
\label{sec:model-misspec}

\Copy{misspecification1}{Another potential benefit of the proposed approach is that it is agnostic to the underlying model, and thus less sensitive to model misspecification.} We assess this scenario in the next experiment with synthetic instances. 
We consider $\numcustomers=50$ customers and $\numproducts=10$ products in two experiment sets, also representing low and high utilities.
The historical prices $\histprice_{\customer\product}$ are drawn uniformly at random from $[2.5,4.5]$ and $[5.5,8.5]$ in the two experiments.
For each customer, we generate their choices using a mixed logit model with two classes \citep{train2009discrete}. More precisely, given $(\histprice_{\customer 1},\dots,\histprice_{\customer \numproducts})$, the probability of customer $\customer$ choosing product $\product$ is
\begin{equation*}
    \frac{1}{2}\cdot\frac{\exp(\alpha_{1\product}-\beta_1 \histprice_{\customer\product})}{1+\sum_{k=1}^\numproducts \exp(\alpha_{1k}-\beta_1 \histprice_{\customer k})}+\frac{1}{2}\cdot\frac{\exp(\alpha_{2\product}-\beta_2 \histprice_{\customer\product})}{1+\sum_{k=1}^\numproducts \exp(\alpha_{2k}-\beta_2 \histprice_{\customer k})}.
\end{equation*}
Here we set $\beta_1=0.5$ and $\beta_2=2$.
In each of the 200 instances, we randomly draw  $(\alpha_{1\product},\alpha_{2\product})$ independently from $[-2,0]$ (for $\histprice_{\customer\product}\in [2.5,4.5]$) and $[1,3]$ (for $\histprice_{\customer \product}\in[5.5,8.5]$). The range choices of $\histprice$ and $\alpha$ guarantee that the historical price ranges cover the optimal MNL with $\beta=0.5$ prices (Table~\ref{tab:scenariosmis}). 

We investigate the case that the model is misspecified.
In particular, we fit the MNL model instead of the mixed logit model using BIOGEME to the historical data.
We calculate the optimal prices for the fitted MNL model and compare its expected revenue to our data-driven approaches (optimal solution and cut-off pricing) under the mixed logit model.


Table~\ref{alg:compare:miss} suggests a better performance of our model-free approach as the optimal prices of the MNL model are designed for a misspecified model.
\Copy{misspecification2}{Note that in this case, the misspecification error between the MNL model and the mixed logit model with two classes is arguably mild.
We may expect the benefit of model-free assortment pricing to be more substantial when there are strong irregular patterns in the data that cannot be captured by the assumed model.}
Combined with the results in \S\ref{sec:small-data-num},
this further suggests that our proposed pricing models could be beneficial when the data size is small or the firm has little confidence in the modeling of the demand.



\ctable[
  notespar,
  caption={Revenue between data-driven pricing and misspecified MNL model.}
  \label{alg:compare:miss},
doinside =\small
]{l c @{\hspace*{6mm}} c @{\hspace*{6mm}} c@{\hspace*{6mm}} c}
{{} }{
\toprule
  &\multicolumn{2}{c}{Expected revenue in the mixed logit model}  & \\\cmidrule(r){2-4}
\vtop{\hbox{\strut \text{Pricing method}}} & \text{Low-utility experiment}  & \text{High-utility experiment}   \\
\midrule
Data-driven optimal         &  $0.725$ ($0.005$)   & $2.393$  ($0.007$)       \\

Cut-off pricing & $0.734$ (0.005)  &   $2.415$ (0.007)  \\

MNL optimal pricing &       $0.635$ (0.014) &       $2.113$ (0.035) \\

\bottomrule}

\subsection{Real Datasets}\label{sec:real-data}

\label{pg:whyrazor}\Copy{whyrazor}{In this section, we apply model-free assortment pricing to the IRI Academic Dataset \citep{Kruger2009}.
The IRI data collects weekly transaction data from 47 U.S. markets from 2001 to 2012, covering more than 30 product categories.
Each transaction includes the week and store of the purchase, the universal product code (UPC) of the purchased item, the number of units purchased and the total paid dollars.
We investigate the product category of razors and the transactions from the first two weeks in 2001. \textcolor{mycolor}{We focus on this category primarily because it forms a proper assortment, i.e., customers are unlikely to purchase more than one unit, if they purchase any.} To construct the assortments, we focus on the top ten (out of $45$) purchased products from all stores during the two weeks, that is, $\numproducts=10$.}
The purchases of all other products are treated as ``no-purchase.''
An assortment is thus defined as the products of the same store in the same week when the customer visits.
Moreover, we follow the procedures in \cite{van2015market} and \cite{ csimcsek2018expectation}: for each purchase record, four no-purchase records of the same assortment are added to the dataset.
Although model-free pricing doesn't depend on the censored demand, the benchmark model-based approaches do.
The goal is to create a reasonable fraction of customers who do not buy any products.
{After data pre-processing, there are in total 18,217 transactions with 2,460 unique sets of assortment price vectors (store/week combinations).}

For the performance evaluation, we resort to a model (estimated from the data) that describes how consumers choose products and calculate the expected revenue under this model.
We fit three models to the data to highlight the model-free or model-insensitive nature of our approach.
More precisely, in the first MNL model we estimate $\{\alpha_\product\}_{\product=1}^{10}$ and $\beta$ in the choice probability
\begin{equation*}
    \frac{\exp(\alpha_\product-\beta \histprice_{\customer\product})}{1+\sum_{k=1}^{10} \exp(\alpha_k-\beta \histprice_{\customer k})}.
\end{equation*}
In the second mixed logit model, we estimate $\left\{w_l,\beta_l, \alpha_{l1},\dots,\alpha_{l,10}\right\}_{l=1}^2$ in the choice probability
\begin{equation*}
    w_1\frac{\exp(\alpha_{1\product}-\beta_1 \histprice_{\customer\product})}{1+\sum_{k=1}^{10} \exp(\alpha_{1k}-\beta_1 \histprice_{ik})}+w_2\frac{\exp(\alpha_{2j}-\beta_2 \histprice_{\customer\product})}{1+\sum_{k=1}^{10} \exp(\alpha_{2k}-\beta_2 \histprice_{ik})}.
\end{equation*}
Both models are estimated using BIOGEME.
We also fit a linear demand model with parameters $\left\{\alpha_\product,\beta_\product\right\}_{\product=1}^{10}$ and $\{\beta_{jk},\gamma_{jk}\}_{j\neq k}$.
The choice probability of product $\product$ in the linear model is
\begin{equation*}
    \alpha_{\product}-\beta_\product \histprice_{\customer\product} + \sum_{k\neq j} (\beta_{jk} I_{ik} \histprice_{ik} + \gamma_{jk} (1-I_{\customer k})),
\end{equation*}
fitted using the ordinary least squares.
Here $I_{\customer k}\in\{0,1\}$ is the indicator for whether product $k$ is included in the assortment seen by customer $i$.
Note that we fit the choice probability separately for product $\product=1,\dots,10$ and the no-purchase probability is one minus their sum.

We generate three datasets for the estimated models considering $\numcustomers=50$ and $\numproducts=10$.
The historical prices $\histprice_{ij}$ are drawn from $\{0.9\bar{\histprice_\product}, 0.95\bar{\histprice_\product}, \bar{\histprice_\product}, 1.05\bar{\histprice_\product}, 1.1\bar{\histprice_\product}\}$, where $\bar{\histprice_\product}$ is the average price of product $\product$ in the IRI data.
The customer choices are then generated using one of the three models.
We calculate the prices using the optimal solution to \eqref{model:OPMIP} and cut-off pricing, and plug them into the three models to evaluate their expected revenues.
We also compute the expected revenue of the incumbent prices, which is the average product prices in the IRI data, under the three models.

Note that we apply our proposed approach to the datasets generated from the estimated models and then evaluate the expected revenues of the model-based prices, the optimal solution to \eqref{model:OPMIP} or cut-off pricing, under the corresponding models.
A seemingly more straightforward approach is to directly apply model-free assortment pricing to the original IRI dataset, and then obtain the optimal prices.
However, the major concern is the conflation of the model misspecification error and the efficacy of the pricing schemes.
Suppose the model-free prices calculated using the IRI data perform poorly under, say, the MNL model.
There might be two reasons: (1) the estimated MNL model accurately captures the pattern in the IRI dataset but the data-driven approach fails to approximate the optimal prices of the MNL model,
or (2) the MNL model does not fit the data and the data-driven approach, which is solely based on the data, cannot possibly approximate the optimal prices of the MNL model.
We lack a reliable way to disentangle the two factors, the goodness of fit versus the performance of our approach.
By the simulated datasets using the three models,
we can control the goodness of fit and isolate the performance of the data-driven pricing schemes.
Nevertheless, since the three models are fitted using the IRI data, they are expected to capture the choice patterns in reality to a large degree.

We also do not compare model-free assortment pricing with the optimal prices under the three models, but only the incumbent prices. 
This is because the optimal prices of the three models are not realistic.
For example, the optimal prices under the MNL model are $\$59.45$ for all products, and the average optimal price under the linear demand model is $\$99.70$. 
However, the price range of the top-ten purchased products in the IRI razors data is $\$3.29$ to $\$7.51$, and the model-free assortment pricing recommends prices  (both optimal and cut-off) between $\$5.68$ to $\$6.64$. 
That is, the estimated demand models cannot extrapolate the demand outside the price range, although they may approximate the demand patterns inside the price range well.
Thus, the resulting optimal prices are not implementable, further suggesting the stability of the model-free approach.


Table~\ref{experiment:razor} shows our results for this experiment with $200$ instances, for each of the three estimated models.
The linear model generates significantly lower expected revenues than the other two, possibly because the estimated demand is lower in the region around the incumbent prices in the linear model.
Compared with the incumbent prices, data-driven assortment pricing significantly improves the expected revenue under all three demand models.
It suggests that our approach offers a robust improvement in this setting.

\ctable[
  notespar,
  caption={Expected revenues of data-driven assortment pricing and the incumbent prices under the three models fitted using the IRI data.}
  \label{experiment:razor},
doinside =\small
]{l c @{\hspace*{6mm}} c @{\hspace*{6mm}} c@{\hspace*{6mm}} c}
{}{
\toprule
  &\multicolumn{3}{c}{Fitted demand model}  & \\\cmidrule(r){2-5}
\text{Pricing method} & \text{MNL} & \text{Mixed logit}  & \text{Linear Demand}   \\
\midrule
Data-driven optimal  &  $1.507$ ($0.007$)  & $1.513$ ($0.008$)   & $1.197$  ($0.018$)      \\

Data-driven cut-off pricing & $1.730$ ($0.007$)   &   $1.731$ (0.007)  & $1.477$ (0.019) \\



Incumbent prices &       $1.464$ (-) &       $1.468$ (-) & $0.944$ (-) \\
\bottomrule}

\section{Practical Considerations and Concluding Remarks}\label{sec:conclusion}
We point out a few practical considerations when our approach is applied to real-world problems. 
The first issue is censoring, i.e., when a customer walks away without purchasing and thus cannot be observed in the data.
Many pricing approaches struggle to handle censoring and completely ignoring data censoring may result in price distortion.
Fortunately, our approach can incorporated data censoring.
Indeed, as mentioned in Remark~\ref{rem:nopurchasecustomers}, customers who do not buy anything can be removed from the dataset without affecting the resulting prices.
As a result, the data-drive prices do not depend on the censored customers.


One implicit assumption we made is that consumers have price sensitivity identical to one, reflected in the quasilinear utilities.
In fact, this assumption can be easily relaxed.
The IC polyhedron can be constructed in the same fashion, as long as an individual customer has the same price sensitivity for all products.
In this case, dividing \eqref{eq:i-valuation} by the same factor results in the same polyhedron and the theoretical results still hold.

Another assumption we make is that the firm does not observe the new incoming customer's information and hence assumes she behaves similar to one of the previously seen customers with equal probability. 
In practice, the firm may potentially know that the customer is a returning customer. Then, one can analytically solve model~\eqref{model:OPMIP} for that customer. If the firm has seen the customer beforehand, it is optimal under the worst-case valuations of the customer to set the prices equal to the ones the customer observed in her historical purchase. Moreover, using consumer features to predict their shopping behavior is a popular practice in modern retailing. For example, an arriving consumer may have a similar background to a segment of past customers. To incorporate consumer features, we may put different weights (as opposed to equal weights) on the IC polyhedra $\valuationset_1, \valuationset_2, \dots, \valuationset_\numcustomers$ in the formulation, based on how similar the arriving consumer is to a past one. Computationally, our approach can still accommodate this case. It remains an exciting research direction to capture the consumer features and properly reflect them in the weights. 

Incorporating product features is another exciting direction.
Consumers form valuations for the base model of a product and certain features. 
The valuation for a product configuration could be obtained as the sum of valuations of the base model and the added features. 
In this setting, the IC polyhedron of valuations can still be formulated 
for customers in the transaction data. 
In contrast to our problem, the pricing should concern the configurations instead of separate base model and features.
Yet it remains an open problem if nonlinear pricing can be analyzed and  efficiently solved.

While we are mainly interested in a situation where the number of customers is low relative to number of products, scalability is still important from a practical point of view. We have demonstrated in Section~\ref{sec:numerics} that our mixed-integer programming reformulations can handle mid-scale problems with hundreds of samples.
For large-scale data, our best approximation (cut-off pricing) is of low complexity and scales linearly with the number of samples. While it provides interpretable and intuitive prices, without the need for a commercial solver, such an approximation also has robust theoretical guarantees and good empirical performance as suggested by our numerical study.


%
%
%
{
\bibliographystyle{informs2014} 
\bibliography{references} 
}

\newpage

\ECSwitch

\ECHead{E-Companion to \textit{Model-Free Assortment Pricing with Transaction Data}}
\section{Additional Tables and Figures}\label{app:tables}
Table~\ref{tab:approx-perf} shows the percentage performance ratio of the conservative, LP relaxation, and cut-off pricing as explained in Section~\ref{sec:approximationPerformance}. All the instances were solved to optimality using Gurobi 9.0.0 in Python with a desktop computer (Intel Core i7-8700, 3.2 GHz). Similarly, Table~\ref{tab:approx-guar} portrays the average performance guarantee of conservative and cut-off pricing approximation algorithms.

\begingroup
\begin{center}

\ctable[
  notespar,
  caption={The relative performance of the approximation strategies. Standard errors are reported in parentheses.}
  \label{tab:approx-perf}, pos=h,
doinside =\small
]{r c @{\hspace*{6mm}} c @{\hspace*{6mm}} c@{\hspace*{6mm}} c}
{{} }{
\toprule
  &\multicolumn{3}{c}{Performance relative to the optimal solution}  & \\\cmidrule(r){2-5}
$(\numcustomers,\numproducts)$ & \text{Conservative} & \text{LP Relaxation}  & \text{Cut-off}   \\
\midrule

$(50,10)$        &  $12.5\%$  ($0.6$\%)  & $79.8\%$   ($0.6$\%)   & $97.6\%$ ($0.1$\%)      \\

$(50,15)$ & $10.8\%$ (0.4\%)  &   $73.3\%$ (0.6\%)  & $97.0\%$ (0.1\%)\\

$(50,20)$ &       $9.1\%$ (0.4\%)&       $70.4\%$ (0.5\%)    & $96.5\%$  (0.1\%) \\

$(50,25)$ &   $9.1\%$ (0.5\%)  &    $72.1\%$ (0.5\%) & $96.0\%$ (0.1\%)  \\

$(100,10)$ & $7.4\%$ (0.3\%)  &   $82.4\%$ (0.4\%)  & $99.0\%$ (0.1\%)\\

$(150,10)$ &       $5.3\%$ (0.2\%) &       $85.0\%$ (0.3\%)    & $99.3\%$ (0.1\%) \\

$(200,10)$ &   $4.4\%$ (0.2\%)  &    $85.4\%$ (0.3\%)  & $99.6\%$  (0.1\%)  \\

\bottomrule}

\ctable[
  notespar,
  caption={The performance guarantee of approximation strategies. Standard errors are reported in parentheses.}
  \label{tab:approx-guar}, pos=h,
doinside =\small
]{r c @{\hspace*{6mm}} c @{\hspace*{6mm}} c}
{{} }{
\toprule
  &\multicolumn{2}{c}{Theoretical performance guarantee}  & \\\cmidrule(r){2-4}
$(\numcustomers,\numproducts)$ & \text{Conservative}  & \text{Cut-off}   \\
\midrule

$(50,10)$ &  $2.2\%$  ($0.2$\%)  & $50.1\%$   ($3.8$\%) \\

$(50,15)$  &  $1.9\%$  ($0.1$\%)  & $49.5\%$   ($3.8$\%) \\

$(50,20)$ &  $1.7\%$  ($0.1$\%)  & $49.7\%$   ($4.0$\%) \\

$(50,25)$  &  $1.9\%$  ($0.2$\%)  & $50.4\%$   ($3.7$\%) \\

$(100,10)$  &  $1.0\%$  ($0.1$\%)  & $49.7\%$   ($2.9$\%) \\

$(150,10)$  &  $0.6\%$  ($0.0$\%)  & $50.1\%$   ($2.5$\%) \\

$(200,10)$  &  $0.5\%$  ($0.0$\%)  & $49.9\%$   ($2.1$\%) \\
\bottomrule}
\end{center}

\ctable[
  notespar,
  caption={The implied performance guarantee for cut-off pricing using the IRI data.}
  \label{tab:ratios}, pos=h,
doinside =\small
]{l c @{\hspace*{6mm}} c@{\hspace*{6mm}} c@{\hspace*{6mm}} r}
{{} }{
\toprule
  &\multicolumn{2}{c}{Ratios}  & \\\cmidrule(r){2-3}
\text{Category} & \text{Median/Mean} & \text{$\ubar{P}/\bar{P}$}  & \text{Implied performance guarantee}   \\
\midrule
Beer & 0.91	& 0.025	& 45.7\% \\

Blades & 0.772	& 0.02 &	38.6\% \\

Carbonated beverages & 0.811 	& 0.026 &	40.6\% \\

Cigarettes & 0.466	& 0.028 &	23.3\% \\

Coffee & 0.98  & 0.033 & 49.0\% \\

Cold cereals & 0.986 &	0.069 &	49.3\% \\

Deodorants & 0.97  &	0.072	& 48.5\% \\

Dippers & 0.769	& 0.07 &	38.5\% \\

Facial tissue & 0.812	& 0.04	& 40.6\%\\

\text{Frozen dinners/entrees} & 0.842	& 0.023	& 42.1\%\\

Frozen pizzas & 0.919 &	0.053 &	46.0\%\\

Household cleaner & 0.964 &	0.033 &	48.2\%\\

Hotdogs & 0.989 & 0.037 &	49.5\%\\

Laundry detergent & 0.914	& 0.031 & 45.7\%\\

Margarine/spreads/butter blends & 0.924 &	0.07	& 46.2\%\\

Mayonnaise & 0.951 &	0.09	& 47.5\%\\

Milk & 1.01  &	0.05 &	50.5\%\\

\text{Mustard \& ketchup} & 0.953  & 0.063 & 47.7\%\\

Paper towels & 0.831 &	0.03 &	41.5\%\\

Peanut butter & 0.861  &	0.046 &	43.1\%\\

Razors & 1.02  &	0.145 &	51.1\%  \\

Photography supplies & 0.92 &	0.057 &	46.0\%\\

Salty snacks & 0.982  &	0.023 &	49.0\%\\

Shampoo & 0.891 & 0.023 & 44.5\%\\

Soup & 0.942 & 0.021 & 47.1\%\\

Spaghetti/Italian sauce & 0.887 & 0.072	 & 44.4\%\\

Sugar substitutes & 0.91  &	0.058 &	45.5\% \\

Toilet tissue & 0.935  &	0.032 & 46.8\%\\

Toothbrush & 0.733 & 0.01	& 36.7\%\\

Toothpaste & 0.831  & 0.013	& 41.6\%\\

Yogurt & 0.711 &	0.037 &	35.5\%\\
\bottomrule}

\begin{figure}
\centering
\caption{The sensitivity of the performance of the data-driven optimal pricing \eqref{model:OPMIP} and cut-off pricing relative to the optimal MNL prices estimated from the data with 100 historical customers. The difference in the revenues is converted to percentage by dividing it by the optimal revenue of the model and then averaged.}
\label{fig:dispersion}
\vskip1em

\scalebox{.8}{
    \begin{tikzpicture}
    \pgfplotsset{compat=newest, ymax=.3}
    \begin{axis}[width=10cm, scale only axis,  xlabel={Standard deviation of historical price distribution}, ylabel={Avg. improvement over the estimated MNL prices}, xmin = .4, xmax=1, yticklabel=\pgfmathparse{100*\tick}\pgfmathprintnumber{\pgfmathresult}\,\%, clip=false,
            yticklabel style={/pgf/number format/.cd,fixed,precision=2}]
    
    \addplot+[mark=none, dotted, thin]
    plot [error bars/.cd, y dir = both, y explicit]
    table[x =sd, y=value3, y error =error3]{disperssion.txt};
    \addlegendentry{Optimal pricing for low utility}
    
        \addplot+[mark=none,thin]
    plot [error bars/.cd, y dir = both, y explicit]
    table[x =sd, y=value4, y error =error4]{disperssion.txt};
    \addlegendentry{Cut-off pricing for low utility}
    
    \addplot+[loosely dotted, mark=none] 
    plot [error bars/.cd, y dir = both, y explicit]
    table[x =sd, y=value, y error =error]{disperssion.txt};
    \addlegendentry{Optimal pricing for high utility}
    
      \addplot+[densely dashed, mark=none, thin] 
    plot [error bars/.cd, y dir = both, y explicit]
    table[x =sd, y=value2, y error =error2]{disperssion.txt};
    \addlegendentry{Cut-off pricing for high utility}
    
    \addplot[ultra thin, dashed, black] coordinates {(.4,0) (1,0)};
    \end{axis}
    \end{tikzpicture}
    }
\end{figure}

\begin{figure}
\centering
\caption{The sensitivity of the performance of the data-driven optimal pricing \eqref{model:OPMIP} and cut-off pricing relative to the optimal MNL prices estimated from the data with 50 historical customers. The data is generated from a mixed logit model.}
\label{fig:dispersionmixd}
\vskip1em
\scalebox{.8}{
    \begin{tikzpicture}
    \pgfplotsset{compat=newest}
    \begin{axis}[width=12cm, scale only axis,  xlabel={Standard deviation of historical price distribution}, ylabel={Avg. revenue in the mixed logit model}, xmin = .4, xmax=1, xtick = {0.4, .5,.6,.7,.8,.9,1}, ymin=0, ymax=3.9]
    
     \addplot+[mark=none, dotted, thin]
    plot [error bars/.cd, y dir = both, y explicit]
    table[x =sd, y=value4, y error =error4]{dispersionmixed.txt};
    \addlegendentry{Optimal pricing for low utility}
    
      \addplot+[mark=none,thin]
    plot [error bars/.cd, y dir = both, y explicit]
    table[x =sd, y=value5, y error =error5]{dispersionmixed.txt};
    \addlegendentry{Cut-off pricing for low utility}
    
       \addplot+[loosely dashed, mark=none, thin]
    plot [error bars/.cd, y dir = both, y explicit]
    table[x =sd, y=value6, y error =error6]{dispersionmixed.txt};
    \addlegendentry{MNL Estimated pricing for low utility}
    
    \addplot+[mark=none, dotted, thin]  
    plot [error bars/.cd, y dir = both, y explicit]
    table[x =sd, y=value, y error =error]{dispersionmixed.txt};
    \addlegendentry{Optimal pricing for high utility}
    
       \addplot+[mark=none, thin]  
    plot [error bars/.cd, y dir = both, y explicit]
    table[x =sd, y=value2, y error =error2]{dispersionmixed.txt};
    \addlegendentry{Cut-off pricing for high utility}
    
     \addplot+[loosely dashed, mark=none, thin]  
    plot [error bars/.cd, y dir = both, y explicit]
    table[x =sd, y=value3, y error =error3]{dispersionmixed.txt};
    \addlegendentry{MNL Estimated pricing for high utility}
    
    \end{axis}
    \end{tikzpicture}
    }
\end{figure}
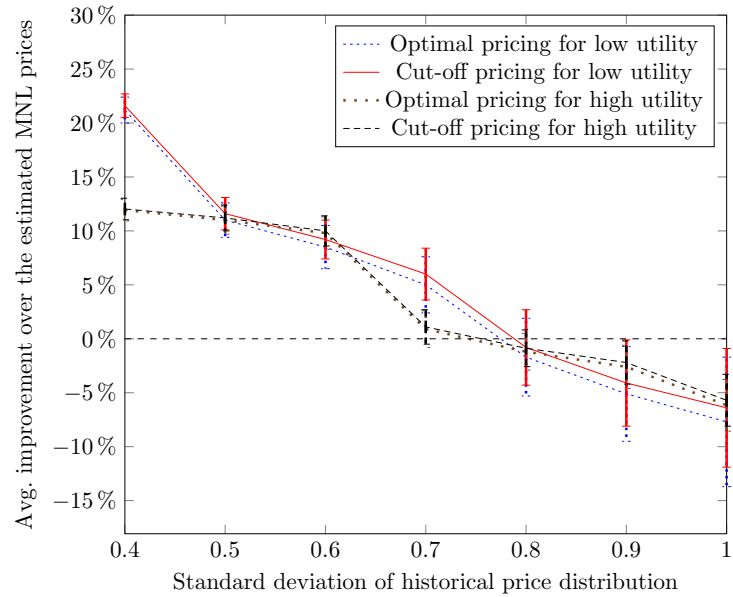

\begin{figure}
\centering
\caption{The sensitivity of the performance of the data-driven cut-off pricing relative to the optimal MNL prices estimated from censored data with 360 historical customers. The difference in the revenues is converted to percentage by dividing it by the optimal revenue of the model and then averaged.}
\label{fig:censor}
\vskip1em
\scalebox{.8}{
    \begin{tikzpicture}
    \pgfplotsset{compat=newest, ymax=.3}
    \begin{axis}[width=10cm, scale only axis,  xlabel={Fraction of non-purchasing customers observed}, ylabel={Avg. improvement over the estimated MNL prices}, xmin = .2, xmax=1, yticklabel=\pgfmathparse{100*\tick}\pgfmathprintnumber{\pgfmathresult}\,\%, clip=false,
            yticklabel style={/pgf/number format/.cd,fixed,precision=2}]
    
            \addplot+[mark=none,thin]
    plot [error bars/.cd, y dir = both, y explicit]
    table[x =sd, y=value2, y error =error2]{censor.txt};
    \addlegendentry{Cut-off pricing for low utility}
    
    \addplot+[densely dashed, mark=none, thin] 
    plot [error bars/.cd, y dir = both, y explicit]
    table[x =sd, y=value, y error =error]{censor.txt};
    \addlegendentry{Cut-off pricing for high utility}
    
    \addplot[ultra thin, dashed, black] coordinates {(.2,0) (1,0)};
    \end{axis}
    \end{tikzpicture}
    }
\end{figure}
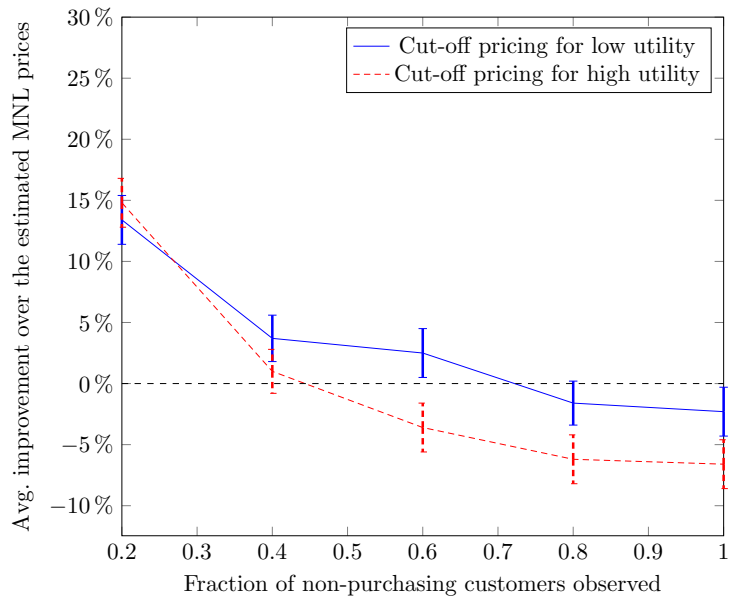

\color{mycolor}
\label{pg:ratiostable}\Copy{ratiostable}{Table~\ref{tab:ratios} contains the calculated performance bound of cut-off pricing using 
the $31$ product categories in the IRI academic data set. 
We consider the ratios of median purchase price to mean purchase price and $\ubar{P}/\bar{P}$ after proper data cleaning. We note that $\ubar{P}/\bar{P}$ could be very small in many categories, but this is primarily due to data entry errors.
For example, some products are purchased at \$0.01.
We remove the top and bottom $0.001$ price quantiles and present the ratio in Table~\ref{tab:ratios}. The implied performance guarantee column in Table~\ref{tab:ratios} contains the implied theoretical performance guarantee by these ratios for each category, obtained through Proposition~\ref{prop:cutoff-pricing-bound}.}

\rev{We also study the sensitivity of the performance of our data-driven model-free optimal and cut-off pricing algorithms with respect to the dispersion of prices in the data. To that avail, we first revisit the low-utility and high-utility experiments with a limited number of customers, explored in Sections~\ref{sec:small-data-num}. For the high-utility experiments,  
in each of the 200 instances, $\{\alpha_\product\}_{\product=1}^{10}$ are independent and drawn uniformly at random from the interval $[1,3]$, while historical prices $\histprice_{\customer\product}$ are drawn uniformly at random from the interval $7\pm \sqrt{3}\sigma$, where $\sigma$ is the standard deviation of historical prices. In the low-utility experiment, in each of the 200 instances, $\{\alpha_\product\}_{\product=1}^{10}$ are independent and drawn uniformly at random from the interval $[-2,0]$ and historical prices are drawn uniformly at random from the interval $3.5 \pm \sqrt{3}\sigma$. }

\rev{Figure~\ref{fig:dispersion} portrays our results for both experiments when the number of historical customers is fixed at $100$. We observe that in both experiments, the performance of the model-free methodologies compared to the estimated MNL optimal prices deteriorate with the increase of $\sigma$. This is intuitive because given that the average historical prices are kept fixed, with the increase of $\sigma$ (while customers make decisions based on the MNL model), more purchases happen at lower prices, making the model-free approaches more conservative. In the meantime, due to a higher price dispersion in the historical data, the MNL estimation which is correctly specified improves, leading to better performing prices.}

\rev{We further study the situation where the MNL model is misspecified and revisit the low-utility and high-utility experiments in Section~\ref{sec:model-misspec}. Here, like before, we set $\beta_1=0.5$ and $\beta_2=2$, while in each of the 200 instances, we randomly draw $(\alpha_{1\product},\alpha_{2\product})$ independently from $[1,3]$ (for high-utility) and $[-2,0]$ (for low-utility). In the high-utility experiments, historical prices $\histprice_{\customer\product}$ are drawn uniformly at random from the interval $7\pm \sqrt{3}\sigma$, while in the low-utility experiment, they are drawn uniformly at random from the interval $3.5 \pm \sqrt{3}\sigma$.}

\rev{Figure~\ref{fig:dispersionmixd} displays our results for both experiments when the number of historical customers is fixed at $50$ and the customers in the historical data made decisions based on the mixed logit model described in Section~\ref{sec:model-misspec}. We observe that like before, with the increase of $\sigma$ the performances of the model-free approach slightly deteriorate, as there will be a larger portion of low price purchases in the historical data which tend to be reflected in the model-free pricing. However, what is interesting is that the performance of the estimated MNL optimal prices does not improve with the increase of $\sigma$, which can be attributed to the fact that the MNL model is misspecified here.}

\rev{Finally, Figure~\ref{fig:censor} shows our results for the case when the data does not record all the non-purchasing customers. In other words, the data is censored and we revisit the low-utility and high-utility experiments with a limited number of customers, explored in Sections~\ref{sec:small-data-num}. As it was mentioned in Remark~\ref{rem:nopurchasecustomers}, censoring does not affect the prices prescribed by our data-driven methodologies, however, it can significantly affect the estimated MNL optimal prices, even when the the MNL choice model is correctly specified. The result in  Figure~\ref{fig:censor} suggests that for a relatively large number of historical customers, $360$, in the high-utility experiment from Section~\ref{sec:model-misspec}, as long as only $40\%$ or less of the non-purchasing customers are recorded in the data, the data-driven prices prescribed by cut-off pricing outperform estimated MNL optimal prices in terms of the expected revenue. We note that in the low-utility experiment from Section~\ref{sec:model-misspec}, this threshold on the fraction of recorded non-purchase customers becomes $60\%$. }

\color{black}
\newpage
\color{mycolor}
\section{Proofs}\label{app:proofs}

\begin{proof}{Proof of Proposition \ref{prop:compare-opt-price}.}
By Proposition~\ref{prop:same-price}, the optimal price output by our framework \eqref{model:OPMIP} satisfies
$\price^*\in \argmax_{\price\ge 0} \price \sum_{i=1}^\numcustomers \indicator(\histprice_{i1}\ge \price)$. Equivalently, we may write $\price^*\in \argmax_{\price\ge 0} \price \frac{\sum_{i=1}^\numcustomers \indicator(\histprice_{i1}\ge \price)}{\numcustomers}$.
By the Law of Large Numbers, we have $\lim_{\numcustomers\to\infty}\frac{\sum_{i=1}^\numcustomers \indicator(\histprice_{i1}\ge \price)}{\numcustomers}=\pr(\histprice_{i1}\ge p) =\int_p^{+\infty} (1- \distfunc(x)) \mathrm{d}\pdist(x)$.
Therefore,
\[ p^* \in \argmax_{\price\ge 0} \lim_{\numcustomers\to\infty} \price \frac{\sum_{i=1}^\numcustomers \indicator(\histprice_{i1}\ge \price)}{\numcustomers} = \argmax_{\price\ge 0}  \price\int_\price^{+\infty} g(x)(1-\distfunc(x))\mathrm{d}x.\]

To prove the last part of Proposition~\ref{prop:compare-opt-price}, we replace $g(x)$ in \eqref{eq:opt-price-single-prod} with an appropriate scaling of $f(x)/(1-\distfunc(x))$ using $\Lambda>0$. Hence, when $g(\price) \propto f(\price)/(1-\distfunc(\price))$,
\[\argmax_{\price\ge 0} \price\int_\price^{+\infty} g(x)(1-\distfunc(x))\mathrm{d}x= \argmax_{\price\ge 0} \Lambda \price \int_\price^{+\infty} f(\price)\mathrm{d}x=\argmax_{\price\ge 0} \revfunc(\price).\] 
\hfill $\square$
\end{proof}

\begin{proof}{Proof of Theorem~\ref{prop:rev-bound}.}
The proof will proceed in several steps. 
\begin{enumerate}
    \item We will show that the solution to \eqref{eq:focsing} on $[0,a]$ with $\price\int_\price^{a} (1-\distfunc(x))\mathrm{d}x>0$ is the unique maximizer of \eqref{eq:opt-price-single-prod}, $p^*$.
    \item We will show that we can assume without any loss of generality that $a=1$.
    \item We will show that if $p^*\leq \hat{p}$ then $\frac{\revfunc(p^*)}{\revfunc(\hat p)}\ge \frac{1}{2}$ and for any $\epsilon>0$ construct an example such that $\frac{\revfunc(p^*)}{\revfunc(\hat p)}\le \frac{1}{2}+\epsilon$.
    \item We will show that if $p^*>\hat{p}$ then $\frac{\revfunc(p^*)}{\revfunc(\hat p)}\ge \frac{\hat p}{p^*}$ and for any $\epsilon>0$ construct an example such that $\frac{\revfunc(p^*)}{\revfunc(\hat p)}\le \frac{\hat{p}^{1-\epsilon}}{p^{*1-\epsilon}}$.
\end{enumerate}

Step 1. When $g(p)=1/b$, $p^*=\argmax_{\price\ge 0} \price\int_\price^{b} g(x)(1-\distfunc(x))\mathrm{d}x=\argmax_{\price\ge 0} \price\int_\price^{a}\frac{1}{b}(1-\distfunc(x))\mathrm{d}x$. Thus
\begin{equation}\label{eq:our-single-price}
p^*=\argmax_{\price\ge 0} \price\int_\price^{a}(1-\distfunc(x))\mathrm{d}x=\argmax_{\price\ge 0} \price \surfunc(p).
\end{equation}
Note that the function $H(p)\triangleq p\surfunc(p)$ on the right hand side of \eqref{eq:our-single-price} is not zero everywhere and is a continuous, bounded and differentiable function, defined over the compact region of $[0,a]$.
On the two ends, it is easy to see $0\surfunc(0)= a\surfunc(a)=0$.
Therefore, the maximizer of $H(p)$ must be in the interior of $[0,a]$ and satisfy the first-order condition for \eqref{eq:our-single-price}:
\[\int_{p^*}^a (1-\distfunc(x))\mathrm{d}x-p^*(1-\distfunc(p^*))= \surfunc(p^*)-\revfunc(p^*)=0.\]
This proves that $p^*$ must satisfy \eqref{eq:focsing}.

Step 2. Assume $a\neq 1$. Then, we could define $y=\frac{p}{a}$, $y\in[0,1]$ and define a new distribution $\bar{\distfunc}(y)=\distfunc(ay)$. We have \[\bar{\revfunc}(y)=y(1-\bar{\distfunc}(y))=\frac{a}{a}y(1-\distfunc(ay))=\frac{1}{a}p(1-\distfunc(p))=\frac{1}{a}\revfunc(p).\]
Moreover, 
\[\bar{\surfunc}(y)=\int_{y}^{1} (1-\bar{\distfunc}(x))\mathrm{d}x=\frac{a}{a}\int_{y}^{1}(1-\distfunc(ax))\mathrm{d}x=\frac{a}{a}\int_{ay}^{a}\frac{1}{a} (1-\distfunc(z))\mathrm{d}z=\frac{1}{a}\int_{p}^{a}(1-\distfunc(z))\mathrm{d}z=\frac{1}{a}\surfunc(p).\]
Where the third equation comes from a change of variable in the integral as $z=ax$. Hence, if $a\neq 1$, we can simply rescale the domain of customer valuations by $a$ and define a new distribution over the new domain. Thus, we can assume $a=1$.

Step 3. Assume $p^*\leq\hat{p}$. Notice that
\begin{equation}\label{eq:welfare}
p(1-\distfunc(p))+\int_{p}^{1} (1-\distfunc(x))\mathrm{d}x
\end{equation}
is decreasing in $p$ since
\[\bigg(p( 1-\distfunc(p))+\int_{p}^{1} (1-\distfunc(x))\mathrm{d}x\bigg)'=1-\distfunc(p)-pf(p)-1+\distfunc(p)=-pf(p)\leq0.\] 
Thus, by assumption of $p^*\leq\hat{p}$, we have
\[R(\hat{p})\leq \hat{p}(1-\distfunc(\hat{p}))+\int_{\hat{p}}^{1} (1-\distfunc(x))\mathrm{d}x\leq p^*(1-\distfunc(p^*))+\int_{p^*}^{1} (1-\distfunc(x))\mathrm{d}x=2p^*(1-\distfunc(p^*))=2R(p^*).\]
The first inequality comes from the fact that $\int_{\hat{p}}^{1} (1-\distfunc(x))\mathrm{d}x\geq 0$. The second inequality follows since \eqref{eq:welfare} is decreasing in $p$ and $p^*\leq\hat{p}$. The prior to the last equality comes from the fact that $p^*$ must satisfy \eqref{eq:focsing}.

To prove that this bound is tight, assume $0<\epsilon_1<\epsilon$ that is sufficiently small and consider the following distribution $F$ such that:
\begin{align*}
F(p) =
    \begin{cases}
        0,                   & \textnormal{if $0\leq p\leq 1-\epsilon_1$},\\
        \frac{1}{\epsilon_1}(p-(1-\epsilon_1)),  & \textnormal{if $1-\epsilon_1< p\leq 1$.}
    \end{cases}
\end{align*}

\label{pg:explicit} \Copy{explicit}{This leads to:

\begin{minipage}{0.45\textwidth}
\begin{align*}
R(p) =
    \begin{cases}
        p,                   & \textnormal{if $0\leq p\leq 1-\epsilon_1$},\\
       \frac{p(1-p)}{\epsilon_1},  & \textnormal{if $1-\epsilon_1< p\leq 1$.}
    \end{cases}
\end{align*}
\end{minipage}
\begin{minipage}{0.45\textwidth}
\begin{align*}
S(p) =
    \begin{cases}
        1-p-\frac{\epsilon_1}{2},                   & \textnormal{if $0\leq p\leq 1-\epsilon_1$},\\
        \frac{(1-p)^2}{2\epsilon_1},  & \textnormal{if $1-\epsilon_1< p\leq 1$.}
    \end{cases}
\end{align*}
\end{minipage}}

\vskip1em
\color{black}
To prove tightness, we need to show that $F$ is a proper CDF that satisfies Assumptions~\ref{asp:cdf-F-G},~\ref{asp:uniform} and~\ref{asp:unique-max}. It is easy to check that $F(p)$ is a proper CDF as it is nonnegative and increasing in $p$ while $F(0)=0$ and $F(1)=1$, and moreover, for any $b\ge1$, Assumptions~\ref{asp:cdf-F-G} and~\ref{asp:uniform} will be satisfied. 

Hence, it remains to show that the bound is tight for $F$ and that $F$ satisfies Assumption~\ref{asp:unique-max}. Notice that for all $0 \leq p\leq 1-\epsilon_1$, $\revfunc'(p)=(p)'=1>0$ and $\surfunc'(p)=-1<0$.
Meanwhile, $\surfunc(0)>\revfunc(0)=0$ and  \[\frac{\epsilon_1}{2}=\surfunc(1-\epsilon_1)<\revfunc(1-\epsilon_1)=1-\epsilon_1\]
when $\epsilon_1<2/3$.
Hence, there is a single point on $[0,1-\epsilon_1]$ that can satisfy \eqref{eq:focsing}.
In fact, this point can be calculated as  $\surfunc(\frac{1}{2}-\frac{\epsilon_1}{4})=\revfunc(\frac{1}{2}-\frac{\epsilon_1}{4})$. 
Furthermore, $\surfunc(1)=\revfunc(1)=0$ and for all $1-\epsilon_1\leq p \leq 1$, if $\epsilon_1<{1}/{3},$ \[\revfunc'(p)=(p-\frac{1}{\epsilon_1}(p^2-(p-p\epsilon_1)))'=1-\frac{2p-1+\epsilon_1}{\epsilon_1}<\surfunc'(p)=\frac{p-1+\epsilon_1}{\epsilon_1}-1\leq 0.\]
Hence $\surfunc(p)$ and $\revfunc(p)$ cannot cross at any interior point on $[1-\epsilon_1, 1]$ when $\epsilon_1<{1}/{3}$. Thus the only candidate for $p^*$ if $\epsilon_1$ is sufficiently small is $\frac{1}{2}-\frac{\epsilon_1}{4}$ (and clearly  $(\frac{1}{2}-\frac{\epsilon_1}{4})\int_{(\frac{1}{2}-\frac{\epsilon_1}{4})}^{1} \big(1-F(z)\big)\mathrm{d}z>0$) and it is the unique maximizer of $p\int_{p}^{1} (1-F(z))\mathrm{d}z$ (we showed this in Step 1.). 

Moreover, $\revfunc(\frac{1}{2}-\frac{\epsilon_1}{4})=\frac{1}{2}-\frac{\epsilon_1}{4}$ while if $\epsilon_1\leq \frac{1}{3}$, $\hat{p}=1-\epsilon_1$ and $\revfunc(\hat{p})=1-\epsilon_1$, resulting in $\revfunc(p^*)/\revfunc(\hat p)=\frac{\frac{1}{2}-\frac{\epsilon_1}{4}}{1-\epsilon_1}\leq \frac{1}{2}+\epsilon$ for $\epsilon_1$ sufficiently small, thus proving the bound is tight (while $\revfunc(p)$ is unimodal and has a unique maximizer $\hat p$). Hence, we have shown the bound is tight for $F$ and that $F$ satisfies Assumption~\ref{asp:unique-max}.

Step 4. Assume $p^*>\hat{p}$. Notice that by assumption, $p^*$ is the first and only point on $[0,1]$ that satisfies \eqref{eq:focsing} with $p^*\int_{p^*}^{1} (1-\distfunc(x))\mathrm{d}x>0$. Moreover since $p^*>\hat{p}$, \[\int_{\hat{p}}^{1} (1-F(z))\mathrm{d}z\ge \int_{p^*}^{1} (1-F(z))\mathrm{d}z,\] while $\hat p>0$. Hence $\hat{p}\int_{\hat{p}}^{1} (1-\distfunc(x))\mathrm{d}x>0$ and by assumption we have $\revfunc(\hat p)\neq \surfunc(\hat p)$. Therefore, since \[\int_{0}^{1}(1-F(z))\mathrm{d}z>0(1-F(0)),\]
then we have that 
\[\revfunc(\hat{p})< \int_{\hat{p}}^{1} (1-F(z))\mathrm{d}z= \frac{\hat{p}}{\hat{p}}\int_{\hat{p}}^{1} (1-F(z))\mathrm{d}z < \frac{p^*}{\hat{p}}\int_{p^*}^{1} (1-F(z))\mathrm{d}z=\frac{p^*}{\hat{p}}\revfunc(p^*),\]
where the last inequality follows from Assumption~\ref{asp:unique-max} and that $p^*$ is the unique maximizer of \eqref{eq:opt-price-single-prod} and the last equality comes from the fact that $p^*$ must satisfy \eqref{eq:focsing}. Thus, \[\frac{\revfunc(p^*)}{\revfunc(\hat{p})}\geq \frac{\hat{p}}{p^*}. \]

It remains to be shown that the bound is tight when $p^*>\hat{p}$. Consider any $0<x_1<x_2<1/e$. Notice that by Assumption~\ref{asp:unique-max}, $\revfunc(p)$ is unimodal on $[0,1]$. Moreover, by assumption we have $p^*>\hat{p}$, hence, for all $x>p^*,\ x(1-F(x))\leq p^*(1-F(p^*))$. Therefore, for all $x>p^*,\ 1-F(x)\leq \frac{\revfunc(p^*)}{x}$. Thus, we claim $p^*<\frac{1}{e}$. This is because, \[\revfunc(p^*)=\int_{p^*}^{1}(1-F(z))\mathrm{d}z<\int_{p^*}^{1} \frac{\revfunc(p^*)}{z}\mathrm{d}z=\revfunc(p^*)\log(\frac{1}{p^*}),\] 
where the first equality comes from the fact that $p^*$ must satisfy \eqref{eq:focsing} while the strict inequality follows form the fact that $F$ is assumed to be continuous and $F(1)=0$, proving that $p^*< \frac{1}{e}$. Let $\bar{x}=ex_2$, we know $\bar{x}<1$.
 Assume $0<\epsilon_2<\epsilon$ is sufficiently small and consider the following distribution $F$ such that:
\begin{align*}
1-F(p) =
    \begin{cases}
        1-\frac{1-\frac{\delta}{x_1^{2-\epsilon_2}}}{\epsilon_2}p,                   & \textnormal{if $0\leq p\leq \epsilon_2$},\\
       \frac{\delta}{x_1^{2-\epsilon_2}},  & \textnormal{if $\epsilon_2\leq p\leq x_1$},\\
        \frac{\delta}{p^{2-\epsilon_2}},  & \textnormal{if $x_1\leq p\leq x_2$},\\
         \frac{\delta}{px_2^{1-\epsilon_2}},  & \textnormal{if $x_2\leq p\leq \bar{x}$},\\
         \frac{\delta}{\bar{x}x_2^{1-\epsilon_2}}-\frac{ \frac{\delta}{\bar{x}x_2^{1-\epsilon_2}}}{\epsilon_2}(p-\bar{x}),  & \textnormal{if $\bar{x}\leq p\leq \bar{x}+\epsilon_2$},\\
         0,                   & \textnormal{if $\bar{x}+\epsilon_2\leq p\leq 1$.}
    \end{cases}
\end{align*}

\label{pg:construction} \Copy{construction}{\rev{The intuition behind the construction of $F$ is as follows. We want to ensure that $\hat{p}=x_1$ while $p^*$ can be arbitrarily close to $x_2$. If $F$ is such that at these points $R(\cdot)$ is proportional to $p$, we would be very close to our goal. It would remain to ensure $p^*$ happens at the unique point at which $R(\cdot)$ and $S(\cdot)$ cross each other. Thus, $F$ must be such that from $\hat{p}$ to $x_2$, the difference between $R(\cdot)$ and $S(\cdot)$ shrinks and reaches zero at a point slightly larger than $x_2$.}}

To formally prove tightness, we need to show that $F$ is a proper CDF that satisfies Assumptions~\ref{asp:cdf-F-G},~\ref{asp:uniform} and~\ref{asp:unique-max}. It is easy to check that $F(p)$ is nonnegative and increasing in $p$ while $F(0)=0$ and $F(1)=1$.
Moreover, for any $b\ge1$, Assumptions~\ref{asp:cdf-F-G} and~\ref{asp:uniform} are satisfied. 

Hence, it remains to show that the bound is tight for $F$ and it satisfies Assumption~\ref{asp:unique-max}. Notice that for all $0\leq p\leq \epsilon_2$, $\revfunc(p)=p-\frac{1-\frac{\delta}{x_1^{2-\epsilon_2}}}{\epsilon_2}p^2$ and $\revfunc'(p)=1-2\frac{1-\frac{\delta}{x_1^{2-\epsilon_2}}}{\epsilon_2}p>0$ if $\delta>\frac{1}{2}x_1^{2-\epsilon_2}$. Moreover, for all $\epsilon_2\leq p\leq x_1$ we have $\revfunc'(p)=\frac{\delta}{x_1^{2-\epsilon_2}}>0$ while $\surfunc(p)$ is decreasing for all $0\leq p \leq x_1$. Furthermore, notice that \[\surfunc(x_2)> \int_{x_2}^{\bar{x}} \frac{\delta}{p 
x_2^{1-\epsilon_2}}\mathrm{d}p=\frac{\delta}{x_2^{1-\epsilon_2}} 
\log \left(\frac{\bar{x}}{x_2} \right)=\frac{\delta}{x_2^{1-\epsilon_2}}=\revfunc(x_2),\] while \begin{eqnarray*}\surfunc(x_1)&=&\int_{x_1}^{x_2}\frac{\delta}{p^{2-\epsilon_2}}\mathrm{d}p +\surfunc(x_2)>\frac{\delta}{(1-\epsilon_2)}(\frac{1}{x_1^{1-\epsilon_2}}-\frac{1}{x_2^{1-\epsilon_2}})+\frac{\delta}{x_2^{1-\epsilon_2}}\\&=&\frac{\delta}{(1-\epsilon_2)}\frac{1}{x_1^{1-\epsilon_2}}-\frac{\epsilon_2 \delta}{(1-\epsilon_2)}\frac{1}{x_2^{1-\epsilon_2}}>\frac{\delta}{x_1^{1-\epsilon_2}}=\revfunc(x_1),\end{eqnarray*} 
thus $\surfunc(p)$ and $\revfunc(p)$ do not cross on $0\leq p\leq x_1$. This is because $\revfunc(p)$ is strictly increasing and $\surfunc(p)$ is decreasing in this region while $\surfunc(0)>\revfunc(0)$ and $\surfunc(x_1)>\revfunc(x_1)$. 

Moreover, if $x_1\leq p\leq x_2$, $\revfunc(p)=\frac{\delta}{p^{1-\epsilon_2}}$ and \[\revfunc'(p)=-(1-\epsilon_2)\frac{\delta}{p^{2-\epsilon_2}}>-\frac{\delta}{p^{2-\epsilon_2}}=\surfunc'(p),\] while $\surfunc(x_2)>\revfunc(x_2)=\frac{\delta}{x_2^{1-\epsilon_2}}$ thus $\surfunc(p)$ and $\revfunc(p)$ do not cross on $x_1\leq p\leq x_2$ either. This is because if they cross at any point in this regions, since $\revfunc'(p)>\surfunc'(p)$ in this region, it must be that $\surfunc(x_2)\leq \revfunc(x_2)$ which is a contradiction since we showed earlier that $\surfunc(x_2)>\revfunc(x_2)$. 

When $x_2\leq p\leq \bar{x}$, $\revfunc'(p)=0$ while $\surfunc'(p)<0$, thus $\surfunc(p)$ and $\revfunc(p)$ can cross at most once in this region, as $\surfunc(x_2)>\frac{\delta}{x_2^{1-\epsilon_2}}=\revfunc(x_2)$ while $\surfunc(\bar{x})<\revfunc(\bar{x})=\frac{\delta}{x_2^{1-\epsilon_2}}$ for $\epsilon_2$ sufficiently small. 

Finally, if $\bar{x}\leq p\leq \bar{x}+\epsilon_2$ and $\epsilon_2<\frac{\bar{x}}{2}$, \[\revfunc'(p)=\bigg(\frac{p\delta}{\bar{x}x_2^{1-\epsilon_2}}-\frac{ \frac{\delta}{\bar{x}x_2^{1-\epsilon_2}}}{\epsilon_2}(p^2-p\bar{x})\bigg)'= \frac{\delta}{\bar{x}x_2^{1-\epsilon_2}}-\frac{ \frac{\delta}{\bar{x}x_2^{1-\epsilon_2}}}{\epsilon_2}(2p-\bar{x})<\surfunc'(p)= -\frac{\delta}{\bar{x}x_2^{1-\epsilon_2}}+\frac{\frac{\delta}{\bar{x}x_2^{1-\epsilon_2}}}{\epsilon_2}(p-\bar{x}),\]
while $\surfunc(\bar{x}+\epsilon_2)=\revfunc(\bar{x}+\epsilon_2)=0$. Therefore, $\surfunc(p)$ and $\revfunc(p)$ do not cross on $[\bar{x},\bar{x}+\epsilon_2)$ since otherwise, because $\surfunc'(p)>\revfunc'(p)$ in this region when $\epsilon_2<\frac{\bar{x}}{2}$, it would imply $0=\surfunc(\bar{x}+\epsilon_2)>\revfunc(\bar{x}+\epsilon_2)$ which is a contradiction.

Thus, $\surfunc(p)$ and $\revfunc(p)$ will cross once and that will happen on $[x_2,\bar{x}]$. This crossing point $p^*$ will satisfy $0<p^*-x_2<\epsilon$ for some $\epsilon_2$ sufficiently small. This is because $\revfunc(x_2)=\frac{\delta}{x_2^{1-\epsilon_2}}$ while $\surfunc(x_2)$ will get arbitrarily close to $\frac{\delta}{x_2^{1-\epsilon_2}}$ for a sufficiently small $\epsilon_2$, and $\surfunc'(p)\leq -\frac{\delta}{\bar{x}x_2^{1-\epsilon_2}}$ in this region while $\revfunc(p)$ is constant (while clearly we have $p^*\int_{p^*}^{1} (1-\distfunc(x))\mathrm{d}x>0$).

Therefore we can conclude that the unique maximizer of $\max_{p\geq0}p\int_{p}^{1}(1-F(z))\mathrm{d}z$, $p^*$ satisfies $0<p^*-x_2<\epsilon$ for some $\epsilon_2$ sufficiently small. Moreover, it is easy to see that $\revfunc(p)$ is unimodal and has a unique maximizer at $\hat{p}=x_1$. We have
\[\frac{\revfunc(p^*)}{\revfunc(\hat{p})}=\frac{\frac{\delta}{p^{*1-\epsilon}}}{\frac{\delta}{\hat{p}^{1-\epsilon}}}=\frac{\hat{p}^{1-\epsilon_2}}{p^{*1-\epsilon_2}}\leq \frac{\hat{p}^{1-\epsilon}}{p^{*1-\epsilon}} .\] 

Hence, we have shown the bound is tight for $F$ and it satisfies Assumption~\ref{asp:unique-max}. Finally, since we already showed $p^*<\frac{1}{e}$, we can also see that $\frac{\revfunc(p^*)}{\revfunc(\hat{p})}> e\hat{p}$.

\hfill $\square$
\end{proof}

\begin{proof}{Proof of Lemma~\ref{lem:IGFR}.}
First we will prove that \eqref{eq:opt-price-single-prod} has a  maximizer that is the unique solution to \eqref{eq:focsing} on $[0,a)$. Since $p\int_{p}^{a} (1-F(x))\mathrm{d}x$ is a continuous, bounded and differentiable function, not zero everywhere and defined on a compact set, while \[0\int_{0}^{a} (1-F(x))\mathrm{d}x=a\int_{a}^{a} (1-F(x))\mathrm{d}x=0,\] it must have a maximizer in the interior of $[0,a]$ and this maximizer must satisfy the first-order condition for \eqref{eq:opt-price-single-prod}. Hence, we must have

\begin{equation}\label{eq:crossing}
\bigg(p\int_{p}^{a} (1-F(x))\mathrm{d}x\bigg)'=\int_{p}^{a} (1-F(x))\mathrm{d}x-p(1-F(p))=0.
\end{equation}
Notice that  
\[\int_{0}^{a} (1-F(x))\mathrm{d}x-0(1-F(0))>0,\] 
while if $p\geq \frac{a}{2}$
\[\int_{p}^{a} (1-F(x))\mathrm{d}x-p(1-F(p))\leq (1-F(p))(a-p)-p(1-F(p))\leq (1-F(p))(a-2p)\leq 0.\]
Thus \eqref{eq:crossing} must have at least one root in $(0,\frac{a}{2}]$. We will show that it cannot have more than one root in $[0,a)$. Notice that 
\[\int_{a}^{a} (1-F(x))\mathrm{d}x-a(1-F(a))=0\]
while
\[\bigg(\int_{p}^{a} (1-F(x))\mathrm{d}x-p(1-F(p))\bigg)'=2F(p)-2+pf(p).\]

By assumption, $\frac{pf(p)}{1-F(p)}$ is strictly increasing, thus, $\frac{pf(p)}{1-F(p)}=2$ can only have at most one root. Let $\hat{x}$ be the single root of $\frac{pf(p)}{1-F(p)}=2$. Then, \eqref{eq:crossing} is strictly decreasing for all $p<\hat{x}$ and strictly increasing for all $p>\hat{x}$ and $0\leq \hat{x}<a$ (otherwise we will have $\int_{a}^{a} (1-F(x))\mathrm{d}x-a(1-F(a))<0$ which is a contradiction). Then by noticing that $\int_{0}^{a} (1-F(x))\mathrm{d}x-0(1-F(0))>0$ while $\int_{a}^{a} (1-F(x))\mathrm{d}x-a(1-F(a))=0$, it turns out that \eqref{eq:crossing} must have exactly one root in $[0,a)$. Hence, $p\int_{p}^{a} (1-F(x))\mathrm{d}x$ has a unique point on $[0,a)$ that satisfies \eqref{eq:focsing} and is its maximizer. If we denote this point with $p^*$, it is clear that $p^*\int_{p^*}^{a} (1-F(x))\mathrm{d}x>0$.

Now we prove $\revfunc(p)$ is unimodal in $[0,a]$ and has a unique maximizer. $\revfunc(p)$ is a continuous, bounded and differentiable function defined on a compact set and not zero everywhere while $\revfunc(0)=\revfunc(a)=0$. Thus, $\revfunc(p)$ must have a maximizer in the interior of $[0,a]$ and it must satisfy the first-order condition for $\revfunc(p)$. Thus at optimality we must have
\[\bigg(p(1-F(p))\bigg)'=1-F(p)-pf(p)=0,\] 
and since $\frac{pf(p)}{1-F(p)}$ is assumed to be strictly increasing, $\frac{pf(p)}{1-F(p)}=1$ (equivalent to the first-order condition) can have at most one root. Thus there is a unique point on $[0,a]$ that satisfies the first-order condition for $\revfunc(p)$ and is its unique maximizer, proving that $\revfunc(p)$ is unimodal as well. 

\hfill $\square$
\end{proof}

\begin{proof}{Proof of Lemma~\ref{lem:fixedbound}.}
First, notice that $\revfunc(p)$ is a continuous, bounded and differentiable function defined on a compact set and not zero everywhere while $\revfunc(0)=\revfunc(a)=0$. Thus, $\revfunc(p)$ must have a maximizer in the interior of $[0,a]$ and it must satisfy the first-order condition for $\revfunc(p)$. Thus by the definition of $\hat p$, it must be that $\revfunc'(\hat p)=1-\distfunc(\hat p)-\hat p f(\hat p)=0$, suggesting that $\frac{1}{\hat p}=\frac{f(\hat p)}{1-\distfunc(\hat p)}$.

Now, for the sake of contradiction, let us assume that $\hat{p}<p^*$. Notice that by assumption, $p^*$ is the first and only point on $[0,a]$ that satisfies \eqref{eq:focsing} with $p^*\int_{p^*}^{a} (1-\distfunc(x))\mathrm{d}x>0$. Moreover since $p^*>\hat{p}$, \[\int_{\hat{p}}^{a} (1-F(z))\mathrm{d}z\ge \int_{p^*}^{a} (1-F(z))\mathrm{d}z,\] while $\hat p>0$. Hence $\hat{p}\int_{\hat{p}}^{a} (1-\distfunc(x))\mathrm{d}x>0$ and by assumption we have $\revfunc(\hat p)\neq \surfunc(\hat p)$, and since $\int_{0}^{a}(1-F(z))\mathrm{d}z>0(1-F(0))$, then we have that 
\[\revfunc(\hat{p})= \hat{p} (1-\distfunc(\hat p))< \int_{\hat{p}}^{a} (1-F(z))\mathrm{d}z,\]
suggesting that 
\[\frac{1-\distfunc(\hat p)}{\int_{\hat{p}}^{a} (1-F(z))\mathrm{d}z}<\frac{1}{\hat p}=\frac{f(\hat p)}{1-\distfunc(\hat p)}.\]
However, this is a contradiction since by assumption we have $\frac{h(p)}{1-H(p)}\ge \frac{f(p)}{1-\distfunc(p)}$ for all $p\in [0,a]$, which suggests:
\[\frac{1-\distfunc(\hat p)}{\int_{\hat{p}}^{a} (1-F(z))\mathrm{d}z}=\frac{\frac{1-\distfunc(\hat p)}{\int_{0}^{a} (1-F(z))\mathrm{d}z}}{\frac{\int_{\hat{p}}^{a} (1-F(z))\mathrm{d}z}{\int_{0}^{a} (1-F(z))\mathrm{d}z}}=\frac{h(\hat p)}{1-H(\hat p)}\ge \frac{f(\hat p)}{1-\distfunc(\hat p)}.\]
Finally, given that we have $\hat p \ge p^*$, it follows directly from Theorem~\ref{prop:rev-bound} that $\frac{\revfunc(p^*)}{\revfunc(\hat p)}\ge \frac{1}{2}$.
\hfill $\square$
\end{proof}

\smallskip


\begin{proof}{Proof of Proposition \ref{prop:convrate}.}

Let $\hat{R}(p)\equiv p \int_\price^{b} g(x)(1-\distfunc(x))\mathrm{d}x$ and let $\hat{R}_m(\price)= \price \frac{\sum_{i=1}^\numcustomers \indicator(\histprice_{i1}\ge \price)}{\numcustomers}=p(1-\hat{\digamma}_m(p))$ where $\hat{\digamma}_m(p)=\frac{\sum_{i=1}^\numcustomers \indicator(\histprice_{i1}\ge \price)}{\numcustomers}$ is the empirical CDF corresponding to $\digamma(p)=\int_0^{\price} g(x)(1-\distfunc(x))\mathrm{d}x$. We note that this is rigorous as $\histprice_{i1}$ for all $i\in \customerset$ is generated independently from the distribution $\pr(\histprice_{i1}\ge p) =\int_p^{+\infty} (1- \distfunc(x)) \mathrm{d}\pdist(x)$.

By applying Hoeffding's inequality and noticing that for all $i\in \customerset$ and $p^*$ we have $ 0\le \frac{p^* \indicator(\histprice_{i1}\ge p^*)}{\numcustomers}\le \frac{b}{m}$, we can write:
\begin{equation}\label{eq:ec-bound}
   \pr(|\hat{R}_m(p^*)-\hat{R}(p^*)|\ge t)\le 2e^{\frac{-2t^2 m}{b^2}}\ \forall \ t>0. 
\end{equation}
Define the event $E\triangleq \{|\hat{R}_m(p^*)-\hat{R}(p^*)|< t\}$.
By assumption, for all $p\in [0,b]$, we have: \[p \int_\price^{b} g(x)(1-\distfunc(x))\mathrm{d}x\le p^* \int_{p^*}^{b} g(x)(1-\distfunc(x))\mathrm{d}x-\alpha (p-p^*)^2.\] 
Therefore, for any $p$ such that $|p-p^*|\ge \epsilon$, we have: 
\[\hat{R} (p)\le \hat{R} (p^*)-\alpha \epsilon^2.\]
Consider any $t\in [0, \alpha \epsilon^2]$. 
On event $E$, if $|p^*_m-p^*|\ge \epsilon$, then $\hat{R}(p^*_m)\le \hat{R}(p^*)-\alpha\epsilon^2$ and we have
\[\hat{R}_m(p^*_m)\ge \hat{R}_m(p^*)\ge \hat{R}(p^*)-t\ge \hat{R}(p^*_m)+\alpha\epsilon^2-t.\]
Where the first inequality follows from the definition of $p^*_m$ while the second inequality happens due to event $E$. Therefore, $|\hat{R}_m(p^*_m)-\hat{R}(p^*_m)|\ge\hat{R}_m(p^*_m)-\hat{R}(p^*_m)\ge \alpha\epsilon^2-t$.

In other words, we know that $E\cap \{|p^*_m-p^*|\ge \epsilon\}\subseteq E\cap \{|\hat{R}_m(p^*_m)-\hat{R}(p^*_m)|\ge \alpha\epsilon^2-t\}$.

Moreover, since $|\hat{R}_m(p^*_m)-\hat{R}(p^*_m)|\ge \alpha\epsilon^2-t$ implies that $\max_{p\ge0}(|\hat{R}_m(p)-\hat{R}(p)|)\ge \alpha\epsilon^2-t$, we have
\begin{equation}\label{eq:eevent-step1}
\pr\big(E\cap \{|p^*_m-p^*|\ge \epsilon\}\big)\le \pr(E\cap\{\max_{p\ge0}(|\hat{R}_m(p)-\hat{R}(p)|)\ge \alpha\epsilon^2-t \}\big)\le \pr(\max_{p\ge0}(|\hat{R}_m(p)-\hat{R}(p)|)\ge \alpha\epsilon^2-t).
\end{equation}
Next, we provide an upper bound for 
$\pr(\max_{p\ge0}(|\hat{R}_m(p)-\hat{R}(p)|)\ge \alpha\epsilon^2-t)$. This can be achieved by applying Dvoretzky–Kiefer–Wolfowitz inequality \citep[page 268]{van2000asymptotic}: since $\hat{\digamma}_m(p)$ is the empirical CDF of $\digamma(p)$ for $p\in [0,b]$, we have:
\begin{align}
   \pr(\max_{p\ge0}(|\hat{R}_m(p)-\hat{R}(p)|)&\ge \alpha\epsilon^2-t)\le \pr\left(b \left\|\frac{\sum_{i=1}^\numcustomers \indicator(\histprice_{i1}\ge \price)}{\numcustomers}-\int_\price^{b} g(x)(1-\distfunc(x))\mathrm{d}x\right\|_{\infty}\ge \alpha\epsilon^2-t\right)\notag\\
   &\le 2e^{\frac{-2(\alpha\epsilon^2-t)^2 m}{b^2}}.  \label{eq:eevent-step2}
\end{align}
Thus we have:
\begin{align*}
   \pr(|p^*_m-p^*|\ge \epsilon)&= \pr(E\cap \{|p^*_m-p^*|\ge \epsilon\})+\pr(E^c\cap\{|p^*_m-p^*|\ge \epsilon\})\\
   &\le \pr(E\cap \{|p^*_m-p^*|\ge \epsilon\})+\pr(E^c)\\
   &\le 2e^{\frac{-2(\alpha\epsilon^2-t)^2 m}{b^2}}+2e^{\frac{-2t^2 m}{b^2}},
\end{align*}
where the last inequality follows from \eqref{eq:ec-bound}, \eqref{eq:eevent-step1}, and~\eqref{eq:eevent-step2}.

To finalize the proof, we need to choose $0\le t\le \alpha\epsilon^2$ such that the quasi-convex function (when $m$ is large enough) \[2e^{\frac{-2t^2 m}{b^2}}+ 2e^{\frac{-2(\alpha\epsilon^2-t)^2 m}{b^2}}\] is minimized, which is achieved, when $t=\frac{\alpha \epsilon^2}{2}$. Therefore, we have:

\[\pr(|p^*_m-p^*|\ge \epsilon)\le 4e^{\frac{-\alpha^2\epsilon^4 m}{2b^2}}.\]

Hence, if we want $|p_m^*-p^*|\ge \epsilon$ with probability at most $\delta$, we need up to $m=\big\lceil(\frac{2b^2}{\alpha^2\epsilon^4})\log(\frac{4}{\delta})\big\rceil$ historical customers, suggesting a sample complexity of $O\big((\frac{b^2}{\alpha^2\epsilon^4})\log(\frac{1}{\delta})\big)$.

\hfill $\square$
\end{proof}

\color{black}
\smallskip

\begin{proof}{Proof of Proposition \ref{prop:revenueLP}.}
    The linear program that results from removing the disjunctive constraint \eqref{eq:disjConditions} from \eqref{model:DP}, i.e.,
    \begin{align*}
        \min_{\revenue, \vecvaluation_\customer} \{ \revenue \colon \revenue \ge 0, \vecvaluation_\customer \in \valuationset_\customer \},
    \end{align*}
    has the trivial finite optimum $\revenue = 0$. Thus, the model \eqref{model:DPLP} is obtained by applying directly Corollary 2.1.2 by \cite{balas1998disjunctive}.
    \hfill $\square$
\end{proof}

\begin{proof}{Proof of Lemma \ref{lem:feasibleOptions}.}

We begin with statement \ref{it:feasibleNopurchase}. If $\price_{\choice{\customer}} \ge \histprice_{\customer \choice{\customer}}$, then the valuation $\vecvaluation^\nopurchase_\customer$ defined by $\novalVar_{\customer \choice{\customer}} = \histprice_{\customer \choice{\customer}}$ and $\novalVar_{\customer \product} = 0$ for all $\product \in \productset \setminus \{ \choice{\customer} \}$ belongs to $\valoptionset^{\nopurchase}_{\customer}(\vecprice)$. Conversely,
if $\vecvaluation^\nopurchase_\customer \in \valoptionset^{\nopurchase}_{\customer}(\vecprice)$,
the inequalities $\novalVar_{\customer \product} \le \price_\product$ and $\novalVar_{\customer \product} \ge \histprice_{\customer \choice{\customer}}$ from \eqref{eq:i-valuation} together imply $\price_{\choice{\customer}} \ge \histprice_{\customer \choice{\customer}}$.

For statement \ref{it:feasibleChoice}, the valuation $\vecvaluation^{\choice{\customer}}_\customer$ defined by $\valuation^{\choice{\customer}}_{\customer \choice{\customer}} = \max\{ \histprice_{\customer \choice{\customer}}, \price_{\choice{\customer}}\}$ and $\valuation^{\choice{\customer}}_{\customer \product} = 0$ for all $\product \in \productset \setminus \{ \choice{\customer} \}$ belongs to $\valoptionset^{\choice{\customer}}_{\customer}(\vecprice)$.

Finally, for statement \ref{it:feasibleH}, consider the set of constraints that are satisfied by  points in $\valoptionset^{\product}_{\customer}(\vecprice)$ after re-arranging the constant terms to the right-hand side of the inequalities:
\begin{align}
\valVar_{\customer \product}
&\ge
\price_\product,
    \label{eq:lemfeasProof-1} \\
\valVar_{\customer \product} - \valVar_{\customer \product'}
&\ge
\price_\product - \price_{\product'},
    &&\forall \product' \in \productset, \label{eq:lemfeasProof-2} \\
\valVar_{\customer \choice{\customer}}
&\ge
\histprice_{\customer \choice{\customer}}, \label{eq:lemfeasProof-3} \\
\valVar_{\customer \choice{\customer}} - \valVar_{\customer \product'}
&\ge
\histprice_{\customer \choice{\customer}} - \histprice_{\customer \product'},
    &&\forall \product' \in \productset \setminus \{ \choice{\customer} \}. \label{eq:lemfeasProof-4}
\end{align}
If $\price_{\product} - \price_{\choice{\customer}} \le \histprice_{\customer \product} - \histprice_{\customer \choice{\customer}}$, we construct a valuation
$\vecvaluation^{\product}_\customer \in \valoptionset^{\product}_{\customer}(\vecprice)$
where
$\valuation^{\product}_{\customer \product}
=
\max\{ \histprice_{\customer \product}, \price_\product \}$,
$\valuation^{\product}_{\customer \choice{\customer}}
=
\max\{ \histprice_{\customer \choice{\customer}}, \histprice_{\customer \choice{\customer}} - \histprice_{\customer \product} + \price_\product \}$,
and
$\valuation^{\product}_{\customer \product'} = 0$ for all $\product' \in \productset \setminus \{ \choice{\customer}, \product \}$. In particular, \eqref{eq:lemfeasProof-1} and \eqref{eq:lemfeasProof-3} are satisfied by construction. Assume now $\price_\product > \histprice_{\customer \product}$. For \eqref{eq:lemfeasProof-2} and \eqref{eq:lemfeasProof-4} with $\product' = \choice{\customer}$, we have
\begin{align*}
    \valVar_{\customer \product} - \valVar_{\customer \choice{\customer}}
    =
    \price_\product - \histprice_{\customer \choice{\customer}} + \histprice_{\customer \product} - \price_\product
    =
    \histprice_{\customer \product} - \histprice_{\customer \choice{\customer}}
    \ge
    \price_\product - \price_{\choice{\customer}},
\end{align*}
where the last inequality follows from the statement hypothesis. For $\product' \neq \choice{\customer}$, note that
$\valVar_{\customer \product} - \price_\product$ is zero while
$\valVar_{\customer \product'} - \price_{\product'}$ is non-positive in \eqref{eq:lemfeasProof-2}, and analogously $\valVar_{\customer \choice{\customer}} - \histprice_{\customer \choice{\customer}}$ is positive while $\valVar_{\customer \product'} - \histprice_{\customer \product'} = -\histprice_{\customer \product'}$ is non-positive. If
$\price_\product \le \histprice_{\customer \product}$, note that $\valVar_{\customer \product} = \histprice_{\customer \product}$ and $\valVar_{\customer \choice{\customer}} = \histprice_{\customer \choice{\customer}}$, and the same derivations above apply.

Finally, the sufficient conditions of \ref{it:feasibleH} follow directly from \eqref{eq:lemfeasProof-2} and \eqref{eq:lemfeasProof-4} with $\product' = \choice{\customer}$ in \eqref{eq:lemfeasProof-2} and $\product' = \product$ in \eqref{eq:lemfeasProof-4}.
\hfill $\square$
\end{proof}

\medskip

\begin{proof}{Proof of Proposition \ref{prop:compactDPLP}.}
    Let $\product \in \productset$ and denote by $G$ the set of feasible solutions to \eqref{model:DPLP}. The projection of $G$ onto variable $\assignVar_\product$ is
    \begin{align*}
        \proj_{\assignVar_\product} G
        &=
        \{
            \assignVar_\product \in [0,1]
            \colon
            \exists\ ((\vecvaluation^1_\customer, \dots, \vecvaluation^{\numproducts}_\customer, \vecnovalVar_\customer), (\assignVar_1, \dots, \assignVar_\product, \dots, \assignVar_{\numproducts}, \assignVar_{\nopurchase})) \in G
        \}\\
        &=
        \{
            \assignVar_\product \in [0,1]
            \colon
            \not \exists \vecvaluation^\product_\customer \in \valoptionset^{\product}_{\customer}(\vecprice)
            \Rightarrow \assignVar_\product = 0
        \} \\
        &=
        \{
            \assignVar_\product \in [0,1]
            \colon
            (\price_{\product} - \price_{\choice{\customer}} > \histprice_{\customer \product} - \histprice_{\customer \choice{\customer}}) \Rightarrow \assignVar_\product = 0
        \}\\
        &=
        \{
            \assignVar_\product \in [0,1]
            \colon
            \assignVar_\product \le \indicator(\price_{\product} - \price_{\choice{\customer}} \le \histprice_{\customer \product} - \histprice_{\customer \choice{\customer}})
        \},
    \end{align*}
    where the second-to-last equality follows from Lemma \ref{lem:feasibleOptions}. The same arguments follow for $\assignVar_{\nopurchase}$. Since the objective is defined only in terms of $\vecassignVar$, we can replace the inequalities of $G$ by the projections depicted above, which results in the equivalent formulation
    \begin{align}
        \wcfun_\customer(\vecprice)
        =
        \min_{\vecassignVar \ge 0}
            &\quad
            \sum_{\product \in \productset} \price_\product \assignVar_\product
            \\
        \textnormal{s.t.}
            &\quad
                \sum_{\product \in \productset} \assignVar_{\product} + \assignVar_{\nopurchase} = 1,
                \label{eq:propCompact1} \\
            &\quad
                \assignVar_{\nopurchase} \le \indicator(\price_{\choice{\customer}} \ge \histprice_{\customer \choice{\customer}}),
                \label{eq:propCompact2} \\
            &\quad
                \assignVar_{\product} \le \indicator(\price_{\product} - \price_{\choice{\customer}} \le \histprice_{\customer \product} - \histprice_{\customer \choice{\customer}}),
                &\forall \product \in \productset.
                \label{eq:propCompact3}
    \end{align}
    Finally, note from \eqref{eq:propCompact1} and \eqref{eq:propCompact2} that
    \begin{align*}
        \sum_{\product \in \productset} \assignVar_{\product}
        =
        1 - \assignVar_{\nopurchase}
        \ge
        1 - \indicator(\price_{\choice{\customer}} \ge \histprice_{\customer \choice{\customer}})
        =
        \indicator(\price_{\choice{\customer}} < \histprice_{\customer \choice{\customer}})
    \end{align*}
    must hold at any feasible solution, and particularly tight at optimality since it is the only constraint that bounds $\vecassignVar$ from below besides the non-negativity conditions. \hfill $\square$
\end{proof}

\medskip

\begin{proof}{Proof of Proposition \ref{prop:propOptimalityDuals}.}
    Statement \ref{it:boundOptPrice} follows from the fact that product $\choice{\customer}$ is always feasible (Lemma \ref{lem:feasibleOptions}-\ref{it:feasibleChoice}) and that the inequality \eqref{eq:feasibleDualOriginal} for $\product=\choice{\customer}$ holds with $\dualJ^{*}_{i\choice{\customer}} = 0$ at optimality, while from inequality \eqref{eq:onechoiceCustomer_compact}, we have that if $\price_{\choice{\customer}}\geq \histprice_{\customer \choice{\customer}}$, $\wcfun_\customer(\vecprice)=0$.
    For statement \ref{it:boundPrices}, if $\wcfun_\customer(\vecprice) = 0$, the condition trivially holds. Assume $\wcfun_\customer(\vecprice) > 0$. Let $P'=\max_{\{\product \in \productset:p_\product<\maxprice\}} p_\product$. It can be easily shown that $P'$ always exists. Moreover, let $\product' \in \productset$ be a product such that $\price_{\product'} \geq \maxprice$. If $\product' = \choice{\customer}$, then $\wcfun_\customer(\vecprice) = 0$ which is a contradiction. Thus, $\price_{\choice{\customer}} < \maxprice$. If
    $\product' \neq \choice{\customer}$, $\wcfun_\customer(\vecprice) < \price_{\product'}$ because of \ref{it:boundOptPrice} and $\price_{\choice{\customer}} < \price_{\product'}$. Reducing the price of $\product'$ to $P'$ therefore does not change the optimal value of \eqref{model:DPDual}, and the same argument can be repeated for other products.
    \hfill $\square$
\end{proof}

\medskip

\begin{proof}{Proof of Proposition \ref{prop:optimalityReformulation}.}
    It suffices to show that both optimal solution values match when conditioned to a fixed $\vecprice \ge 0$. First,
    by Proposition \ref{prop:propOptimalityDuals}-\ref{it:boundPrices}, we can restrict our analysis to $\price_\product < \maxprice$ for all $\product \in \productset$ without loss of generality.

    Consider any customer $\customer \in \customerset$. If $\price_{\choice{\customer}} \geq \histprice_{\customer \choice{\customer}}$, then we must have $\allocVar_{\customer \choice{\customer}} = 0$ in \eqref{model:OPB}; otherwise, we can assume $\allocVar_{\customer \choice{\customer}} = 1$ since that can only be benefitial to the objective. Thus, the objective functions of both models match.

    Suppose now that $\price_{\product} - \price_{\choice{\customer}} \le \histprice_{\customer \product} - \histprice_{\customer \choice{\customer}}$ for $\product \neq \choice{\customer}$, i.e., product $\product$ is feasible to purchase by customer $\customer$. We necessarily must have $\allocVar_{\customer \product} = 1$ because of \eqref{eq:indicatorB} and, thus, \eqref{eq:ctDualCondition} and \eqref{eq:OPBdualfeasible} match. If otherwise $\price_{\product} - \price_{\choice{\customer}} > \histprice_{\customer \product} - \histprice_{\customer \choice{\customer}}$, then we may have either $\allocVar_{\customer \product} = 0$ or
    $\allocVar_{\customer \product} = 1$. Since the objective of \eqref{model:OPB} maximizes $\dualS_\customer$, we can assume
    $\allocVar_{\customer \product} = 0$, which can only relax the bound on $\dualS_\customer$ in \eqref{eq:OPBdualfeasible}. Thus, \eqref{eq:ctDualCondition} and \eqref{eq:OPBdualfeasible} also match in this case, i.e., at optimality, the values of $\dualS_\customer$ (and hence the optimal values of both models) are the same. \hfill $\square$
\end{proof}

\medskip

\begin{proof}{Proof of Theorem \ref{thm:precisionOpt}.}
    We first show \ref{it:thmFinite}. First, \eqref{model:OPCL} with $\eps = 0$ is always feasible. That is, for any $\vecprice \ge 0$ and customer $\customer \in \customerset$, we set $\allocVar_{\customer \choice{\customer}} = 1$ if and only if $\price_{\choice{\customer}} \le \histprice_{\customer \choice{\customer}}$, for all $\product \in \productset \setminus \{\choice{\customer}\},\ \allocVar_{\customer \product} = 0$ if and only if $\price_\product - \price_{\choice{\customer}} \ge \histprice_{\customer \product} - \histprice_{\customer \product}$, and $\dualS_\customer$ appropriately to satisfy \eqref{eq:OPBdualfeasible-cl}. Thus, it remains to show that $\approxfun(0)$ is bounded from above. This follows from noting that, for all $\customer \in \customerset$, $\price_{\choice{\customer}}$ is bounded by $\histprice_{\customer \choice{\customer}} + \maxprice$ in inequality \eqref{eq:indicatorA-cl} and that $\dualS_\customer$ is bounded by $\price_{\choice{\customer}}$ in inequality \eqref{eq:OPBdualfeasible-cl} with $\product = \choice{\customer}$.

    \medskip
    For \ref{it:thmPositive}, let $(\vecprice^0, \bm{\vecdualS}^0, \vecallocVar^0)$ be an optimal solution tuple and consider any $\product' \in \productset$ such that $\price^0_{\product'} = 0$. We show an alternative feasible solution with the same (optimal) value after increasing $\price^0_{\product'}$ as in the statement. If $\allocVar^0_{\customer \choice{\customer}} = 0$ for some customer $\customer$, then increasing $\price^0_{\product'}$ does not affect $\dualS_\customer$ nor the final solution value, given that $\allocVar^0_{\customer \product}$ is adjusted appropriately for $\product \neq \choice{\customer}$ to ensure feasibility. If otherwise $\allocVar^0_{\customer \choice{\customer}} = 1$ for some customer $\customer$, we have two cases:
    \begin{enumerate}
        \item Case 1, $\dualS^0_\customer = 0$. In such a scenario, we can equivalently set $\allocVar^0_{\customer \choice{\customer}} = 0$ and apply the same adjustments to $\allocVar^0_{\customer \product}$ for all $\product \neq \choice{\customer}$ as above, preserving the solution value.
        \item Case 2, $\dualS^0_\customer > 0$. Due to inequality \eqref{eq:OPBdualfeasible-cl}, we must have
        $\product' \neq \choice{\customer}$ and $\allocVar^0_{\customer \product'} = 0$. Thus, $\price^0_{\product'} - \price^0_{\choice{\customer}} \ge \histprice_{\customer \product'} - \histprice_{\customer \choice{\customer}}$ from inequality \eqref{eq:indicatorB-cl}. Increasing $\price^0_{\product'}$ thefore just increases the left-hand side of such inequality, and hence does not impact feasibility nor the solution value.
    \end{enumerate}

    \medskip
    We now show \ref{it:thmBound}. Since the feasible set of \eqref{model:OPCL} relaxes and restricts that of \eqref{model:OPB} for $\eps = 0$ and $\eps > 0$, respectively, the inequality
    $$0 \le \approxfun(0) - \numcustomers \optrevenue \le \approxfun(0) - \sum_{\customer \in \customerset} \wcfun_{\customer}(\vecprice')$$
    follows directly. We will now show that
    $$\approxfun(0) - \sum_{\customer \in \customerset} \wcfun_{\customer}(\vecprice') \le \precision' \le \precision.$$

    Consider the ordered vector of prices $\vecprice^0 > 0$ in the statement and the associated tuple $(\vecprice^0, \bm{\vecdualS}^0, \vecallocVar^0)$. Let $\customer \in \customerset$ be a customer such that $\allocVar^0_{\customer \choice{\customer}}\dualS^0_\customer > 0$. The inequality \eqref{eq:OPBdualfeasible-cl} is tight for some $\product' \le \choice{\customer}$, i.e., we can have $\allocVar^0_{\customer \product} = 0$ for all $\product < \product'$. This implies that $\price^0_{\product} - \price^0_{\choice{\customer}} \ge \histprice_{\customer \product} - \histprice_{\customer \choice{\customer}}$ for those indices due to \eqref{eq:indicatorB-cl}.

    Next, for ease of notation, let $\sigma = \precision'/(\numcustomers\numproducts)$ so that $\vecprice' = (\price^0_1 - \sigma, \price^0_2 - 2\sigma, \dots, \price^0_\numproducts - \numproducts \sigma)$. For $\product < \product'$,
    \begin{align*}
        \price'_{\product} - \price'_{\choice{\customer}} = \price^0_{\product} - \price^0_{\choice{\customer}} + (\choice{\customer} - \product)\sigma > \histprice_{\customer \product} - \histprice_{\customer \choice{\customer}},
    \end{align*}
    since $\product < \choice{\customer}$. Thus, when evaluating $\wcfun_\customer(\vecprice')$, the constraints \eqref{eq:feasibleDual} in \eqref{model:DPD} for $\product < \product'$ remains non-binding, i.e.,
    \begin{align*}
        \dualS_\customer \le \price'_\product + \indicator(\price'_{\product} - \price'_{\choice{\customer}} > \histprice_{\customer \product} - \histprice_{\customer \choice{\customer}}) \maxprice
        =
        \price'_\product + \maxprice
    \end{align*}
    for $\product < \product'$. Thus, since the price of any product $\product \ge \product'$ is decreased by $\product \sigma \le \numproducts \sigma$, evaluating the new price vector $\vecprice'$ in \eqref{model:DPD} yields: for all $\customer \in \customerset$,
    \begin{align*}
        \wcfun_\customer(\vecprice')
        \ge
        \indicator(\price'_{\choice{\customer}} < \histprice_{\customer \choice{\customer}}) (\dualS^0_{\customer} - \numproducts \sigma)
        =
        \indicator(\price^0_{\choice{\customer}} - \choice{\customer} \sigma < \histprice_{\customer \choice{\customer}}) (\dualS^0_{\customer} - \numproducts \sigma)
        \ge \allocVar^0_{\customer \choice{\customer}} \dualS^0_{\customer} - \numproducts \sigma,
    \end{align*}
    where the last inequality follows from the fact that inequality \eqref{eq:indicatorA-cl} implies $\price^0_{\choice{\customer}} \le \histprice_{\customer \choice{\customer}}$, which implies $\price^0_{\choice{\customer}} - \choice{\customer} \sigma < \histprice_{\customer \choice{\customer}}$.
    Finally, summing the above inequality over all $\customer$, we have:
    \begin{align*}
        \sum_{\customer \in \customerset} \wcfun_\customer(\vecprice') \ge \sum_{\customer \in \customerset} (\allocVar^0_{\customer \choice{\customer}}\dualS^0_\customer - \numproducts \sigma) = \approxfun(0) - \precision',
    \end{align*}
    concluding the proof. \hfill $\square$
\end{proof}

\smallskip

\begin{proof}{Proof of Proposition \ref{prop:same-price}.}
    We first show by contradiction that all prices are the same at optimality. 
    To this end, note from inequality \eqref{eq:indicatorB-cl} with $\eps = 0$ that, for any $\customer \in \customerset$, if a product $\product \neq \choice{\customer}$ is \textit{not} eligible for purchase (i.e., $\allocVar_{\customer \product} = 0$) then $\price_\product \ge \price_{\choice{\customer}}$, since all historical prices for $\customer$ are the same. Consider now an optimal solution $(\vecprice^*, \bm{\vecdualS}^*, \vecallocVar^*)$ and let $\price^{\min} \equiv \min_{\product \in \productset} \price^*_\product$ be the minimum optimal price. Suppose there exists some customer $\customer \in \customerset$ who selects a product $\product$ such that $\dualS^*_\customer = \price^*_\product > \price^{\min}$. But this implies that 
    $\price^{\min} < \price^*_\product \le \price^*_{\choice{\customer}}$ (since $\choice{\customer}$ is always feasible), i.e., we must have $\allocVar^*_{\customer \product} = 1$ which by inequality \eqref{eq:OPBdualfeasible-cl} implies that $\dualS^*_\customer \le \price^{\min}$, a contradiction.

    Finally, let $\price^*$ be the optimal (scalar) price of all products, and suppose $\customer' \in \customerset$ is the smallest customer index such that 
        $\price^* \le \histprice_{\customer'}.$ It follows that any customer $\customer < \customer'$ does not purchase any product (since $\price^* = \price^*_{\choice{\customer}} > \histprice_\customer$), while all customers 
    $\customer \ge \customer'$ yield a revenue of $\price^*$ (since $\price^* = \price^*_{\choice{\customer}} \le \histprice_\customer$). Thus, we must have 
    $\price^* = \histprice_{\customer'}$ at optimality, and the total revenue is 
    $(\numcustomers - \customer' + 1)\histprice_{\customer'}$. The element $\customer^*$ in the proposition statement is the index that maximizes this revenue.
    \hfill $\square$
\end{proof}

\smallskip

\begin{proof}{Proof of Proposition \ref{prop:constant-price}.}
    The solution value of the proposed solution (i.e., a lower bound to the problem) is $\sum_{\customer \in \customerset} \histprice_{\choice{\customer}}$, since all customers would purchase their choice $\choice{\customer}$. Due to inequality \eqref{eq:OPBdualfeasible-cl}, this is also an upper bound, and hence is optimal.
    \hfill $\square$
\end{proof}

\smallskip

\begin{proof}{Proof of Proposition \ref{prop:consv-pricing}.}
By definition \eqref{eq:conservative-price}, the total revenue is at least $\numcustomers \ubar{\histprice}$. Conversely, by inequalities \eqref{eq:OPBdualfeasible-cl} and \eqref{eq:indicatorA-cl}, the total revenue is at most $\sum_{\customer \in \customerset}\histprice_{\customer \choice{\customer}} \le \numcustomers \bar{\histprice}$. The ratio hence follows from dividing the lower bound by the upper bound.

We now construct an instance where this ratio is asymptotically tight. 
Consider $\numproducts = 1$ product, $\numcustomers > 1$ customers, and fix 
any $\histprice_1, \histprice_2$ such that $\histprice_1 < \histprice_2$. Customer $1$ purchases the product at price $\histprice_1$, while the remaining customers $2, \dots, \numcustomers$ purchase it at $\histprice_2$. The conservative pricing \eqref{eq:conservative-price} will set the product's price at $\histprice_1$, yielding a total revenue of $\numcustomers \histprice_1$. For any sufficiently large $\numcustomers$, the optimal pricing strategy sets $\histprice_2$ as the optimal price, yielding a total revenue of $(\numcustomers-1)\histprice_2$ (since customer 1 will not purchase any product). The performance ratio is therefore $\numcustomers \histprice_1 / (\numcustomers-1)\histprice_2$. Taking the limit $\numcustomers \rightarrow +\infty$ with respect to this ratio completes the proof.
\hfill $\square$
\end{proof}

\smallskip

\begin{proof}{Proof of Proposition \ref{prop:cutoff-pricing-bound}.}
Without loss of generality, suppose the customer index set $\customerset$ is ordered according to historical purchase prices, i.e.,
$0 < \ubar{\histprice} 
= 
\histprice_{1 \choice{1}} \le \histprice_{2 \choice{2}} \le \dots \le \histprice_{\numcustomers \choice{\numcustomers}} = \bar{\histprice}$. Let $\dualS^{\textnormal{OPT}}$ and $\dualS^{\textnormal{CP}}$ be the optimal solution value of \eqref{model:OPMIP} and the total revenue obtained by the cut-off price \eqref{eq:cutoff-price}, respectively. It follows from \eqref{eq:OPBdualfeasible-cl} for $\product = \choice{\customer}$ that $\dualS^{\textnormal{OPT}} \le \sum_{\customer \in \customerset} \histprice_{\customer \choice{\customer}}$. Furthermore, given the price ordering, notice that \eqref{eq:modelOCPCsimplification} evaluates to 
$\max_{\customer \in \customerset} (\numcustomers - \customer + 1)\histprice_{\customer \choice{\customer}}$. By the cut-off pricing definition \eqref{eq:cutoff-price}, we therefore have $\dualS^{\textnormal{CP}} \ge \max_{\customer \in \customerset} (\numcustomers - \customer + 1)\histprice_{\customer \choice{\customer}}$. 
This is because due to the incentive-compatibility constraint \eqref{eq:indicatorB-cl} some customers might not purchase the products priced at the cut-off price $p^*$ and opt for a product that is priced higher than $p^*$. Thus,
\begin{align}\label{proof:cutoff1}
    \frac{\dualS^{\textnormal{CP}}}{\dualS^{\textnormal{OPT}}} 
    \ge 
    \frac{
        \max_{\customer \in \customerset} (\numcustomers - \customer + 1)\histprice_{\customer \choice{\customer}}
    }
    {
        \sum_{\customer \in \customerset} \histprice_{\customer \choice{\customer}}
    }.
\end{align}
Next, if we divide both the numerator and denominator of the right-hand side ratio above by the number of customers $\numcustomers \ge 1$, we obtain 
\begin{align} \label{proof:cutoff2}
    \frac{
        \max_{\customer \in \customerset} (\numcustomers - \customer + 1)\histprice_{\customer \choice{\customer}}
    }
    {
        \sum_{\customer \in \customerset} \histprice_{\customer \choice{\customer}}
    }
    =
    \cfrac{
        \cfrac{
            \max_{\customer \in \customerset} (\numcustomers - \customer + 1)\histprice_{\customer \choice{\customer}}
        }
        {
            \numcustomers
        }
    }
    {
        \cfrac{
            \sum_{\customer \in \customerset} \histprice_{\customer \choice{\customer}}
        }
        {
            \numcustomers
        }
    }
    =
    \frac{
        \max_{\customer \in \customerset} [1 - \cdf(\histprice_{\customer \choice{\customer}})] \histprice_{\customer \choice{\customer}}
    }
    {
        \mathbb{E}[\rv]
    },
\end{align}
where $\rv$ is a non-negative discrete random variable uniformily distributed on the set $\{\histprice_{1 \choice{1}}, \dots, \histprice_{\numcustomers \choice{\numcustomers}}\}$, and $\cdf(\cdot)$ and $\mathbb{E}[\rv]$ denote the left continuous c.d.f. (i.e.,$\cdf(x)\equiv P(\rv<x)$) and the expectation of $\rv$, respectively. This problem now bears resemblance to the personalized pricing problem
studied in \cite{elmachtoub2018value}.

Let $\maxnum \equiv \max_{\customer \in \customerset} [1 - \cdf(\histprice_{\customer \choice{\customer}})] \histprice_{\customer \choice{\customer}}$ be the numerator of the ratio above. We have $\ubar{\histprice} \le \maxnum \le \bar{\histprice}$; in particular, the left-hand side inequality holds since 
$[1 - \cdf(\ubar{\histprice})] \ubar{\histprice}=\ubar{\histprice}$. We can hence rewrite the expectation term as 
\begin{align*}
    \mathbb{E}[\rv] 
    =
    \sum_{\customer \in \customerset} \frac{1}{\numcustomers} \histprice_{\customer \choice{\customer}}
    &=
    \sum_{\customer \in \customerset} \frac{1}{\numcustomers} 
    \left[
        \int_{0}^{\histprice_{\customer \choice{\customer}}} 1\,\mathrm{d}x    
    \right]\\
    &=
    \sum_{\customer \in \customerset} \frac{1}{\numcustomers} 
    \left[
        \int_{0}^{\histprice_{\customer \choice{\customer}}} 1\,\mathrm{d}x    
        +
        \int_{\histprice_{\customer \choice{\customer}}}^{+\infty} 0\,\mathrm{d}x    
    \right]\\
    &=
    \sum_{\customer \in \customerset} \frac{1}{\numcustomers} 
    \left[
        \int_{0}^{+\infty} \indicator(x \leq \histprice_{\customer \choice{\customer}})\,\mathrm{d}x    
    \right]\\
    &=
    \int_{0}^{+\infty}
    \sum_{\customer \in \customerset} 
    \left[
        \frac{1}{\numcustomers} 
        \indicator(x \leq \histprice_{\customer \choice{\customer}})
    \right]\,\mathrm{d}x\\
    &=
    \int_{0}^{+\infty} [1 - \cdf(x)]\,\mathrm{d}x \\
    &=
    \int_{0}^{\bar{\histprice}} [1 - \cdf(x)]\,\mathrm{d}x 
    \\
    &=
    \int_{0}^{\maxnum} [1 - \cdf(x)]\,\mathrm{d}x 
    +
    \int_{\maxnum}^{\bar{\histprice}} [1 - \cdf(x)]\,\mathrm{d}x \\
    &\le 
    \maxnum + \int_{\maxnum}^{\bar{\histprice}} [1 - \cdf(x)]\,\mathrm{d}x 
    \\
    &\le 
    \maxnum + \int_{\maxnum}^{\bar{\histprice}} \frac{\maxnum}{x} \,\mathrm{d}x 
    \\
    &=
    \maxnum + \maxnum \log \left( \frac{\bar{\histprice}}{\maxnum} \right)
    \\
    &\le 
    \maxnum + \maxnum\log \left( \frac{\bar{\histprice}}{\ubar{\histprice}} \right),
\end{align*}
where the previous-to-the-last inequality follows because $1 - \cdf(x) \le \maxnum/x$ for any $x \in [\ubar{\histprice}, \bar{\histprice}]$ by definition. From the inequality above, we obtain
\begin{align*}
\frac{\mathbb{E}[\rv]}{\maxnum} \le 1 + \log \left( \frac{\bar{\histprice}}{\ubar{\histprice}} \right)
\Leftrightarrow
\frac{\dualS^{\textnormal{CP}}}{\dualS^{\textnormal{OPT}}} 
\ge
\frac{\maxnum}{\mathbb{E}[\rv]} 
\ge 
\frac{1}{
    1 + \log \left( \frac{\bar{\histprice}}{\ubar{\histprice}} \right)
}.
\end{align*}
\color{mycolor}
Finally, from \eqref{proof:cutoff1} and \eqref{proof:cutoff2} we have that
\begin{align*}
   \frac{\dualS^{\textnormal{CP}}}{\dualS^{\textnormal{OPT}}} 
    \ge 
    \frac{
        \max_{\customer \in \customerset} [1 - \cdf(\histprice_{\customer \choice{\customer}})] \histprice_{\customer \choice{\customer}}
    }
    {
        \mathbb{E}[\rv]
    }
    \ge \frac{
        (1 - \cdf(\histprice_{k \choice{k}})) \histprice_{k \choice{k}}
    }
    {
        \mathbb{E}[\rv]
    }\ge \frac{\histprice_{k \choice{k}}}{2\mathbb{E}[\rv]} \ge \frac{med(P)}{2\mathbb{E}[\rv]}
\end{align*}
where $med(P)$ denotes the median of $\rv$, while $\histprice_{k \choice{k}}=\min_{\customer\in\customerset}\{\histprice_{\customer\choice{\customer}}:\histprice_{\customer\choice{\customer}}\ge med(P)\}$. The prior to the last inequality follows from the fact that $1-\cdf(\histprice_{k \choice{k}})\ge \frac{1}{2}$ (as $\cdf(\cdot)$ is defined to be left continuous), while the last inequality follows from $\histprice_{k \choice{k}}\ge med(P)$, proving the performance bound.
\color{black}

It remains to show that the ratio is asymptotically tight. Consider an instance with any number $\numcustomers > 0$ of customers  and $\numcustomers - k + 1$ products, where
$k \in \mathbb{N}$ is any positive integer such that $k \le \numcustomers$. The historical observed prices $\vechistprice_{\customer}$ by customer $\customer$, in turn, are defined by the following vectors, where $\delta$ is any scalar such that $0 < \delta < 1/(\numcustomers-1)$:
\begin{align*}
    \vechistprice_{1} 
    &= 
    \left(
        \frac{\numcustomers}{\numcustomers},\frac{\numcustomers}{\numcustomers-1}, 
        \frac{\numcustomers}{\numcustomers-2},\dots, \frac{\numcustomers}{k+1}, 
        \frac{\numcustomers}{k}
    \right),
    \\
    \vechistprice_{2} 
    &= 
    \left(
        \frac{\numcustomers}{\numcustomers} + \delta,\frac{\numcustomers}{\numcustomers-1}, 
        \frac{\numcustomers}{\numcustomers-2},\dots, \frac{\numcustomers}{k+1}, \frac{\numcustomers}{k}
    \right),
    \\
    \vechistprice_{3} 
    &= 
    \left(
        \frac{\numcustomers}{\numcustomers} + \delta,\frac{\numcustomers}{\numcustomers-1} + \delta, 
        \frac{\numcustomers}{\numcustomers-2},\dots, \frac{\numcustomers}{k+1}, 
        \frac{\numcustomers}{k}
    \right),
    \\
    \dots
    \\
    \vechistprice_{\numcustomers-k} 
    &= 
    \left(
        \frac{\numcustomers}{\numcustomers} + \delta,\frac{\numcustomers}{\numcustomers-1}+ \delta, 
        \frac{\numcustomers}{\numcustomers-2} + \delta, 
        \dots, \frac{\numcustomers}{k+1}, \frac{\numcustomers}{k}
    \right),
    \\
    \vechistprice_{\numcustomers-k+1} 
    &= 
    \left(
        \frac{\numcustomers}{\numcustomers} + \delta,\frac{\numcustomers}{\numcustomers-1}+ \delta, 
        \frac{\numcustomers}{\numcustomers-2} + \delta, 
        \dots, \frac{\numcustomers}{k+1} + \delta, \frac{\numcustomers}{k}
    \right),
    \\
    \vechistprice_{\numcustomers-k+2} 
    &= 
    \left(
        \frac{\numcustomers}{\numcustomers} + \delta,\frac{\numcustomers}{\numcustomers-1}+ \delta, 
        \frac{\numcustomers}{\numcustomers-2} + \delta, 
        \dots, \frac{\numcustomers}{k+1} + \delta, \frac{\numcustomers}{k}
    \right),
    \\
    \dots
    \\
    \vechistprice_{\numcustomers} 
    &= 
    \left(
        \frac{\numcustomers}{\numcustomers} + \delta,\frac{\numcustomers}{\numcustomers-1}+ \delta, 
        \frac{\numcustomers}{\numcustomers-2} + \delta, 
        \dots, \frac{\numcustomers}{k+1} + \delta, \frac{\numcustomers}{k}
    \right).
\end{align*}

For the historical purchase choices, each customer $\customer \in \{1,\dots,\numcustomers-k\}$ purchased product $\choice{\customer} = \customer$, while all remaining customers $\customer' \in \{\numcustomers-k+1, \dots, \numcustomers\}$ pick the same product $\choice{\customer'} = \numcustomers-k+1$.

An upper bound on the optimal revenue $\dualS^{\textnormal{OPT}}$ for this instance is:
\begin{align*}
    \dualS^{\textnormal{OPT}} 
    \le 
    \sum_{\customer \in \customerset} \histprice_{\customer \choice{\customer}} 
    =  \sum_{\customer = 1}^{\numcustomers - k+1} \histprice_{\customer \customer} 
    +
     \sum_{\customer = \numcustomers - k+2}^{\numcustomers} \histprice_{\customer (m-k+1)}
    =  \sum_{\customer = 0}^{\numcustomers - k} \frac{\numcustomers}{\numcustomers - \customer} 
    +
     \sum_{\customer = \numcustomers - k+1}^{\numcustomers-1} \frac{m}{k}
    =
    \sum_{\customer = 0}^{\numcustomers - k} \frac{\numcustomers}{\numcustomers - \customer} 
    + 
    \numcustomers \frac{k-1}{k}.
\end{align*}
We now define prices whose resulting total revenue is arbitrarily close to $\dualS^{\textnormal{OPT}}$. Specifically, consider the price vector
$$
    \vecprice = 
    \left(
        \frac{\numcustomers}{\numcustomers},
        \frac{\numcustomers}{\numcustomers-1} - \delta,
        \frac{\numcustomers}{\numcustomers-2} - 2\delta,
        \dots,
        \frac{\numcustomers}{k} - (\numcustomers-k)\delta
    \right).
$$
From the definition above, since $\price_{\choice{\customer}} \leq \histprice_{\customer \choice{\customer}}$ for all $\customer \in \customerset$, by inequality \eqref{eq:indicatorA-cl} every customer purchases a product. We next show that customer $\customer$ purchases product $\choice{\customer}$. For any $\product < \choice{\customer}$,
\begin{align*}
    \price_\product - \price_{\choice{\customer}}
    =
    \frac{m}{m-j+1} - \frac{m}{m-\choice{\customer}+1}
    - \product \delta + \choice{\customer}\delta 
    \geq 
    \histprice_{\customer \product} - \histprice_{\customer \choice{\customer}},
\end{align*}
i.e., such products $\product$ violates incentive-compatibility constraints (inequality \eqref{eq:indicatorB-cl} with $\epsilon=0$) and will not be purchasable by customer $\customer$. Moreover, from our choice of $\delta< 1/(\numcustomers-1)$, it can be easily verified that $\price_1 < \price_2 < \dots < \price_{\numcustomers-k+1}$ and therefore $\price_\product > \price_{\choice{\customer}}$ for all $\customer$ and $\product > \choice{\customer}$. Product $\price_{\choice{\customer}}$ must be necessarily chosen by customer $\customer$ under her worst-case valuation ($v_{\customer \choice{\customer}}=\histprice_{\customer \choice{\customer}}$ and $v_{\customer k}=0$ for all $k\in \productset\setminus\{\choice{\customer}\}$) and the revenue from this pricing is hence 
\begin{align*}
    \sum_{\customer \in \customerset} \price_{\choice{\customer}}
    =
    \sum_{\customer = 1}^{\numcustomers - k+1} \price_{\customer}
    +
    (k-1) \price_{m-k+1}
    =
    \sum_{\customer = 0}^{\numcustomers - k} 
    \left(
        \frac{\numcustomers}{\numcustomers - i} 
        -
        \customer\delta
    \right)
    + 
    \numcustomers \frac{k-1}{k}- (k-1)(\numcustomers - k)\delta,
\end{align*}
thus, as $\delta \rightarrow 0$, the total revenue obtained from $\vecprice$ approximates that of $\dualS^{\textnormal{OPT}}$.

We now show that the revenue obtained from the cut-off pricing \eqref{eq:cutoff-price} is $\numcustomers$. Suppose that the customer index that solves \eqref{eq:modelOCPCsimplification} is $\customer'$, and hence $\price^* = \histprice_{\customer' \choice{\customer'}}$. If $\choice{\customer'}<m-k+1$, then the cut-off prices are $p^{CP}_j=\frac{m}{k}$ for all $j<\choice{\customer'}$, $p^{CP}_j=\frac{m}{m-j+1}$ for all $\choice{\customer'}\leq j<m-k+1$ and $p^{CP}_{m-k+1}=\frac{m}{k}$. However, if $\choice{\customer'}=m-k+1$, the cut-off prices are $p^{CP}_j=\frac{m}{k}$ for all $j\in \productset$. Suppose that $\choice{\customer'}<m-k+1$. By the construction of historical prices, 
$$
\histprice_{\customer' \choice{\customer'}} 
= 
\frac{\numcustomers}{\numcustomers - \choice{\customer'} + 1},
$$
and therefore $\numcustomers - \choice{\customer'} + 1$ customers would purchase a product since $\histprice_{\customer \choice{\customer}} \ge \histprice_{\customer' \choice{\customer'}}$ for all $\customer \ge \customer'$. Analogously, since by construction for all customers $\customer < \customer'$, $\histprice_{\customer \choice{\customer}}<\histprice_{\customer' \choice{\customer'}}$,
none of the worst-case historical customers $\customer < \customer'$ would purchase any product.
Moreover, for any of the worst-case customer $\customer \ge \customer'$, the cut-off price \eqref{eq:cutoff-price} of its chosen product $\choice{\customer}$ is 
$$
\price^{\textnormal{CP}}_{\choice{\customer}}
= 
\frac{\numcustomers}{\numcustomers - \choice{\customer} + 1}.
$$
This implies that, for any $\customer > \customer'$,
\begin{align*}
    \price^{\textnormal{CP}}_{\choice{\customer'}} - \price^{\textnormal{CP}}_{\choice{\customer}}
    =
    \frac{\numcustomers}{\numcustomers - \choice{\customer'} + 1}
    -
    \frac{\numcustomers}{\numcustomers - \choice{\customer} + 1}
    <
    \left( 
        \frac{\numcustomers}{\numcustomers - \choice{\customer'} + 1}
        + \delta
    \right)
    -
    \frac{\numcustomers}{\numcustomers - \choice{\customer} + 1}
    =
    \histprice_{\customer' \choice{\customer'}} - \histprice_{\customer \choice{\customer}},
\end{align*}
i.e., product $\choice{\customer'}$ is incentive-compatible with all customers $\customer > \customer'$. Because $\price^{\textnormal{CP}}_{\choice{\customer'}}$ is the lowest price across all the products, the total revenue of the cut-off solution  
\begin{align*}
    \dualS^{\textnormal{CP}} 
    = 
    (\numcustomers - \choice{\customer'} + 1) \histprice_{\customer' \choice{\customer'}}
    =
    (\numcustomers - \choice{\customer'} + 1)\frac{\numcustomers}{\numcustomers - \choice{\customer'} + 1}
    =
    \numcustomers
\end{align*}
If $\choice{\customer'}=m-k+1$, as mentioned before the cut-off prices are $p^{CP}_j=\frac{m}{k}$ for all $j\in \productset$ and it can also be easily confirmed that $\dualS^{\textnormal{CP}}=m$ as only the $k$ historical customers who historically purchased product $m-k+1$ will make a purchase, under price $\frac{m}{k}$, in the worst-case.

Finally, as $\delta \rightarrow 0$,
\begin{align*}
    \frac{\dualS^{\textnormal{CP}}}{\dualS^{\textnormal{OPT}}}
    \rightarrow
    \cfrac{\numcustomers}
    {
        \sum_{\customer = 0}^{\numcustomers - k} 
        \cfrac{\numcustomers}{\numcustomers - \customer}         
        + 
        \numcustomers \cfrac{k-1}{k}
    }
    &=
    \cfrac{\numcustomers}
    {
        \numcustomers 
        \left(
            \sum_{\customer = 0}^{\numcustomers - k} 
            \cfrac{1}{\numcustomers - \customer}         
            + 
            \cfrac{k-1}{k}
        \right)
    }
    \\
    &=
    \cfrac{1}
    {
        \cfrac{1}{\numcustomers} + \dots + \cfrac{1}{k} - \cfrac{1}{k} + \cfrac{k}{k}
    }
    \\
    &=
    \cfrac{1}
    {
        \sum_{\customer = k+1}^{\numcustomers} \cfrac{1}{\customer} + \cfrac{k}{k}
    }.
\end{align*}

For a sufficiently large $\numcustomers$ and $k$ (e.g., by multiplying both by the same constant), the ratio above can be made sufficiently close to   notespar,

\begin{align*}
    \cfrac{1}{\log \numcustomers - \log k + \cfrac{1}{2\numcustomers} - \cfrac{1}{2k} + 1}
    =
    \cfrac{1}{\log\left(\cfrac{\numcustomers}{k}\right) + \cfrac{1}{2\numcustomers} - \cfrac{1}{2k} + 1},
\end{align*}

where the last equality follows from the fact that for large enough $n$, $\sum_{\I=1}^{n} \frac{1}{\I}=\log(n)+\gamma+\frac{1}{2n}$ where $\gamma$ is the Euler–Mascheroni constant. The ratio above can approximate $1/(1+\log(\bar{\histprice}/\ubar{\histprice}))$ at any desired precision since $\numcustomers/k$ can be made sufficiently close to $\bar{\histprice}/\ubar{\histprice}$, while both numbers are also sufficiently large. Thus, the ratio is asymptotically tight for the constructed instance.
\hfill $\square$

\end{proof}

\section{Extending External Validity to the MNL Model}
\label{exvalidmnl}
When there are multiple products, the external validity analysis of our approach becomes significantly more complicated,
primarily because both the model-free and model-based optimal prices are not tractable.
Given such hurdles, we focus on a simple yet relevant special case, where for some $\beta>0$, the probability of customer $\customer$ choosing product $\product$ from the assortment is given by
\begin{equation}\label{probmnl}
    \frac{\exp(\alpha-\beta \histprice_{\customer\product})}{1+\sum_{k=1}^\numproducts \exp(\alpha-\beta \histprice_{\customer k})}.
\end{equation}
That is, we limit our scope to the situation where customers make decisions based on the MNL choice model and the average attractiveness of the products is symmetric (the products are horizontally differentiated). 
Moreover, we assume the historical prices satisfy $\histprice_{\customer\product}\equiv P_{\customer}$ for all $\customer\in \customerset$.
Such a uniform price is common in some settings, e.g., as noted by \cite{draganska2006consumer}, it is a well established policy to offer all flavors in a product line of yogurt at the same price. 
In this case, implied by Proposition~\ref{prop:same-price}, model-free optimal prices are the same for all the products and can be expressed in a closed form for uniformly distributed historical prices. 
We denote the model-free optimal price of all the products as $p^*$. 
Moreover, given that customers make decisions based on the MNL choice model, the model-based optimal prices are equal for all the products, denoted by $\hat{\price}$. 
We also define the expected revenue of the firm when all the products are priced at $p$ under the MNL model~\eqref{probmnl} as $\mrevfunc(p)$.
It turns out the insights developed for a single product in Theorem~\ref{prop:rev-bound} can be extended to this special case with multiple products.

\begin{proposition}\label{prop:mnl}
Suppose the firm sets prices $P_i$ from a distribution with PDF $g(\price) = 1/b$ for $\price\in [0,b]$, where $b>\hat p$. We have
$
    \frac{\mrevfunc(p^*)}{\mrevfunc(\hat p)}\ge  \frac{1}{2}.
$
Moreover, the bound is asymptotically tight: for any $\epsilon>0$ that is sufficiently small, a case can be constructed in which there exists $\bar{\alpha}$ such that when $\alpha > \bar{\alpha}$, we have  $\mrevfunc(p^*)/\mrevfunc(\hat p)\leq \frac{1}{2}+\epsilon$.
\end{proposition}

We note that the bound is asymptotically tight since the worst-case holds when $\alpha \to \infty$.

\begin{proof}{Proof of Proposition~\ref{prop:mnl}.}

First, we will show that $p^*$ must be unique. By Proposition~\ref{prop:same-price}, the optimal price output by our framework \ref{model:OPMIP} satisfies
$\price^*\in \argmax_{\price\ge 0} \price \sum_{i=1}^\numcustomers \indicator(\histprice_{i\choice{i}}\ge \price)$. Dividing the latter part of this quantity by the number of historical customers $\numcustomers$ does not affect the value of $\price^*$. In other words, $\price^*\in \argmax_{\price\ge 0} \price \frac{\sum_{i=1}^\numcustomers \indicator(\histprice_{i\choice{i}}\ge \price)}{\numcustomers}$.

By the Law of Large Numbers we have 
\[\lim_{\numcustomers\to\infty}\frac{\sum_{i=1}^\numcustomers \indicator(\histprice_{i\choice{i}}\ge \price)}{\numcustomers}=\pr(\histprice_{i\choice{i}}\ge p) =\int_{p}^{b}\frac{1}{b}(1-\frac{1}{ne^{\alpha -\beta x}+1})\mathrm{d}x.\]

Therefore 
\[ p^*\in \argmax_{\price\ge 0} \price \lim_{\numcustomers\to\infty} \frac{\sum_{i=1}^\numcustomers \indicator(\histprice_{i\choice{i}}\ge \price)}{\numcustomers}= \argmax_{0\leq p\leq b} p\int_{p}^{b}(1-\frac{1}{ne^{\alpha -\beta x}+1})\mathrm{d}x\]

Now, notice that since \[0\int_{0}^{b}(1-\frac{1}{ne^{\alpha -\beta x}+1})\mathrm{d}x=b\int_{b}^{b}(1-\frac{1}{ne^{\alpha -\beta x}+1})\mathrm{d}x=0,\]
and 
\begin{equation}\label{eq:ourMNL}
p\int_{p}^{b}(1-\frac{1}{ne^{\alpha -\beta x}+1})\mathrm{d}x
\end{equation}
is a continuous and differentiable function, not zero everywhere in $[0,b]$ and defined on the compact region of $[0,b]$, it has a maximizer and its maximizer must be in the interior of $[0,b]$, and must satisfy the first-order condition for \eqref{eq:ourMNL}.

Now we will show there is a unique $p^*\in[0,b]$ that satisfies the first-order condition for \eqref{eq:ourMNL}, and then we will show that $p^*\leq \hat{p}$.

To that avail, notice that the first-order condition for \eqref{eq:ourMNL} with respect to $p$ is equivalent to:
\footnotesize
\[\bigg(p\int_{p}^{b}(1-\frac{1}{ne^{\alpha -\beta x}+1})\mathrm{d}x\bigg)'=\int_{p}^{b}(1-\frac{1}{ne^{\alpha -\beta x}+1})\mathrm{d}x-p(1-\frac{1}{ne^{\alpha -\beta p}+1})=\frac{\log(\frac{ne^{\alpha -\beta p}+1}{ne^{\alpha -\beta b}+1})}{\beta}-p(1-\frac{1}{ne^{\alpha -\beta p}+1})=0.\]
\normalsize
To show that $p^*$ must be unique, we need to show 
\[q(p)=\frac{\log(\frac{ne^{\alpha -\beta p}+1}{ne^{\alpha -\beta b}+1})}{\beta}-p(1-\frac{1}{ne^{\alpha -\beta p}+1})=0\]
has only one root in $[0,b]$. Notice that 
\[q(b)=\frac{\log(\frac{ne^{\alpha -\beta b}+1}{ne^{\alpha -\beta b}+1})}{\beta}-b(1-\frac{1}{ne^{\alpha -\beta b}+1})<0.\]
Moreover, \[q(0)=\frac{\log(\frac{ne^{\alpha -\beta 0}+1}{ne^{\alpha -\beta b}+1})}{\beta }-0(1-\frac{1}{ne^{\alpha -\beta 0}+1})>0,\]
thus there must be at least one root. 

Now notice that 
\[q'(p)=\bigg(\frac{\log(\frac{ne^{\alpha -\beta p}+1}{ne^{\alpha -\beta b}+1})}{\beta}-p(1-\frac{1}{ne^{\alpha -\beta p}+1})\bigg)'=\frac{n\beta e^\alpha((\beta p-2)e^{\beta p}-2ne^\alpha)}{(e^{\beta p}+ne^\alpha)^2},\]

which is $<0$ if $p<\frac{w_0(2ne^{\alpha -2})+2}{\beta}$ and is $>0$ if $p>\frac{w_0(2ne^{\alpha -2})+2}{\beta}$. If $\frac{w_0(2ne^{\alpha-2})+2}{\beta}>b$, then $q(p)$ is strictly decreasing in $[0,b]$ and has exactly one root in $[0,b]$. If $\frac{w_0(2ne^{\alpha-2})+2}{\beta}\leq b$, then since $q(b)<0$ and $q(p)$ is strictly increasing on $(\frac{w_0(2ne^{\alpha -2})+2}{\beta},b]$, $q(p)$ can have only one root, $p^*$, on $[0,b]$. 

Now, we show that $p^*\leq \hat{p}$. Assume otherwise, that for some $\alpha$, $b$, $n$ and $\beta>0$, $p^*>\hat{p}$. It is trivial to see that if \begin{equation}\label{eq:altmnl}
\frac{\log(ne^{\alpha -\beta p}+1)}{\beta}-p(1-\frac{1}{ne^{\alpha -\beta p}+1})=0
\end{equation}
has a root, its root must be larger than or equal to $p^*$. To show that \eqref{eq:altmnl} also has a unique root it suffices to notice that \[\frac{\log(ne^{\alpha -\beta 0}+1)}{\beta }-0(1-\frac{1}{ne^{\alpha -\beta 0}+1})>0,\] and \[\lim_{p\to\infty}\big(\frac{\log(ne^{\alpha -\beta p}+1)}{\beta }-p(1-\frac{1}{ne^{\alpha -\beta p}+1})\big)=0,\]
while \[\big(\frac{\log(ne^{\alpha -\beta p}+1)}{\beta }-p(1-\frac{1}{ne^{\alpha -\beta p}+1})\big)'=\frac{n\beta e^\alpha((\beta p-2)e^{\beta p}-2ne^\alpha)}{(e^{\beta p}+ne^\alpha)^2},\] which is $<0$ if $p<\frac{w_0(2ne^{\alpha-2})+2}{\beta}$ and is $>0$ if $p>\frac{w_0(2ne^{\alpha-2})+2}{\beta}$.

Now let $\bar{p}$ be the single root of \eqref{eq:altmnl}. We know $\bar{p}\geq p^*>\hat{p}$. Notice that the derivative of \eqref{eq:altmnl} with respect to $\alpha$ is 
\[\frac{ne^\alpha(ne^\alpha+(1-\beta p)e^{\beta p})}{\beta (e^{\beta p}+ne^\alpha)^2},\] 
which is $<0$ if $p>\frac{w_0(ne^{\alpha-1})+1}{\beta }=\hat{p}$. Thus, \eqref{eq:altmnl} is strictly decreasing in $\alpha$ at $\bar{p}$. Moreover, $\frac{w_0(ne^{\alpha-1})+1}{\beta}$ is increasing in $\alpha$, hence
\[\lim_{\alpha\to -\infty} \big(\frac{\log(ne^{\alpha-\beta \bar{p}}+1)}{\beta }-\bar{p}(1-\frac{1}{ne^{\alpha-\beta \bar{p}}+1})\big)>0,\] 
however this is a contradiction since
\[\frac{\log(ne^{\alpha-\beta p}+1)}{\beta }-p(1-\frac{1}{ne^{\alpha-\beta p}+1})\leq0\]
for all $p\geq0$ when $\alpha\to -\infty$. The latter claim follows from the fact that $\frac{\log(ne^{\alpha-\beta p}+1)}{\beta }$ is decreasing in $p$ and $\frac{\log(ne^{\alpha}+1)}{\beta }=0$ as $\alpha\to -\infty$ while $p(1-\frac{1}{ne^{\alpha-\beta p}+1})\geq0$ if $p\geq 0$. As a result, by contradiction we can conclude that $\bar{p}\leq \hat{p}$ which implies $p^*\leq \hat{p}$ since we know $p^*\leq \bar{p}$.
\vskip 1em

Now we show the revenue from $p^*$ is at least $\frac{1}{2}$ times the revenue from $\hat{p}$. We have $p^*\leq \hat{p}\leq b$. Notice that 

\begin{equation}\label{eq:mnlwelfderiv}
    \bigg(\frac{\log(\frac{ne^{\alpha-\beta p}+1}{ne^{\alpha-\beta b}+1})}{\beta}+p(1-\frac{1}{ne^{\alpha-\beta p}+1})\bigg)'=-\frac{n\beta pe^{\beta p+\alpha}}{(e^{\beta p}+ne^\alpha)^2}\leq 0
\end{equation}    
for all $p\geq0$. Moreover, for all $0\leq p\leq b$, 
\[\frac{\log(\frac{ne^{\alpha-\beta p}+1}{ne^{\alpha-\beta b}+1})}{\beta}\geq0,\]
thus \footnotesize\\
\[\mrevfunc(\hat p)=\hat{p}(1-\frac{1}{ne^{\alpha-\beta \hat{p}}+1})\leq \frac{\log(\frac{ne^{\alpha-\beta \hat{p}}+1}{ne^{\alpha-\beta b}+1})}{\beta}+\hat{p}(1-\frac{1}{ne^{\alpha-\beta \hat{p}}+1})\leq \frac{\log(\frac{ne^{\alpha-\beta p^*}+1}{ne^{\alpha-\beta b}+1})}{\beta}+p^*(1-\frac{1}{ne^{\alpha-\beta p^*}+1})\\=2\mrevfunc(p^*).\]
\normalsize
Where he second inequality follows from \eqref{eq:mnlwelfderiv} and the fact that $0\le p^*\le \hat{p}$. The last equality follows from the fact that $p^*$ must satisfy the first-order condition for \eqref{eq:ourMNL}. Thus, we have shown that $\mrevfunc(p^*)\geq \frac{1}{2} \mrevfunc(\hat{p})$.

It remains to show the bound is tight. Consider a case where $n=2$, $\beta=1$, and $b\gg \alpha$. Let $0<\epsilon_1\ll\epsilon$. For any $\sigma>0$, there exists $\tilde{\alpha}$ such that if $\alpha>\tilde{\alpha}$, for any $0\leq p\leq (1-\epsilon_1)\alpha$, $|\log(\frac{2e^{\alpha-p}+1}{2e^{\alpha-b}+1})-\log(2e^{\alpha-p})|<\sigma$, and $|p(1-\frac{1}{2e^{\alpha-p}+1})-p|<\sigma$. Thus, for any $\sigma>0$, there exists $\tilde{\alpha}$ such that if $\alpha>\tilde{\alpha}$, then $|p^*-\frac{\alpha+\log(2)}{2}|<\sigma$ (i.e., $p^*$ can be arbitrarily close to the root of $\alpha-p+log(2)-p=0$). Hence, for any $\sigma>0$, there exists $\tilde{\alpha}$ such that if $\alpha>\tilde{\alpha}$, $|\mrevfunc(p^*)-\frac{\alpha+\log(2)}{2}|<\sigma$. Now, we notice that for any $\sigma>0$, there exists $\tilde{\alpha}$ such that if $\alpha>\tilde{\alpha}$, then $|\mrevfunc((1-\epsilon_1)\alpha)-(1-\epsilon_1)\alpha|<\sigma$. Finally, \[\lim_{\alpha\to \infty} \frac{\mrevfunc(p^*)}{\mrevfunc(\hat{p})}\le \lim_{\alpha\to \infty}\frac{\frac{\alpha+\log(2)}{2}}{(1-\epsilon_1)\alpha}=\frac{1}{2(1-\epsilon_1)}\leq \frac{1}{2}+\epsilon,\]
proving the bound is asymptotically tight.\hfill $\square$
\end{proof}
\section{The Effect of Data Censoring on Pricing}\label{sec:censresults}
In this section we study the effect of data censoring on pricing. First, we characterize how censoring distorts the estimated customer demand model from the data, in the case of single product pricing. Next, we study a simple, yet illuminating case of demand and we provide a sufficient condition under which the asymptotic model-free optimal price outperforms the asymptotic optimal price estimated from the censored data. Finally, we extend the insights derived for the problem of single product pricing to that of multiple products pricing, when customers are making decisions based on the symmetric MNL choice model.

\subsection{Single Product}
\begin{lemma}\label{lemma:censlemma}
Suppose Assumption~\ref{asp:cdf-F-G} holds. Then, if at each historical price only a fraction $z$ of non-purchasing customers are observed, the asymptotic estimated purchase probability of a customer at price $p$ would be $\hat{F}(\cdot)=\frac{1-F(p)}{1-F(p)+z F(p)}$. 
\end{lemma}

\begin{proof}{Proof of Lemma \ref{lemma:censlemma}.}
Let $B(p)$ be the number of customers that purchase the product at price $p$ in the data, and let $N(p)$ indicate the number of non-purchasing customers at price $p$. By assumption, only $z N(p)$ of these non-purchasing customers are recorded in the data. By definition, at any price $p$, we have $1-\hat{F}(p)=\frac{B(p)}{B(p)+z N(p)}$, while $1-F(p)=\frac{B(p)}{B(p)+N(p)}$. Therefore, \[B(p)=(1-\hat{F}(p))(B(p)+z N(p))=(1-F(p))(B(p)+N(p)).\] 

Thus we have
\[1-\hat{F}(p)=(1-F(p))\frac{B(p)+N(p)}{B(p)+z N(p)}=(1-F(p))\frac{\frac{B(p)+N(p)}{B(p)+N(p)}}{\frac{B(p)+z N(p)}{B(p)+N(p)}}=\frac{1-F(p)}{1-F(p)+z F(p)}.\]
Where the last equality follows from the fact that \[\frac{B(p)+z N(p)}{B(p)+N(p)}=\frac{B(p)}{B(p)+N(p)}+\frac{z N(p)}{B(p)+N(p)}=1-F(p)+z F(p).\]
\hfill $\square$
\end{proof}

\color{black}

\begin{proposition}\label{prop:unif-bound}
Suppose Assumptions~\ref{asp:cdf-F-G} and~\ref{asp:uniform} hold. Moreover, assume $F(\price) = \frac{p}{a}$.
We have:
\begin{equation*}
    \frac{\revfunc(p^*)}{\revfunc(\hat p)}\ge  \frac{8}{9}.
\end{equation*}
\end{proposition}

\begin{proof} {Proof of Proposition~\ref{prop:unif-bound}.}

From Proposition \ref{prop:compare-opt-price}, we have
\footnotesize
\[p^* \in \argmax_{\price\ge 0} \price\int_\price^{+\infty} g(x)(1-\distfunc(x))\mathrm{d}x=\argmax_{p\geq 0}\frac{p}{ab}\int_{p}^{a}(a-x)\mathrm{d}x= \argmax_{p\geq 0}p(ax-\frac{x^2}{2})\Big|_{p}^{a}=\argmax_{p\geq 0}p(\frac{a^2}{2}-ap+\frac{p^2}{2}).\]
\normalsize
By taking the first-order condition with respect to $p$ we obtain
\[\frac{a^2}{2}-2ap+\frac{3p^2}{2}=0,\] Thus, $p^*=\frac{a}{3},$ which results in a revenue of $\frac{2a}{9}$, while the optimal revenue is $\frac{a}{4}$, proving the claimed ratio. We note that $p=a$ is also a root of $\frac{a^2}{2}-2ap+\frac{3p^2}{2}=0$. However, $a\int_{a}^{b}g(x)(1-F(x))\mathrm{d}x=0$, and hence \[a\neq \argmax_{p\geq 0}p\int_{p}^{b}g(x)(1-F(x))\mathrm{d}x.\]
\hfill $\square$
\end{proof}
\begin{corollary}\label{cor:censor}
Suppose Assumptions~\ref{asp:cdf-F-G} and~\ref{asp:uniform} hold. Moreover, let $F(p)=\frac{p}{a}$. Then as long as at each historical price, only a fraction $z<0.25$ of non-purchasing customers are recorded in the data, the asymptotic optimal model-free price $p^*$ generates a higher revenue than that of the model-based optimal price, estimated from the censored data.
\end{corollary}

\begin{proof}{Proof of Corollary \ref{cor:censor}.}
By Lemma~\ref{lemma:censlemma}, the estimated asymptotic purchase probability at price $p$ would be $\frac{a-p}{a-p+z p}$ and the model-based optimal price will be the optimizer of $p(\frac{a-p}{a-p+z p})$. By taking the first order condition, we note that the optimal model-based price would be $\frac{a(1-\sqrt{z})}{1-z}$ leading to a revenue of $\frac{a(1-\sqrt{z})(\sqrt{z}-z)}{(1-z)^2}$. To conclude the proof, we note that from Proposition~\ref{prop:unif-bound}, the revenue from $p^*$ would be $\frac{2a}{9}$ and it is easy to verify that $\frac{(1-\sqrt{z})(\sqrt{z}-z)}{(1-z)^2}<\frac{2}{9}$ as long as $z< 0.25$.
\hfill $\square$
\end{proof}
\subsection{Multiple Products}
We note that with a slight abuse of notation, we can extend the insight from Corollary~\ref{cor:censor} to the setting with multiple products explored in Section~\ref{exvalidmnl}. We assume the historical prices satisfy $\histprice_{\customer\product}\equiv P_{\customer}$ for all $\customer\in \customerset$. Moreover, we assume that customers make decisions based on the MNL choice model specified in Equation~(\ref{probmnl}). Then, the following corollary specifies the estimated MNL optimal price from censored data in terms of $w_0(\cdot)$, where $w_0(\cdot)$ indicates the positive part of the lambert function and is defined as $w_0(x)e^{w_0(x)}=x, \ x\geq 0$. 
\begin{corollary}\label{cor:mnlcens}
Suppose the firm samples the price vector seen by each past customer $i$, $P_i$, from a distribution with PDF $g(\price) = 1/b$ for $\price\in [0,b]$, where $b>\hat p=\frac{w_0(ne^{\alpha -1})+1}{\beta}$. If at each historical price vector, only a fraction $z$ of non-purchasing customers are recorded in the data, the asymptotic estimated MNL optimal price can be expressed as $\bar{p}=\frac{w_0(\frac{ne^{\alpha -1}}{z})+1}{\beta}$. Moreover, there exists $\hat{z}>0$ such that as long as $0\le z<\hat{z}$ the asymptotic optimal model-free price $p^*$ generates a higher revenue than $\bar{p}$, the model-based optimal price, estimated from the censored data.
\hfill $\square$
\end{corollary}

\begin{proof}{Proof of Corollary \ref{cor:mnlcens}.}
With a slight abuse of notation we can notice that Lemma~\ref{lemma:censlemma} is applicable here. Thus, at each given price $p$, the asymptotic estimated probability of purchase can be expressed as:
\[\hat{F}(p)=\frac{ne^{\alpha-\beta p}}{ne^{\alpha-\beta p}+z}\]
and hence $\bar{p}$ will be the optimizer of:
\[\frac{pne^{\alpha-\beta p}}{ne^{\alpha-\beta p}+z}.\]

By taking the first order condition, it can be easily verified that $\bar{p}=\frac{w_0(\frac{ne^{\alpha -1}}{z})+1}{\beta}$. We note that when there is no censoring in the data, the MNL optimal price is $\hat{p}=\frac{w_0(ne^{\alpha -1})+1}{\beta}$. It follows from Proposition~\ref{prop:mnl} that the model-free optimal price $p^*$ will always guarantee at least $50\%$ of the revenue from $\hat{p}$ and is not affected by censoring as mentioned in Remark~\ref{rem:nopurchasecustomers}. However, we note that since $0\le z\le 1$, we have $\bar{p}\ge \hat{p}\ge p^*$ for all $z$, while $\mrevfunc(p)$ is strictly decreasing in $p$ for all $p> \hat{p}$ (which follows from the fact that $\hat{p}$ is the unique point satisfying the first order condition for the pseudo-concave function $\mrevfunc(p)$ while $\mrevfunc(0)=\mrevfunc(\infty)=0$). Moreover, as the positive part of $w_0(\cdot)$ is known to be a strictly increasing function, $\bar{p}$ is strictly decreasing in $z$. Therefore, $\mrevfunc(\bar{p})$ is strictly increasing in $z$, becoming zero when $z\to 0$, as $\lim_{z\to0}\bar{p}=\infty$, leading to $\lim_{z\to 0} \mrevfunc(\bar{p})=0$. Therefore, there must exist $\hat{z}\in(0,1]$ such that for all $0\le z< \hat{z}$, $\mrevfunc(\bar{p})\le \mrevfunc(p^*)$. 
\hfill $\square$
\end{proof}

\color{black}

\section{Supplementary Results and Discussions}\label{supresults}
\color{mycolor}
\begin{proposition}
\label{prop:average prices}
If prices are set to their average historical values, the \textcolor{mycolor}{worst-case} fraction between the revenue obtained with this approach and the optimal value w.r.t. \eqref{model:OPMIP} is asymptotically zero as the number of historical customers grow.
\end{proposition}
\begin{proof}{Proof of Proposition \ref{prop:average prices}.}
Assume we have $n=1$ product and $m>1$ historical customers where the purchase price of customers $1,\cdots, (m-1)$, was $1$, and the purchase price of customer $m$ was 2. Then, the average historic price of the product is $\frac{m+1}{m}$. At this price, customers $1,\cdots,(m-1)$ will not make a purchase under their worst-case valuations while customer $m$ will purchase the product at price $\frac{m+1}{m}$. However, setting the price of the product at $1$ results in a revenue of $m$, since all customers will purchase the product. Thus, the ratio of the revenue from the average price to the optimal value of \eqref{model:OPMIP} could be less than or equal to $\frac{m+1}{m^2}$. Taking the limit $m \rightarrow +\infty$ with respect to this ratio completes the proof.
\hfill $\square$
\end{proof}

\begin{proposition}
\label{prop:random prices}
If prices are set uniformly at random based on the empirical distribution defined by historical prices, the \textcolor{mycolor}{worst-case} fraction between the revenue obtained with this approach and the optimal value w.r.t. \eqref{model:OPMIP} is asymptotically zero as the number of historical customers grow.
\end{proposition}

\begin{proof}{Proof of Proposition \ref{prop:random prices}.}
Assume we have $n$ products and $n$ historical customers. The historical prices $\vechistprice_{\customer}$ observed by customer $\customer$, are defined by the following vectors:
\begin{align*}
    \vechistprice_{1} 
    &= 
   \left(
        1, 
       2,2,\dots, 
       2
    \right),
    \\
    \vechistprice_{2} 
    &= 
   \left(
        2, 
       1,2,\dots, 
       2
    \right),
    \\
    \vechistprice_{3} 
    &= 
    \left(
        2, 
       2,1,2,\dots,
       2
    \right),
    \\
    \dots\\
    \vechistprice_{n} 
    &= 
    \left(
        2, 
       2,\dots,2, 
       1
    \right).
\end{align*}

For the historical purchase choices, each customer $\customer \in \{1,\dots,n\}$ purchased product $\customer$.

Assume the vector of prices $p$ is equal to one of $\vechistprice_{1},\dots, \vechistprice_{n}$ at random (with equal probability). If $p=\vechistprice_{1}$, then for $\customer=1$, we have $P_{1\choice{1}}\leq p_{\choice{1}}$ and hence the revenue from customer $1$ is $1$. For any customer $\customer>1$, we have $P_{\customer\choice{\customer}}< p_{\choice{\customer}}$, hence the revenue form any customer $\customer>1$ is zero. Similarly we can show that if $p=\vechistprice_{i}$ for any $\customer>1$, the total revenue will be $1$, under the worst-case customer valuations. 

However, for $p=(1,1,\dots,1)$, the revenue will be $n$ as every customer will make a purchase at price $1$. Therefore, the ratio of the expected revenue from the price observed by a random historical customer to the optimal value of \eqref{model:OPMIP} could be less than or equal to $\frac{1}{n}$. Taking the limit $n \rightarrow +\infty$ with respect to this ratio completes the proof.
\hfill $\square$
\end{proof}
In the following example we show that the prices obtained through the LP relaxation of the program \eqref{model:OPMIP}, have a worse worst-case performance than the conservative pricing approach.
\smallskip
\begin{example}\label{exp:lp-fail}
    Consider an instance with $\numcustomers=3$ customers and $\numproducts=2$ products.
    The historical prices $\histprice_{\customer \product}$ and the customer choices (in bold) are listed in the table below:

    \medskip
    \begin{center}
        \begin{tabular}{c|cc}
Price            & Product 1 & Product 2 \\
            \hline
            Customer 1 & \textbf{1} & 2 \\
            Customer 2 & 2 & \textbf{3} \\
            Customer 3 & \textbf{1} & 3 \\
        \end{tabular}
    \end{center}

    \medskip
    It can be shown that the LP relaxation yields the price vector
    $\vecprice^{\textnormal{LP}} = (1.2, 2.3)$ and an upper bound of $4.8$ to the optimal value of \eqref{model:OPMIP}. However, when evaluating \eqref{model:OPMIP}
    with variables $\vecprice$ fixed to $\vecprice^{\textnormal{LP}}$, both customers 1 and 3 do not purchase any products, while customer 2 purchases product 1 in the worst-case. The following set of valuations, $v_{ij}$, that are drawn from the IC polyhedra of customers 1, 2 and 3 conform to these customer choices that lead to a worst-case revenue for the firm.
    
            \medskip
    \begin{center}
        \begin{tabular}{c|cc}
     Valuation       & Product 1 & Product 2 \\
            \hline
            Customer 1 & 1 & 0 \\
            Customer 2 & 2 & 3 \\
            Customer 3 & 1 & 0 \\
        \end{tabular}
    \end{center}
    \medskip
    Thus, the total revenue generated from $\vecprice^{\textnormal{LP}}$ is \$1.2. The optimal solution of this instance is \$4.0 certified by prices $\vecprice^* = (1.0, 2.0)$, which leads to an LP solution ratio of 30\%, worse than the conservative price ratio of $\ubar{\histprice} / \bar{\histprice} = 1/3 \approx 33\%$.

    The optimal LP solution associated with variables $\vecallocVar$ provides insights into the poor performance of the heuristic. In particular, consider inequality \eqref{eq:OPBdualfeasible-cl} for customer 2 ($\customer = 2$) and product 1 ($\product = 1$):
    \begin{align*}
        \dualS_2 \le \price_1 + \histprice_{2 2}(1 - \allocVar_{21}).
    \end{align*}
    At the LP optimality, $\allocVar_{2 1}=0.4$ and the above right-hand side is equal to 3. The inequality is also tight, leading to a (relaxed) revenue of $\dualS_2 = 3.0$. However, with integrality constraints, $\allocVar_{2 1}=1.0$ and the constraint is again tight with $\dualS_2 \le \price_1 = 1.2$, which is a significant decrease in revenue. The same issue is identified for the other customers.

    More generally, the big-M structure of inequalities \eqref{eq:OPBdualfeasible-cl} tends to result in the optimal LP solution pricing $\price_\product$ slightly higher than the historical prices for most of the customer $\customer$ with $\choice{\customer} = \product$, in spite of the fact that the indicator $\indicator(\price_\product > \histprice_{\customer \choice{\product}})$ has been encoded in \eqref{model:OPMIP} to represent the no-purchase option. More precisely, $\allocVar_{\customer \product}$ can be set to $1 - \eps$ for any sufficiently small $\eps$ to overcome that condition, and the objective value of the LP is not impacted significantly since the customer is still assumed to purchase the product. However, when the integrality constraint is imposed, either the inequality becomes binding with respect to some price or the no-purchase option is chosen, leading to a smaller expected revenue. \hfill$\square$
\end{example}

\end{document}